\documentclass[12pt]{book}
\usepackage{here,amsfonts,amssymb,pstricks,pst-node,
pst-plot,pst-coil,epsfig,makeidx}
\usepackage{amsmath}
\newtheorem{thm}{Theorem}[section]
\newtheorem{lemma}[thm]{Lemma}

\newtheorem{prop}[thm]{Proposition}

\newtheorem{example}{Example}[chapter]
\newtheorem{defin}[thm]{Definition}

\newbox\sample
\newif\ifproofmode
\newif\ifsymindex
\global\symindexfalse
\newwrite\inx
\def\indsyma#1#2{\ifproofmode\marginpar{$\scriptstyle#1$}\fi%
\ifx#2\empty\write\inx{$\noexpand#1$,\space\thepage}%
\write\inx{\string\newline}\else%
\write\inx{$\noexpand#1$,\space#2,\space\thepage}%
\write\inx{\string\newline}\fi\ignorespaces}%
\def\indsym#1#2{\ifsymindex%
\ifproofmode\marginpar{$\scriptstyle#1$}\fi%
\ifx#2\empty\write\inx{\string\item \space$\noexpand#1$,\space\thepage}%
\else%
\write\inx{\string\item \space$\noexpand#1$,\space#2,\space\thepage}%
\fi\ignorespaces\fi}%
\newskip\dangerskipb
\newskip\dangerskip
\def\hang{\hangindent\dangerskip}

\def\s#1{{\cal #1}}

\def\lag{\left\langle}
\def\rag{\right\rangle}

\def\proof{\noindent{\it Proof\/}.\enspace}

\def\remark{\bigskip\noindent{\bf Remark:}\enspace}
\def\endremark{\bigskip}
\def\remarks{\bigskip\noindent{\bf Remarks:}\enspace}

\font\manual=manfnt at 12pt
\def\danbend{{\manual\char127}}
\def\datanger{\medbreak\begingroup\clubpenalty=10000
 \def\par{\endgraf\endgroup\medbreak} \noindent\hang\hangafter=-2
 \hbox to0pt{\hskip-3.5pc\danbend\hfill}}
\outer\def\danger{\datanger}%
\def\ddatanger{\medbreak\begingroup\clubpenalty=10000
 \def\par{\endgraf\endgroup\medbreak} \noindent\hang\hangafter=-2
 \hbox to0pt{\hskip-3.5pc\danbend\kern1pt%
\danbend\hfill}}

\def\dobdownarrow{\mathop{\vbox{\kern2pt \hbox{$\Big\downarrow$}\kern-16.5pt
                          \nointerlineskip\hbox{$\Big\downarrow$}}}}

\def\lrightarrow{\hbox to 25pt{\rightarrowfill}}

\def\supexp{exp(m,n,p)=m^{m^{m^{\cdot^{\cdot^{\cdot^{m^{p}}}}}}}
\vbox{\hbox{$\Big\}\scriptstyle n$}\kern0pt}}

\def\supexpo#1#2#3{#1^{#1^{\cdot^{\cdot^{\cdot^{#1^{#2}}}}}}
\vbox{\hbox{$\Big\}\scriptstyle #3$}\kern0pt}}

\def\sqr#1#2{{\vcenter{\hrule height .#2pt
         \hbox{\vrule width.#2pt height#1pt \kern#1pt
             \vrule width.#2pt}
         \hrule height.#2pt}}}

\def\square{\mathchoice\sqr64\sqr64\sqr{2.1}3\sqr{1.5}3}

\def\bigsquare{\mathchoice\sqr76\sqr76\sqr{2.1}3\sqr{1.5}3}

\def\lag{\langle}
\def\rag{\rangle}

\def\co{\colon}

\newskip\bogcentering \bogcentering= 0pt plus 1000pt minus 1000pt 

\def\matth{\mathsurround=0pt}
\def\fakrightarrowfill{$\matth \mathord- \mkern-6mu
  \cleaders\hbox{$\mkern-2mu \mathord- \mkern-2mu$}\hfill
 \mkern-6mu \mathord\rightarrow$}

\def\fakoverrightarrow#1{\vbox{\ialign{##\crcr
  \fakrightarrowfill\crcr\noalign{\kern-1pt\nointerlineskip}
 $\hfil\displaystyle{#1}\hfil$\crcr}}}

\def\cases#1{\left\{\,\vcenter{\normalbaselines\matth
  \ialign{$##\hfil$&\quad##\hfil\crcr#1\crcr}}\right.}

\newif\ifdtatp

\def\displaty{%
\global \dtatptrue \openup \jot \matth \everycr{\noalign{\ifdtatp \global 
\dtatpfalse \vskip -\lineskiplimit \vskip \normallineskiplimit \else 
\penalty \interdisplaylinepenalty \fi }}}

\def\displaylignes#1{\displaty
   \halign{\hbox to\displaywidth{$\displaystyle##$}\crcr
   #1\crcr}}
\def\eqaligneno#1{\displaty \tabskip=\bogcentering
 \halign to\displaywidth{\hfil$\displaystyle{##}$\tabskip=0pt
 &$\displaystyle{{}##}$\hfil\tabskip=\bogcentering
 &\llap{$##$}\tabskip=0pt\crcr
 #1\crcr}}
\def\leqaligneno#1{\displaty \tabskip=\bogcentering
 \halign to\displaywidth{\hfil$\displaystyle{##}$\tabskip=0pt
 &$\displaystyle{{}##}$\hfil\tabskip=\bogcentering
 &\kern-\displaywidth\rlap{$##$}\tabskip=\displaywidthpt\crcr
 #1\crcr}}
\def\ligne{\hbox to\hsize}
\newdimen\nouvpagewidth
\newdimen\offwidth
\newdimen\lawidthoui
\nouvpagewidth=\lawidthoui
\def\kboxit#1{\vbox{\hrule\hbox{\vrule\kern3pt
              \vbox{\kern3pt#1\kern3pt}\kern3pt\vrule}\hrule}}

\def\kboxitb#1{\vbox{\hrule\hbox{\vrule\kern3pt
              \vbox{\kern3pt#1\kern3pt}\kern3pt\vrule}\hrule}}

\def\laboxaround#1{
\aboxaround{\hbox to\hsize{\hfill\box2\hfill}}{#1}
}

\def\boxar#1#2{
\aboxaround{\hbox to\hsize{\hfill#1\hfill}}{#2}
}

\def\aboxaround#1#2{
\setbox4=\vbox{\hsize #2\noindent\strut#1\strut}
\kboxitb{\box4}}

\def\kframeit#1{\vbox{\hrule\hbox{\vrule\kern5pt
              \vbox{\kern5pt#1\kern5pt}\kern5pt\vrule}\hrule}}

\newskip\savnormalbaselineskip
\newskip\savnormallineskip
\newdimen\savnormallineskiplimit

\def\square{\bigsquare}

\def\dInt{{\rm Int}\>}
\def\opball#1#2{B_{0}(#1,#2)}
\def\Ispac#1#2{#1^{(#2)}}
\def\flecheabov#1{\buildrel #1\over\longrightarrow}
\def\shortexact#1#2#3#4#5{0\,\flecheabov{}\,#1\,\flecheabov{#2}\,#3
\,\flecheabov{#4}\,#5\,\flecheabov{}\,0}
\def\dist#1#2#3{d_{#1}(#2,\,#3)}
\def\absval#1#2{|#1 - #2|}
\def\eudist#1#2#3{\left(|#1_{1} - #2_{1}|^{2}+ \cdots + 
|#1_{#3} - #2_{#3}|^{2}\right)^{\frac{1}{2}}}
\def\linvec#1#2{(#1_{1},\ldots, #1_{#2})}
\def\cloball#1#2{B(#1,#2)}
\def\eunorme#1#2{\left(|#1_{1}|^{2} + \cdots + 
|#1_{#2}|^{2}\right)^{\frac{1}{2}}}
\def\ncloball#1{B(#1)}
\def\nopball#1{B_{0}(#1)}
\def\reals{\mathbb{R}}
\def\complex{\mathbb{C}}
\def\integs{\mathbb{Z}}
\def\natnums{\mathbb{N}}
\def\rats{\mathbb{Q}}

\def\affs{\s{E}}
\def\mapdef#1#2#3{#1\co #2\rightarrow #3}

\def\famil#1#2#3{(#1_{#2})_{#2\in #3}}

\def\affreal{\mathbb{A}}
\def\libvecbo#1#2{{\bf #1#2}}
\def\novect#1{#1}
\def\interio#1{\buildrel \circ\over #1}
\def\dimm{\mathrm{dim}}
\def\norme#1{\left\|#1\right\|}
\def\smnorme#1{\|#1\|}
\def\adher#1{\overline{#1}}
\def\fr#1{\mathrm{Fr}\>#1}

\def\dBd{\partial}
\def\dstar{{\rm St}\>}
\def\dcstar{\overline{\hbox{\rm St}}\>}

\title{The Classification Theorem for Compact Surfaces\\
And  A Detour On Fractals\\
}
\author{Jean Gallier\\
Department of Computer and Information Science\\
University of Pennsylvania\\
Philadelphia, PA 19104, USA\\
e-mail: {\tt jean@saul.cis.upenn.edu}\\
\ \\
\copyright\ Jean Gallier\\
{\bf Please, do not} reproduce {\bf without permission} of the author}
\begin{document}
\maketitle
\ \vfill\eject
\begin{center}
{\large \bf
The Classification Theorem for Compact Surfaces\\
And  A Detour On Fractals
}
\end{center}

\vspace{1cm}
\begin{center}
Jean Gallier
\end{center}

\vspace{2cm}

\noindent
{\bf Abstract\/}.
In the words of Milnor himself, the classification theorem
for compact surfaces is a formidable result. According to Massey,
this result was obtained in the early
$1920$'s and was the culmination of the work of many.
Indeed, a rigorous proof requires, among other things, 
a precise definition of a surface and of orientability, 
a precise notion of triangulation, and a precise
way of determining whether two 
surfaces are  homeomorphic or not. This requires some
notions of algebraic topology such as, fundamental groups,
homology groups, and the Euler-Poincar\'e characteristic.
Most steps of the proof are rather involved and 
it is easy to loose track.

\medskip
The purpose of these notes is to present a fairly complete
proof of the classification Theorem for compact surfaces.
Other presentations are often quite informal (see the references
in Chapter V) and we have tried to be more rigorous. Our main source
of inspiration is the beautiful book on
Riemann Surfaces by Ahlfors and Sario. However, Ahlfors and Sario's
presentation is very formal and quite compact.
As a result, uninitiated readers will probably have a hard
time reading this book.

\medskip
Our goal is to help the reader reach the top of the mountain
and help him not to get lost or discouraged
too early. This is not an easy task!

\medskip
We provide quite a bit of topological background material
and the basic facts of algebraic topology needed for 
understanding how the proof goes, with  more than
an impressionistic feeling. 
We hope that these notes will be helpful to 
readers interested in geometry, and who still believe in 
the rewards of serious hiking!

\tableofcontents
\vfill\eject
\chapter{Surfaces}
\label{chap1}
\section{Introduction}
\label{intro}
Few things are as rewarding as finally stumbling upon
the view of a breathtaking landscape at the turn of a path
after a long hike. Similar experiences occur in mathematics,
music, art, etc.
When I first read about the classification of the compact surfaces,
I sensed that if I prepared myself for a long hike, I could
probably enjoy the same kind of exhilarating feeling.

\medskip
In the words of Milnor himself, the classification theorem
for compact surfaces is a formidable result. According to Massey
\cite{Massey87}, this result was obtained in the early
$1920$'s, and was the culmination of the work of many.
Indeed, a rigorous proof requires, among other things, 
a precise definition of a surface and of orientability, 
a precise notion of triangulation, and a precise
way of determining whether two 
surfaces are  homeomorphic or not. This requires some
notions of algebraic topology such as, fundamental groups,
homology groups, and the Euler-Poincar\'e characteristic.
Most steps of the proof are rather involved and 
it is easy to loose track.

\medskip
One aspect of the proof that I find particularly 
fascinating is the use of certain kinds of graphs
(cell complexes) and of some kinds of rewrite rules on these
graphs, to show that every triangulated surface is equivalent
to some cell complex {\it in normal form \/}.
This presents a challenge to researchers interested
in rewriting, as the objects are unusual (neither terms
nor graphs), and rewriting is really modulo cyclic
permutations (in the case of boundaries).
We hope that these notes will  inspire some of the researchers
in the field of rewriting to investigate these mysterious
rewriting systems.

\medskip
Our goal is to help the reader reach the top of the mountain
(the classification theorem for compact surfaces, with or without
boundaries (also called borders)), 
and help him not to get lost or discouraged
too early. This is not an easy task!
On the way, we will take a glimpse at fractals defined
in terms of iterated function systems. 

\medskip
We provide quite a bit of topological background material
and the basic facts of algebraic topology needed for 
understanding how the proof goes, with  more than
an impressionistic feeling. 
Having reviewed some material on complete and compact metric
spaces, we indulge in a short digression on the
Hausdorff distance between compact sets, and the definition
of fractals in terms of iterated function systems.
However, this is just a pleasant interlude, 
our main goal being the classification theorem for compact surfaces.

\medskip
We also review abelian groups, and present a proof of the structure
theorem for finitely generated abelian groups due to Pierre Samuel.
Readers with a good
mathematical background should proceed directly to Section
\ref{surf1}, or even to Section \ref{scomplexes}.

\medskip
We hope that these notes will be helpful to 
readers interested in geometry, and who still believe in 
the rewards of serious hiking!

\medskip\noindent
{\it Acknowledgement\/}:
I would like to thank Alexandre Kirillov for inspiring
me to learn about fractals, through his excellent lectures
on fractal geometry given in the Spring of 1995.
Also many thanks to Chris Croke, Ron Donagi, David Harbater,
Herman Gluck,  and Steve Shatz, from whom
I learned most of my topology and geometry.
Finally, special thanks
to Eugenio Calabi and Marcel Berger,
for giving fascinating courses in the Fall of 1994, which
changed my scientific life irrevocably (for the best!).

\bigskip
Basic topological notions are given in Chapter
\ref{chap7}. In this chapter, we simply review
quotient spaces.

\section{The Quotient Topology}
\label{quotop}
Ultimately, surfaces will be viewed as  spaces obtained
by identifying (or gluing) edges of plane polygons
and to define this process rigorously, we need the concept
of quotient topology. This section
is intended as a review and it is far from being complete.
For more details, consult Munkres  \cite{Munkrestop},
Massey \cite{Massey87,Massey}, Amstrong \cite{Amstrong}, or 
Kinsey \cite{Kinsey}.

\begin{defin}
\label{quotopdef}
{\em
Given any topological space $X$ and any set $Y$, 
for any surjective function $\mapdef{f}{X}{Y}$, 
we define the {\it quotient topology on $Y$ determined by $f$\/}
(also called the {\it identification topology on $Y$ determined by $f$\/}),
by requiring a subset $V$ of $Y$ to be open if $f^{-1}(V)$ is an open set
in $X$. Given an equivalence relation $R$ on a topological space
$X$, if $\mapdef{\pi}{X}{X/R}$ is the projection 
sending every $x\in X$ to its equivalence class $[x]$ in $X/R$,
the space $X/R$ equipped with the quotient topology determined by $\pi$
is called the {\it quotient space of $X$ modulo $R$\/}.
Thus a set $V$ of equivalence classes in $X/R$ is open iff
$\pi^{-1}(V)$ is open in $X$, which is equivalent to the fact
that $\bigcup_{[x]\in V} [x]$ is open in $X$.
}
\end{defin}

\medskip
It is immediately verified that 
Definition \ref{quotopdef} defines topologies, and that
$\mapdef{f}{X}{Y}$ and $\mapdef{\pi}{X}{X/R}$ are continuous
when $Y$ and $X/R$ are given these quotient topologies.

\danger
One should be careful that if $X$ and $Y$ are topological spaces
and $\mapdef{f}{X}{Y}$ is a continuous surjective map,
$Y$ {\it does not\/} necessarily have the quotient topology
determined by $f$. 
Indeed, it may not be true that a subset $V$ of $Y$ is open
when $f^{-1}(V)$ is open. However, this will be true in two important
cases. 

\begin{defin}
\label{opclomap}
{\em
A continuous map $\mapdef{f}{X}{Y}$ is an {\it open map\/}
(or simply {\it open\/}) if $f(U)$ is open in $Y$ whenever $U$ is open in $X$,
and similarly, $\mapdef{f}{X}{Y}$ is a {\it closed map\/}
(or simply {\it closed\/}) if $f(F)$ is closed in 
$Y$ whenever $F$ is closed in $X$.
}
\end{defin}

\medskip
Then, $Y$ has the quotient topology induced by the continuous
surjective map $f$ if either $f$ is open or $f$ is closed.
Indeed, if $f$ is open, then assuming that $f^{-1}(V)$ is open
in $X$, we have $f(f^{-1}(V)) = V$ open in $Y$.
Now, since $f^{-1}(Y - B) = X - f^{-1}(B)$, for any subset $B$ of $Y$,
a subset $V$ of $Y$ is open in the quotient topology iff
$f^{-1}(Y - V)$ is closed in $X$. From this, we can deduce that
if $f$ is a closed map, then $V$ is open in $Y$ iff $f^{-1}(V)$
is open in $X$.

\medskip
Among the desirable features of the quotient topology,
we would like  compactness, connectedness,
arcwise connectedness, or the Hausdorff separation
property, to be preserved. 
Since $\mapdef{f}{X}{Y}$ and $\mapdef{\pi}{X}{X/R}$ are continuous,
by Proposition \ref{concont}, its version for arcwise connectedness,
and Proposition \ref{compac7}, compactness, connectedness, and
arcwise connectedness, are indeed preserved. Unfortunately,
the Hausdorff separation property is not necessarily preserved.
Nevertheless, it is preserved in some special important
cases.

\begin{prop}
\label{quothaus1}
Let $X$ and $Y$ be topological spaces, $\mapdef{f}{X}{Y}$
a continuous surjective map, and assume that $X$ is compact,
and that $Y$ has the quotient topology determined by $f$.
Then $Y$ is Hausdorff iff $f$ is a closed map.
\end{prop}

\proof If $Y$ is Hausdorff, because $X$ is compact
and $f$ is continuous, since every closed set $F$ in $X$ is compact,
by Proposition \ref{compac7},
$f(F)$ is compact, and since $Y$ is Hausdorff, $f(F)$ is closed,
and $f$ is a closed map. For the converse, we use Proposition 
\ref{compac4}. Since $X$ is Hausdorff, every set $\{a\}$ consisting
of a single element $a\in X$ is closed, and since $f$ is a closed map,
$\{f(a)\}$ is also closed in $Y$. Since $f$ is surjective,
every set $\{b\}$ consisting of a single element $b\in Y$ is closed.
If $b_1, b_2\in Y$ and $b_1 \not= b_2$, since 
$\{b_1\}$ and $\{b_2\}$ are closed in $Y$ and $f$ is continuous,
the sets $f^{-1}(b_1)$ and
$f^{-1}(b_2)$ are closed in $X$, and thus compact,
and by Proposition \ref{compac4},
there exists some disjoint open sets $U_1$ and $U_2$ such that
$f^{-1}(b_1)\subseteq U_1$ and $f^{-1}(b_2)\subseteq U_2$.
Since $f$ is closed, the sets $f(X - U_1)$ and $f(X - U_2)$ are closed,
and thus the sets
\[\eqaligneno{
V_1 &= Y - f(X - U_1)\cr
V_2 &= Y - f(X - U_2)\cr
}\]
are open, and it is immediately verified that $V_1\cap V_2 = \emptyset$,
$b_1 \in V_1$, and $b_2 \in V_2$. This proves that $Y$ is Hausdorff.
$\square$

\remark
It is easily shown that another equivalent 
condition for $Y$ being Hausdorff is that
\[\{(x_1,x_2)\in X\times X\ |\ f(x_1) = f(x_2)\}\]
is closed in $X\times X$. 
\endremark

Another useful proposition deals with subspaces and the quotient
topology.

\begin{prop}
\label{quothaus2}
Let $X$ and $Y$ be topological spaces, $\mapdef{f}{X}{Y}$
a continuous surjective map, and assume that
$Y$ has the quotient topology determined by $f$.
If $A$ is a closed subset (resp. open subset) of $X$ and $f$ 
is a closed map (resp. is an open map),
then $B = f(A)$ has the same topology considered
as a subspace of $Y$, or as having the quotient topology
induced by $f$.
\end{prop}

\proof Assume that $A$ is open and that $f$ is an open map.
Assuming that $B = f(A)$ has the subspace topology,
which means that the open sets of $B$ are the sets of the form
$B\cap U$, where $U\subseteq Y$ is an open set of $Y$,
because $f$ is open,
$B$ is open in $Y$, and it is immediate that $\mapdef{f|A}{A}{B}$
is an open map. But then, by a previous observation, $B$ has
the quotient topology induced by $f$. 
The proof when $A$ is closed and  $f$ is a closed map is
similar.
$\square$

\medskip
We now define (abstract) surfaces.

\section{Surfaces: A Formal Definition}
\label{surf1}
Intuitively, what distinguishes a surface from an arbitrary
topological space, is that a surface has the property
that for every point on the surface, there is a small
neighborhood  that looks like a little planar region.
More precisely, a surface is a topological space
that can be covered by open sets that can be mapped
homeomorphically  onto open sets of the plane.
Given such an open set $U$ on the surface $S$, there is an open set
$\Omega$ of the plane $\reals^2$, and a homeomorphism
$\mapdef{\varphi}{U}{\Omega}$. The pair $(U, \varphi)$
is usually called a {\it coordinate system\/}, or
{\it  chart\/}, of $S$, and 
$\mapdef{\varphi^{-1}}{\Omega}{U}$ is called a 
{\it parameterization\/} of $U$. We can think of
the maps  $\mapdef{\varphi}{U}{\Omega}$ as defining small planar
maps of small regions on $S$, similar to geographical maps.
This idea can be extended to higher dimensions,
and leads to the notion of a topological manifold.

\begin{defin}
\label{topman}
{\em
For any $m\geq 1$, a {\it (topological) $m$-manifold\/} is
a second-countable, topological Hausdorff space $M$, together
with  an open cover $(U_i)_{i\in I}$ and
a family $(\varphi_i)_{i\in I}$ of homeomorphisms
$\mapdef{\varphi_i}{U_i}{\Omega_i}$, where each $\Omega_i$
is some open subset of $\reals^m$. Each pair $(U_i, \varphi_i)$
is called a {\it coordinate system\/}, or
{\it chart\/} (or local chart) of $M$, 
each homeomorphism $\mapdef{\varphi_i}{U_i}{\Omega_i}$ is called
a {\it coordinate map\/},
and its inverse  $\mapdef{\varphi^{-1}_{i}}{\Omega_i}{U_i}$ is called a 
{\it parameterization\/} of $U_i$.
For any point $p\in M$, for any coordinate system
$(U, \varphi)$ with $\mapdef{\varphi}{U}{\Omega}$,
if $p \in U$, we say that
$(\Omega, \varphi^{-1})$ is a {\it parameterization of $M$ at $p$\/}.
The family  $(U_i, \varphi_i)_{i\in I}$ is often called an
{\it atlas\/} for $M$. 
A {\it (topological) surface\/} is a connected $2$-manifold.
}
\end{defin}

\remarks
\begin{enumerate}
\item[(1)] 
The terminology is not universally agreed upon. For example,
some authors (including Fulton \cite{Fulton95}) call  the maps
$\mapdef{\varphi^{-1}_{i}}{\Omega_i}{U_i}$ charts!
Always check the direction of the homeomorphisms
involved in the definition of a manifold
(from $M$ to $\reals^m$, or the other way around).
\item[(2)]
Some authors define a surface as a $2$-manifold, i.e., they
do not require a surface to be connected. Following 
Ahlfors and Sario \cite{Ahlfors},
we find it more convenient to assume that surfaces are connected.
\item[(3)] 
According to Definition \ref{topman}, 
$m$-manifolds (or surfaces) do not have any differential
structure. This is usually emphasized  by calling such
objects {\it topological\/} $m$-manifolds (or {\it topological\/} surfaces).
Rather than being pedantic, until specified otherwise,
we will simply use the term $m$-manifold (or surface).
A $1$-manifold is also called a curve.
\end{enumerate}

One may wonder whether it is possible that a topological
manifold $M$ be both an $m$-manifold and an $n$-manifold
for $m \not= n$. For example, could a surface also be a curve?
Fortunately, for connected manifolds,
this is not the case. By a deep theorem
of Brouwer (the invariance of dimension theorem), 
it can be shown that a connected
$m$-manifold is not an $n$-manifold
for $n\not= m$.

\medskip
Some readers many find the definition of a surface quite abstract.
Indeed, the definition does not assume that a surface is a subspace
of any given ambient space, say $\reals^n$, for some $n$.
Perhaps, such surfaces should be called ``abstract surfaces''.
In fact, it can be shown that every surface can be
embedded in $\reals^4$, which is somewhat disturbing, since
$\reals^4$ is hard to visualize! Fortunately, orientable
surfaces can be embedded in $\reals^3$.
However, it is not necessary to use these embeddings to 
understand the topological structure of surfaces.
In fact, when it comes to higher-order manifolds, 
($m$-manifolds for $m\geq 3$), 
and such manifolds do arise naturally
in mechanics, robotics and computer vision, even though it can be shown
that an $m$-manifold can be embedded in $\reals^{2m}$
(a hard theorem due to Whitney),
this usually does not help in understanding its structure.
In the case $m = 1$ (curves), it is not too difficult to prove
that a $1$-manifold is homeomorphic to either a circle or
an open line segment (interval).

\medskip
Since an $m$-manifold $M$ has an open cover of  sets homeomorphic
with open sets of $\reals^m$,  an $m$-manifold is locally
arcwise connected and locally compact. By Theorem \ref{arcwiseth},
the connected components of an $m$-manifold are arcwise connected,
and in particular, a surface is arcwise connected.

\medskip
An open subset $U$ on a surface $S$ is called a {\it Jordan region\/}
if its closure $\overline{U}$ can be mapped homeomorphically
onto a closed disk of $\reals^2$, in such a way that
$U$ is mapped onto the open disk, and thus, that
the boundary of $U$ is mapped homeomorphically onto
the circle, the boundary of the open disk. This means that the boundary
of $U$ is a Jordan curve. Since every point in an open set of 
the plane $\reals^2$
is the center of a closed disk contained in that open set, we note that
every surface has an open cover by Jordan regions.

\medskip
Triangulations are a fundamental tool to obtain a deep understanding
of the topology of surfaces. Roughly speaking, a triangulation
of a surface is a way of cutting up the surface into triangular
regions, such that these triangles are the images of
triangles in the plane, and the edges of these
planar triangles form a graph with certain properties.
To formulate this notion precisely, we need to define
simplices and simplicial complexes.
This can be done in the context of any affine space.

\chapter{Simplices, Complexes, and Triangulations}
\label{chap4}
\section{Simplices and Complexes}
\label{scomplexes}
A simplex is just the convex hull of a finite number of
affinely independent points, but we also need to define faces,
the boundary, and the interior, of a simplex.

\begin{defin}
\label{simplexdef}
{\em 
Let $\affs$ be any normed affine space. Given any $n+1$
affinely independent points $a_0,\ldots,a_n$ in $\affs$, the 
{\it $n$-simplex (or simplex) $\sigma$
defined by $a_0,\ldots,a_n$\/} is the convex hull
of the points  $a_0,\ldots,a_n$, that is, the set of all
convex combinations $\lambda_0 a_0 + \cdots + \lambda_n a_n$, where 
$\lambda_0  + \cdots + \lambda_n = 1$, and $\lambda_i \geq 0$
for all $i$, $0\leq i \leq n$. We call $n$ the {\it dimension\/} of 
the $n$-simplex $\sigma$, and the points $a_0,\ldots,a_n$
are the {\it vertices\/} of $\sigma$.
Given any subset $\{a_{i_{0}},\ldots,a_{i_{k}}\}$ of
$\{a_0,\ldots,a_n\}$ (where $0\leq k \leq n$), the $k$-simplex
generated by  $a_{i_{0}},\ldots,a_{i_{k}}$ is called a {\it face\/}
of $\sigma$. A face $s$ of $\sigma$ is a {\it proper face\/}
if $s\not= \sigma$ (we agree that the empty set is a face of any simplex).
For any vertex $a_i$, the face generated by 
$a_0,\ldots,a_{i-1},a_{i+1},\ldots,a_n$ (i.e.,  omitting $a_i$)
is called the {\it face opposite $a_i$\/}. Every face which
is a $(n-1)$-simplex is called a {\it boundary face\/}.
The union of the boundary faces is the {\it boundary of $\sigma$\/},
denoted as $\dBd \sigma$, and the complement of $\dBd \sigma$ in $\sigma$
is the {\it interior\/} $\dInt \sigma = \sigma - \dBd \sigma$ of $\sigma$.
The interior $\dInt \sigma$ of $\sigma$ 
is sometimes called an {\it open simplex\/}.
}
\end{defin}

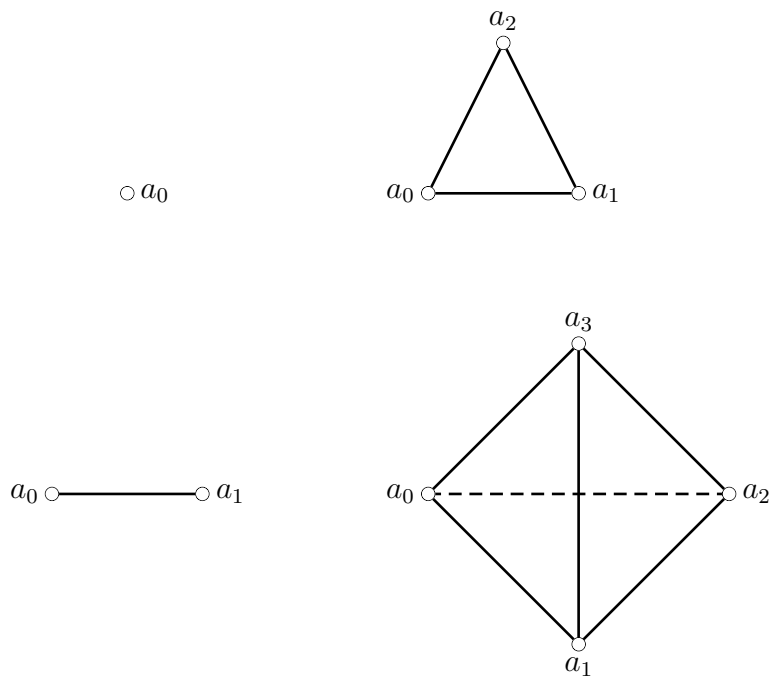
\begin{figure}
  \begin{center}
    \begin{pspicture}(-1,-2)(9,6.2)
    \psline[linewidth=1pt](0,0)(2,0)
    \psline[linewidth=1pt](5,4)(7,4)
    \psline[linewidth=1pt](5,4)(6,6)
    \psline[linewidth=1pt](6,6)(7,4)
    \psline[linewidth=1pt](5,0)(7,-2)
    \psline[linewidth=1pt](5,0)(7,2)
    \psline[linewidth=1pt](9,0)(7,2)
    \psline[linewidth=1pt](9,0)(7,-2)
    \psline[linewidth=1pt](7,-2)(7,2)
    \psline[linewidth=1pt,linestyle=dashed](5,0)(9,0)
    \psdots[dotstyle=o,dotscale=1.5](1,4)
    \psdots[dotstyle=o,dotscale=1.5](0,0)
    \psdots[dotstyle=o,dotscale=1.5](2,0)
    \psdots[dotstyle=o,dotscale=1.5](5,4)
    \psdots[dotstyle=o,dotscale=1.5](7,4)
    \psdots[dotstyle=o,dotscale=1.5](6,6)
    \psdots[dotstyle=o,dotscale=1.5](5,0)
    \psdots[dotstyle=o,dotscale=1.5](7,2)
    \psdots[dotstyle=o,dotscale=1.5](9,0)
    \psdots[dotstyle=o,dotscale=1.5](7,-2)
    \uput[0](1,4){$a_0$}
    \uput[180](0,0){$a_0$}
    \uput[0](2,0){$a_1$}
    \uput[180](5,4){$a_0$}
    \uput[0](7,4){$a_1$}
    \uput[90](6,6){$a_2$}
    \uput[180](5,0){$a_0$}
    \uput[90](7,2){$a_3$}
    \uput[0](9,0){$a_2$}
    \uput[-90](7,-2){$a_1$}
    \end{pspicture}
  \end{center}
  \caption{Examples of simplices}
\label{trangfig1}
\end{figure}

\medskip
It should be noted that for a $0$-simplex consisting of a single
point $\{a_0\}$, $\dBd \{a_0\} = \emptyset$, and
$\dInt \{a_0\} = \{a_0\}$.
Of course, a $0$-simplex is a single point, a $1$-simplex
is the line segment $(a_0, a_1)$, a $2$-simplex is a triangle
$(a_0, a_1, a_2)$ (with its interior),
and a $3$-simplex is a tetrahedron
$(a_0, a_1, a_2, a_3)$ (with its interior), as illustrated
in Figure \ref{trangfig1}.

We now state a number of properties of simplices, 
whose proofs are left as an exercise.
Clearly, a point $x$ belongs to the boundary $\dBd \sigma$ of $\sigma$
iff at least one of its barycentric coordinates
$(\lambda_0,\ldots,\lambda_n)$ is zero, and a point $x$ belongs to the 
interior $\dInt \sigma$ of $\sigma$ iff 
all of  its barycentric coordinates $(\lambda_0,\ldots,\lambda_n)$
are positive, i.e., $\lambda_i > 0$ for all $i, 0\leq i\leq n$.
Then, for every $x\in\sigma$, there is a unique face $s$
such that $x\in \dInt s$, the face generated by those points $a_i$ for which
$\lambda_i > 0$, where $(\lambda_0,\ldots,\lambda_n)$ 
are the barycentric coordinates of $x$.

\medskip
A simplex $\sigma$ is convex,  arcwise connected, compact, and closed.
The interior $\dInt \sigma$ of a simplex is convex, 
arwise connected, open, and
$\sigma$ is the closure of $\dInt \sigma$.

\medskip
For the last property, we recall the following definitions.
The {\it unit $n$-ball\/} $B^n$ is the set of points in $\affreal^n$
such that $x_1^2 + \cdots + x_n^2 \leq 1$.
The {\it unit $n$-sphere\/} $S^{n-1}$ is the set of points in $\affreal^n$
such that $x_1^2 + \cdots + x_n^2 = 1$.
Given a point $a\in \affreal^n$ and a nonnull vector
$\novect{u}\in\reals^n$, the set of all points
$\{a + \lambda\novect{u}\ |\ \lambda\geq 0\}$ is called a
{\it ray emanating from $a$\/}.
Then, every $n$-simplex is homeomorphic to the unit ball 
$B^n$, in such a way that its boundary $\dBd \sigma$ is homeomorphic
to the $n$-sphere $S^{n-1}$.

\medskip
We will prove a slightly more general result about convex sets,
but first, we need a simple proposition.

\begin{prop}
\label{convclo}
Given a normed affine space $\affs$, for any nonempty convex set
$C$, the topological closure $\overline{C}$ of $C$ is also convex.
Furthermore, if $C$ is bounded, then $\overline{C}$ is also
bounded.
\end{prop}

\proof First, we show the following simple inequality.
For any four points $a, b, a', b'\in\affs$, for any $\epsilon > 0$,
for any $\lambda$ such that $0\leq \lambda \leq 1$,
letting  $c = (1 - \lambda) a + \lambda b$ and
$c' = (1 - \lambda) a' + \lambda b'$,
if $\smnorme{\libvecbo{a}{a'}} \leq \epsilon$ and
$\smnorme{\libvecbo{b}{b'}} \leq \epsilon$, then
$\smnorme{\libvecbo{c}{c'}} \leq \epsilon$.

\medskip
This is because 
\[\libvecbo{c}{c'} = (1 - \lambda)\libvecbo{a}{a'} + \lambda \libvecbo{b}{b'},\]
and thus
\[\smnorme{\libvecbo{c}{c'}} \leq 
(1 - \lambda)\smnorme{\libvecbo{a}{a'}} + \lambda \smnorme{\libvecbo{b}{b'}}\leq
(1 - \lambda)\epsilon + \lambda\epsilon = \epsilon.\]
Now, if $a, b\in \overline{C}$, by the definition of closure,
for every $\epsilon > 0$, the open ball $\opball{a}{\epsilon/2}$
must intersect $C$ in some point $a'$,
the open ball $\opball{b}{\epsilon/2}$
must intersect $C$ in some point $b'$,
and by the above inequality,  $c' = (1 - \lambda) a' + \lambda b'$
belongs to the open ball $\opball{c}{\epsilon}$. Since $C$ is convex,
$c' = (1 - \lambda) a' + \lambda b'$ belongs to $C$, and 
$c' = (1 - \lambda) a' + \lambda b'$ also belongs 
to the open ball $B_o(c,\epsilon)$, which shows that
for every $\epsilon > 0$, the open ball $\opball{c}{\epsilon}$ intersects $C$,
which means that $c \in \overline{C}$, and thus that $\overline{C}$ is convex.
Finally, if $C$ is contained in some ball of radius $\delta$,
by the previous discussion, it is clear that $\overline{C}$
is contained in a ball of radius $\delta + \epsilon$,
for any $\epsilon > 0$.
$\square$

\medskip
The following proposition shows that topologically, closed bounded
convex sets in $\affreal^n$ are equivalent to closed balls.
We will need this proposition in dealing with triangulations.

\begin{prop}
\label{convlem1}
If $C$ is any nonempty bounded and convex open set 
in $\affreal^n$, for any point $a\in C$,
any ray emanating from $a$ intersects $\dBd C = \overline{C} - C$
in exactly one point. Furthermore, there is a homeomorphism
of $\overline{C}$ onto the (closed) unit ball $B^n$,
which maps $\dBd C$ onto the $n$-sphere $S^{n-1}$.
\end{prop}

\proof
Since $C$ is convex and bounded,
by Proposition \ref{convclo}, $\overline{C}$ is also convex and bounded.
Given any ray $R=\{a + \lambda\novect{u}\ |\ \lambda\geq 0\}$,
since $R$ is obviously convex, 
the set $R\cap \overline{C}$ is convex, bounded,  and closed in $R$,
which means that $R\cap \overline{C}$ is a closed segment
\[R\cap \overline{C} =
\{a + \lambda\novect{u}\ |\ 0\leq \lambda \leq \mu\},\]for
some $\mu > 0$. Clearly, $a + \mu\novect{u}\in  \dBd C$.
If the ray $R$ intersects
$\dBd C$ in another point $c$, we have
$c = a + \nu\novect{u}$ for some $\nu > \mu$,
and since $\overline{C}$ is convex, 
$\{a + \lambda\novect{u}\ |\ 0\leq \lambda \leq \nu\}$ is
contained in $R\cap \overline{C}$ for $\nu > \mu$,
which is absurd. Thus, every ray emanating from $a$
intersects $\dBd C$ in a single point.

\medskip
Then, the map $\mapdef{f}{\affreal^n - \{a\}}{S^{n-1}}$
defined such that $f(x) = \libvecbo{a}{x}/\smnorme{\libvecbo{a}{x}}$
is continuous. By the first part, the restriction
$\mapdef{f_b}{\dBd C}{S^{n-1}}$ of $f$ to $\dBd C$ is a bijection
(since every point on $S^{n-1}$ corresponds to a unique ray
emanating from $a$). Since $\dBd C$ is a closed and bounded
subset of $\affreal^n$, it is compact, and thus
$f_b$ is a homeomorphism. Consider the inverse
$\mapdef{g}{S^{n-1}}{\dBd C}$ of $f_b$, which is also a
homeomorphism. We need to extend $g$ to a homeomorphism between
$B^n$ and $\overline{C}$. Since $B^n$
is compact, it is enough to extend $g$ to a continuous
bijection. This is done by defining 
$\mapdef{h}{B^n}{\overline{C}}$, such that:
\[h(\novect{u}) = \cases{
    (1 - \smnorme{\novect{u}}) a 
     + \smnorme{\novect{u}} g(\novect{u}/\smnorme{\novect{u}}) & if
     $\novect{u}\not= \novect{0}$;\cr
     a & if $\novect{u} = \novect{0}$.\cr
}\]

\medskip
It is clear that $h$ is bijective and continuous for 
$\novect{u}\not= \novect{0}$.
Since $S^{n-1}$ is compact and $g$ is continuous on $S^{n-1}$,
there is some $M > 0$ such that $\smnorme{\libvecbo{a}{g(\novect{u})}} \leq M$
for all $\novect{u}\in S^{n-1}$, and if $\smnorme{\novect{u}} \leq \delta$, 
then $\smnorme{\libvecbo{a}{h(\novect{u})}} \leq \delta M$, which shows that
$h$ is also continuous for $\novect{u} = \novect{0}$.
$\square$

\remark
It is useful to note that the second part of the proposition
proves that if $C$ is a bounded convex open subset of $\affreal^n$,
then any homeomorphism $\mapdef{g}{S^{n-1}}{\dBd C}$ can be
extended to a homeomorphism 
$\mapdef{h}{B^n}{\overline{C}}$. By Proposition \ref{convlem1},
we obtain the fact that if $C$ is a bounded convex open subset of 
$\affreal^n$, then any homeomorphism $\mapdef{g}{\dBd C}{\dBd C}$ can be
extended to a homeomorphism 
$\mapdef{h}{\overline{C}}{\overline{C}}$.
We will need this fact later on (dealing with
triangulations).
\endremark

We now need to put simplices together to form more complex shapes.
Following Ahlfors and Sario \cite{Ahlfors}, 
we define abstract complexes and
their geometric realizations. This seems easier than defining
simplicial complexes directly, as for example, in Munkres
\cite{Munkresalg}.

\begin{defin}
\label{complexdef}
{\em 
An {\it abstract complex\/} (for short {\it complex\/})
is a pair $K = (V, \s{S})$
consisting of a (finite or infinite) nonempty set $V$
of {\it vertices\/},
together with a family $\s{S}$ of finite subsets of $V$
called {\it abstract simplices\/} (for short {\it simplices\/}),
and satisfying the following conditions:
\begin{enumerate}
\item[(A1)] Every $x \in V$ belongs to at least one and at most a finite
number of simplices in $\s{S}$.
\item[(A2)] Every subset of a simplex $\sigma\in \s{S}$ is also a simplex in $\s{S}$.
\end{enumerate}
If $\sigma\in\s{S}$ is a nonempty simplex of $n+1$ vertices, then
its dimension is $n$, and it is called an {\it $n$-simplex\/}.
A $0$-simplex $\{x\}$ is identified with the vertex $x\in V$.
The {\it dimension of an abstract complex\/} is the maximum dimension of its
simplices if finite, and $\infty$ otherwise.
}
\end{defin}

We will use the abbreviation complex for abstract complex,
and simplex for abstract simplex.
Also, given a simplex $s\in\s{S}$, we will often use the abuse of notation
$s\in K$.
The purpose of condition (A1) is to insure that
the geometric realization of a complex is locally compact. 
Recall that given any set $I$, the real vector space
$\Ispac{\reals}{I}$ freely generated by $I$ is defined 
as the subset of the cartesian product $\reals^I$
consisting of families $\famil{\lambda}{i}{I}$ of elements of $\reals$
with finite support
(where $\reals^I$ denotes the set of all functions
from $I$ to $\reals$).
Then, every abstract complex $(V, \s{S})$ has a geometric 
realization as a topological subspace of the normed vector space
$\Ispac{\reals}{V}$.
Obviously, $\Ispac{\reals}{I}$ can be viewed as a normed affine space
(under the norm $\smnorme{x} = \max_{i\in I} \{x_i\}$)
denoted as $\Ispac{\affreal}{I}$.

\begin{defin}
\label{geomreal}
{\em 
Given an abstract complex $K = (V, \s{S})$, 
its {\it geometric realization\/}
(also called the {\it polytope of $K = (V, \s{S})$\/})
is the subspace $K_g$ 
of $\Ispac{\affreal}{V}$ defined as follows:
$K_g$ is the set of all families $\lambda = (\lambda_a)_{a\in V}$
with finite support, such that:
\begin{enumerate}
\item[(B1)]
$\lambda_a \geq 0$, for all $a\in V$;
\item[(B2)]
The set $\{a\in V\ |\ \lambda_a > 0\}$ is a simplex in $\s{S}$;
\item[(B3)] 
$\sum_{a\in V} \lambda_a = 1$.
\end{enumerate}
}
\end{defin}

For every simplex $s\in \s{S}$, we obtain a subset $s_g$ of $K_g$
by considering those  families $\lambda = (\lambda_a)_{a\in V}$
in $K_g$ such that $\lambda_a = 0$ for all $a\notin s$. Then, by (B2),
we note that 
\[K_g = \bigcup_{s\in \s{S}} s_g.\]
It is also clear that for every $n$-simplex $s$, its geometric
realization $s_g$ can be identified with
an $n$-simplex in $\affreal^n$.

\medskip
Given a vertex $a\in V$, we define the {\it star of $a$\/},
denoted as $\dstar a$, as the finite union of the
interiors $\interio{s_g}$ of the geometric simplices $s_g$ such that $a\in s$.
Clearly, $a\in \dstar a$.
The {\it closed star of $a$\/},
denoted as $\dcstar a$, is the finite union of the
geometric simplices $s_g$ such that $a\in s$.

\medskip
We define a topology on $K_g$ by defining a subset $F$ of $K_g$ to be closed
if $F\cap s_g$ is closed in $s_g$ for all $s\in \s{S}$.
It is immediately verified that the axioms of a topological space are indeed
verified. Actually, we can find a nice basis for this topology,
as shown in the next proposition.

\begin{prop}
\label{startop}
The family of subsets $U$ of $K_g$ such that $U\cap s_g = \emptyset$ for
all  by finitely many $s\in \s{S}$, and such that
$U\cap s_g$ is open in $s_g$ when $U\cap s_g \not= \emptyset$, 
forms a basis of  open sets for the topology of $K_g$. 
For any $a\in V$, the star $\dstar a$ of $a$ is open, the closed star
$\dcstar a$ is the closure of $\dstar a$ and is compact,
and both $\dstar a$ and $\dcstar a$ are arcwise connected.
The space $K_g$ is locally compact, locally arcwise connected, 
and Hausdorff. 
\end{prop}

\proof 
To see that a set $U$ as defined above is open, 
consider the complement $F = K_g - U$ of $U$.
We need to show that $F\cap s_g$ is closed in $s_g$ for all $s\in \s{S}$.
But $F\cap s_g = (K_g - U)\cap s_g = s_g - U$, and if
$s_g \cap U\not= \emptyset$, then $U\cap s_g$ is open in $s_g$, and thus
$s_g - U$ is closed in $s_g$. Next, given any open subset $V$ of $K_g$,
since by $(A1)$, every
$a\in V$ belongs to finitely many simplices $s\in \s{S}$, letting
$U_a$ be the union of the interiors of the finitely many $s_g$ such that
$a\in s$, it is clear that $U_a$ is open in $K_g$, and
that $V$ is the union of the open sets of the form
$U_a\cap V$, which shows that
the sets $U$ of the proposition form a basis of the topology of $K_g$.
For every $a\in V$, the star $\dstar a$ of $a$ has a nonempty
intersection with only finitely many simplices $s_g$,
and $\dstar a\cap s_g$ is the interior of $s_g$ (in $s_g$),
which is open in $s_g$, and $\dstar a$ is open.
That $\dcstar a$ is the closure of $\dstar a$ is obvious,
and since each simplex $s_g$ is compact, and 
$\dcstar a$ is a finite union of compact simplices, it is compact.
Thus, $K_g$ is locally compact.
Since $s_g$ is arcwise connected, for every open set $U$
in the basis, if $U\cap s_g\not= \emptyset$,
$U\cap s_g$ is an open set in $s_g$ that contains
some arcwise connected set $V_s$ containing $a$, and 
the union of these arcwise connected
sets $V_s$ is arcwise connected, and clearly an open set of $K_g$.
Thus, $K_g$ is locally arcwise connected. 
It is also immediate that $\dstar a$ and $\dcstar a$ are arcwise connected.
Let $a, b\in K_g$, and assume that $a\not= b$. If $a, b\in s_g$ for some
$s\in\s{S}$, since $s_g$ is Hausdorff, there are disjoint open sets
$U, V\subseteq s_g$ such that $a\in U$ and $b\in V$.
If $a$ and $b$ do not belong to the same simplex, then
$\dstar a$ and $\dstar b$ are disjoint open sets such that
$a\in\dstar a$ and $b\in\dstar b$. 
$\square$

\medskip
We also note that for any two simplices $s_1, s_2$ of $\s{S}$,
we have 
\[(s_1\cap s_2)_g = (s_1)_g \cap (s_2)_g.\]

\medskip
We say that a complex $K = (V, \s{S})$ is connected if
it is not the union of two  complexes
$(V_1, \s{S}_1)$ and $(V_2, \s{S}_2)$, 
where $V = V_1\cup V_2$
with $V_1$ and $V_2$ disjoint,
and $\s{S} = \s{S}_1\cup\s{S}_2$ with $\s{S}_1$ and $\s{S}_2$ disjoint. 
The next proposition shows that a connected complex 
contains countably many simplices. This is an important fact, since
it implies that if a surface can be triangulated, then its topology must
be second-countable.

\begin{prop}
\label{concomplex}
If $K = (V, \s{S})$ is a connected complex, then $\s{S}$ and $V$
are countable.
\end{prop}

\proof
The proof is very similar to that of the second
part of Theorem \ref{arcwiseth}. 
The trick consists in defining the right notion of arcwise
connectedness. We say that two vertices $a, b\in V$ are 
{\it path-connected, or that there is a path from $a$ to $b$\/}
if there is a sequence $(x_0,\ldots,x_n)$ of vertices $x_i\in V$,
such that $x_0 = a$, $x_n = b$, and $\{x_i, x_{i+1}\}$,
is a simplex in $\s{S}$, for all $i, 0\leq i\leq n-1$.
Observe that every simplex $s\in \s{S}$ is path-connected.
Then, the proof consists in showing that if $(V, \s{S})$ is
a connected complex, then it is path-connected.
Fix any vertex $a\in V$, and let $V_a$ be the set of all
vertices that are path-connected to $a$. We claim that for any
simplex $s\in\s{S}$, if $s\cap V_a \not= \emptyset$, then
$s\subseteq V_a$, which shows that if  $\s{S}_a$ is the subset of $\s{S}$
consisting of all simplices having some vertex in $V_a$,
then $(V_a, \s{S}_a)$ is a complex.
Indeed, if $b\in s\cap V_a$, there is a path from $a$ to $b$.
For any $c\in s$, since $b$ and $c$ are path-connected, then
there is a path from $a$ to $c$, and $c\in V_a$, which
shows that $s\subseteq V_a$. A similar reasoning applies to the 
complement $V - V_a$ of $V_a$, and we obtain a complex
$(V - V_a, \s{S} - \s{S}_a)$. But $(V_a, \s{S}_a)$ and
$(V - V_a, \s{S} - \s{S}_a)$ are disjoint complexes, contradicting 
the fact that $(V, \s{S})$ is connected.
Then, since every simplex $s\in \s{S}$ is finite and
every path is finite, the number of path from $a$ is countable,
and because $(V, \s{S})$ is path-connected,
there are at most countably many vertices in $V$ and 
at most countably many simplices $s\in\s{S}$.
$\square$

\section{Triangulations}
\label{triangsec}
We now return to surfaces and define the notion of triangulation. 
Triangulations are special kinds of complexes of dimension $2$,
which means that the simplices involved are points, line segments,
and triangles.

\begin{defin}
\label{triangul}
{\em 
Given a surface $M$, a {\it triangulation of $M$\/} 
is a pair $(K, \sigma)$ consisting of 
a $2$-dimensional complex $K = (V,\s{S})$
and of a map $\mapdef{\sigma}{\s{S}}{2^{M}}$ 
assigning a closed subset
$\sigma(s)$ of $M$ to every simplex $s\in \s{S}$, 
satisfying the following conditions:
\begin{enumerate}
\item[(C1)] $\sigma(s_1\cap s_2) = \sigma(s_1)\cap\sigma(s_2)$, for all 
$s_1, s_2\in\s{S}$.
\item[(C2)] For every $s\in\s{S}$, there is a homeomorphism $\varphi_s$ from
the geometric realization $s_g$ of $s$ to $\sigma(s)$, 
such that $\varphi_s(s'_g) = \sigma(s')$, for every
$s' \subseteq s$.
\item[(C3)] $\bigcup_{s\in \s{S}} \sigma(s) = M$, that is, the sets
$\sigma(s)$ cover $M$.
\item[(C4)] For every point $x\in M$, there is some neighborhood of $x$ which
meets only finitely many of the $\sigma(s)$.
\end{enumerate}
}
\end{defin}

If $(K, \sigma)$ is a triangulation of $M$, we also
refer to the map $\mapdef{\sigma}{\s{S}}{2^{M}}$ as a 
triangulation of $M$ and we also say that
$K$ is a triangulation $\mapdef{\sigma}{\s{S}}{2^{M}}$ of $M$.
As expected,  given a triangulation $(K, \sigma)$ of a surface $M$,
the geometric realization $K_g$ of $K$ is homeomorphic to the 
surface $M$, as shown by the following proposition.

\begin{prop}
\label{trianglem1}
Given any triangulation $\mapdef{\sigma}{\s{S}}{2^{M}}$
of a surface $M$, there is a homeomorphism
$\mapdef{h}{K_g}{M}$ from the geometric realization $K_g$ of
the complex $K = (V, \s{S})$ onto the surface $M$, such that each 
geometric simplex $s_g$ is mapped onto $\sigma(s)$.
\end{prop}

\proof Obviously, for every vertex $x\in V$,
we let $h(x_g) = \sigma(x)$. If $s$ is a $1$-simplex,
we define $h$ on $s_g$ using any of the homeomorphisms
whose existence is asserted by (C1).
Having defined $h$ on the boundary of each $2$-simplex $s$,
we need to extend $h$ to the entire $2$-simplex $s$.
However, by (C2), there is some homeomorphism $\varphi$ from
$s_g$ to $\sigma(s)$, and if it does not agree with $h$
on the boundary of $s_g$, which is a triangle, by
the remark after Proposition  \ref{convlem1}, since
the restriction of $\varphi^{-1}\circ h$ to the boundary of $s_g$
is a homeomorphism, it can be extended to a homeomorphism
$\psi$ of $s_g$ into itself, and then
$\varphi\circ \psi$ is a homeomorphism of $s_g$ onto $\sigma(s)$
that agrees with $h$ on the boundary of $s_g$.
This way, $h$ is now defined on the entire $K_g$.
Given any closed set $F$ in $M$, for every simplex $s\in \s{S}$,
\[h^{-1}(F)\cap s_g = h^{-1}|_{s_{g}}(F),\]
where $h^{-1}|_{s_{g}}(F)$ is the restriction of $h$ to $s_g$,
which is continuous by construction, and thus,
$h^{-1}(F)\cap s_g$ is closed for all $s\in \s{S}$, which
shows that $h$ is continuous.
The map $h$ is injective
because of (C1), surjective because of (C3),
and its inverse is continuous because of (C4). Thus,
$h$ is indeed a homeomorphism mapping  $s_g$ onto $\sigma(s)$.
$\square$

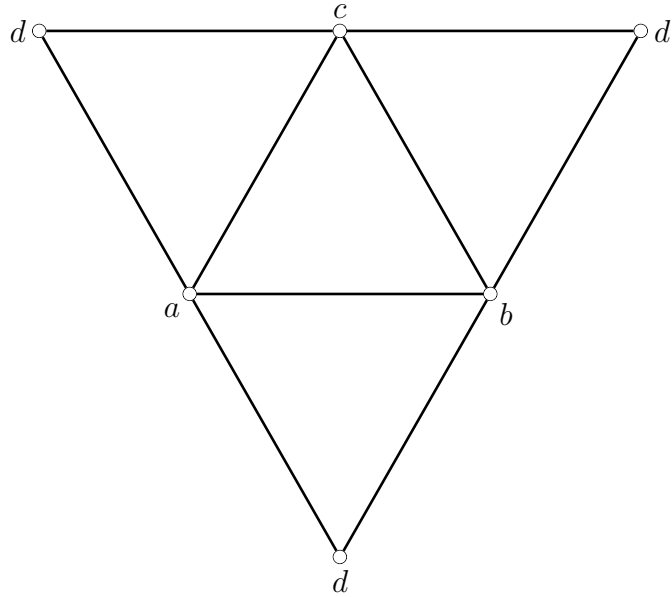
\begin{figure}
  \begin{center}
    \begin{pspicture}(-4,-0.5)(4,7.2)
    \psline[linewidth=1pt](0,0)(-4,7)
    \psline[linewidth=1pt](0,0)(4,7)
    \psline[linewidth=1pt](-4,7)(4,7)
    \psline[linewidth=1pt](-2,3.5)(2,3.5)
    \psline[linewidth=1pt](-2,3.5)(0,7)
    \psline[linewidth=1pt](2,3.5)(0,7)
    \psdots[dotstyle=o,dotscale=1.5](0,0)
    \psdots[dotstyle=o,dotscale=1.5](-4,7)
    \psdots[dotstyle=o,dotscale=1.5](4,7)
    \psdots[dotstyle=o,dotscale=1.5](-2,3.5)
    \psdots[dotstyle=o,dotscale=1.5](2,3.5)
    \psdots[dotstyle=o,dotscale=1.5](0,7)
    \uput[-90](0,0){$d$}
    \uput[180](-4,7){$d$}
    \uput[0](4,7){$d$}
    \uput[-135](-2,3.5){$a$}
    \uput[-45](2,3.5){$b$}
    \uput[90](0,7){$c$}
    \end{pspicture}
  \end{center}
  \caption{A triangulation of the sphere}
\label{spherefig}
\end{figure}

\medskip
Figure \ref{spherefig} shows a triangulation of the {\it sphere\/}.

The geometric realization of the above triangulation is
obtained by pasting together the pairs of 
edges labeled $(a,d)$, $(b,d)$, $(c,d)$. The geometric
realization is a tetrahedron.

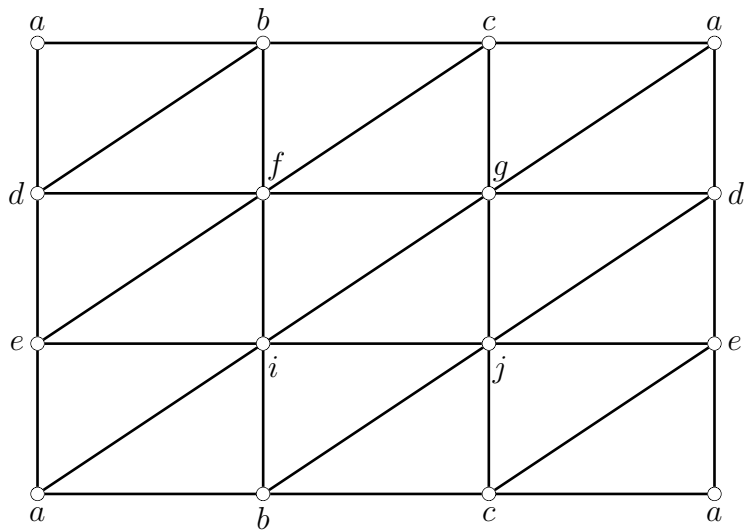
\begin{figure}
  \begin{center}
    \begin{pspicture}(0,-0.5)(9,6.2)
    \psline[linewidth=1pt](0,4)(3,6)
    \psline[linewidth=1pt](6,0)(9,2)
    \psline[linewidth=1pt](0,0)(9,0)
    \psline[linewidth=1pt](0,0)(0,6)
    \psline[linewidth=1pt](9,0)(9,6)
    \psline[linewidth=1pt](0,6)(9,6)
    \psline[linewidth=1pt](0,0)(9,6)
    \psline[linewidth=1pt](0,2)(6,6)
    \psline[linewidth=1pt](3,0)(9,4)
    \psline[linewidth=1pt](0,2)(9,2)
    \psline[linewidth=1pt](0,4)(9,4)
    \psline[linewidth=1pt](3,0)(3,6)
    \psline[linewidth=1pt](6,0)(6,6)
    \psdots[dotstyle=o,dotscale=1.5](0,0)
    \psdots[dotstyle=o,dotscale=1.5](0,2)
    \psdots[dotstyle=o,dotscale=1.5](0,4)
    \psdots[dotstyle=o,dotscale=1.5](0,6)
    \psdots[dotstyle=o,dotscale=1.5](3,0)
    \psdots[dotstyle=o,dotscale=1.5](3,2)
    \psdots[dotstyle=o,dotscale=1.5](3,4)
    \psdots[dotstyle=o,dotscale=1.5](3,6)
    \psdots[dotstyle=o,dotscale=1.5](6,0)
    \psdots[dotstyle=o,dotscale=1.5](6,2)
    \psdots[dotstyle=o,dotscale=1.5](6,4)
    \psdots[dotstyle=o,dotscale=1.5](6,6)
    \psdots[dotstyle=o,dotscale=1.5](9,0)
    \psdots[dotstyle=o,dotscale=1.5](9,2)
    \psdots[dotstyle=o,dotscale=1.5](9,4)
    \psdots[dotstyle=o,dotscale=1.5](9,6)
    \uput[-90](0,0){$a$}
    \uput[180](0,2){$e$}
    \uput[180](0,4){$d$}
    \uput[90](0,6){$a$}
    \uput[-90](3,0){$b$}
    \uput[-60](3,2){$i$}
    \uput[60](3,4){$f$}
    \uput[90](3,6){$b$}
    \uput[-90](6,0){$c$}
    \uput[-60](6,2){$j$}
    \uput[60](6,4){$g$}
    \uput[90](6,6){$c$}
    \uput[-90](9,0){$a$}
    \uput[0](9,2){$e$}
    \uput[0](9,4){$d$}
    \uput[90](9,6){$a$}
    \end{pspicture}
  \end{center}
  \caption{A triangulation of the torus}
\label{torusfig}
\end{figure}

\medskip
Figure \ref{torusfig} shows a triangulation of a surface
called a {\it torus\/}.

The geometric realization of the above triangulation is
obtained by pasting together the pairs of 
edges labeled $(a,d)$, $(d,e)$, $(e,a)$, and the pairs of edges
labeled  $(a,b)$, $(b,c)$, $(c,a)$.

\begin{figure}[H]
  \begin{center}
    \begin{pspicture}(0,-0.5)(9,6.2)
    \psline[linewidth=1pt](0,4)(3,6)
    \psline[linewidth=1pt](6,0)(9,2)
    \psline[linewidth=1pt](0,0)(9,0)
    \psline[linewidth=1pt](0,0)(0,6)
    \psline[linewidth=1pt](9,0)(9,6)
    \psline[linewidth=1pt](0,6)(9,6)
    \psline[linewidth=1pt](0,0)(9,6)
    \psline[linewidth=1pt](0,2)(6,6)
    \psline[linewidth=1pt](3,0)(9,4)
    \psline[linewidth=1pt](0,2)(9,2)
    \psline[linewidth=1pt](0,4)(9,4)
    \psline[linewidth=1pt](3,0)(3,6)
    \psline[linewidth=1pt](6,0)(6,6)
    \psdots[dotstyle=o,dotscale=1.5](0,0)
    \psdots[dotstyle=o,dotscale=1.5](0,2)
    \psdots[dotstyle=o,dotscale=1.5](0,4)
    \psdots[dotstyle=o,dotscale=1.5](0,6)
    \psdots[dotstyle=o,dotscale=1.5](3,0)
    \psdots[dotstyle=o,dotscale=1.5](3,2)
    \psdots[dotstyle=o,dotscale=1.5](3,4)
    \psdots[dotstyle=o,dotscale=1.5](3,6)
    \psdots[dotstyle=o,dotscale=1.5](6,0)
    \psdots[dotstyle=o,dotscale=1.5](6,2)
    \psdots[dotstyle=o,dotscale=1.5](6,4)
    \psdots[dotstyle=o,dotscale=1.5](6,6)
    \psdots[dotstyle=o,dotscale=1.5](9,0)
    \psdots[dotstyle=o,dotscale=1.5](9,2)
    \psdots[dotstyle=o,dotscale=1.5](9,4)
    \psdots[dotstyle=o,dotscale=1.5](9,6)
    \uput[-90](0,0){$d$}
    \uput[180](0,2){$e$}
    \uput[180](0,4){$f$}
    \uput[90](0,6){$a$}
    \uput[-90](3,0){$c$}
    \uput[-60](3,2){$j$}
    \uput[60](3,4){$g$}
    \uput[90](3,6){$b$}
    \uput[-90](6,0){$b$}
    \uput[-60](6,2){$k$}
    \uput[60](6,4){$h$}
    \uput[90](6,6){$c$}
    \uput[-90](9,0){$a$}
    \uput[0](9,2){$f$}
    \uput[0](9,4){$e$}
    \uput[90](9,6){$d$}
    \end{pspicture}
  \end{center}
  \caption{A triangulation of the projective plane}
\label{projplan}
\end{figure}
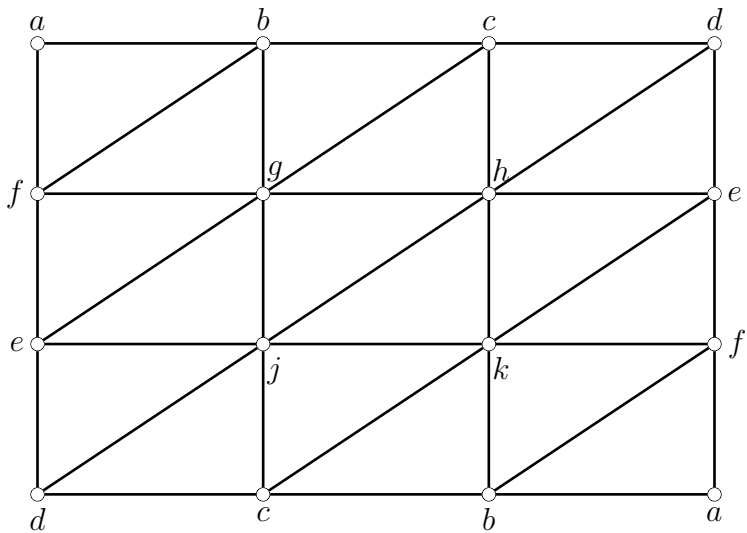

\medskip
Figure \ref{projplan} shows a triangulation of  a surface
called the {\it projective plane\/}.

The geometric realization of the above triangulation is
obtained by pasting together the pairs of 
edges labeled $(a,f)$, $(f,e)$, $(e,d)$, and the pairs of edges
labeled  $(a,b)$, $(b,c)$, $(c,d)$.
This time, the gluing requires a ``twist'', since
the the paired edges have opposite orientation.
Visualizing this surface in $\affreal^3$ is actually
nontrivial.

\begin{figure}
  \begin{center}
    \begin{pspicture}(0,-0.5)(9,6.2)
    \psline[linewidth=1pt](0,4)(3,6)
    \psline[linewidth=1pt](6,0)(9,2)
    \psline[linewidth=1pt](0,0)(9,0)
    \psline[linewidth=1pt](0,0)(0,6)
    \psline[linewidth=1pt](9,0)(9,6)
    \psline[linewidth=1pt](0,6)(9,6)
    \psline[linewidth=1pt](0,0)(9,6)
    \psline[linewidth=1pt](0,2)(6,6)
    \psline[linewidth=1pt](3,0)(9,4)
    \psline[linewidth=1pt](0,2)(9,2)
    \psline[linewidth=1pt](0,4)(9,4)
    \psline[linewidth=1pt](3,0)(3,6)
    \psline[linewidth=1pt](6,0)(6,6)
    \psdots[dotstyle=o,dotscale=1.5](0,0)
    \psdots[dotstyle=o,dotscale=1.5](0,2)
    \psdots[dotstyle=o,dotscale=1.5](0,4)
    \psdots[dotstyle=o,dotscale=1.5](0,6)
    \psdots[dotstyle=o,dotscale=1.5](3,0)
    \psdots[dotstyle=o,dotscale=1.5](3,2)
    \psdots[dotstyle=o,dotscale=1.5](3,4)
    \psdots[dotstyle=o,dotscale=1.5](3,6)
    \psdots[dotstyle=o,dotscale=1.5](6,0)
    \psdots[dotstyle=o,dotscale=1.5](6,2)
    \psdots[dotstyle=o,dotscale=1.5](6,4)
    \psdots[dotstyle=o,dotscale=1.5](6,6)
    \psdots[dotstyle=o,dotscale=1.5](9,0)
    \psdots[dotstyle=o,dotscale=1.5](9,2)
    \psdots[dotstyle=o,dotscale=1.5](9,4)
    \psdots[dotstyle=o,dotscale=1.5](9,6)
    \uput[-90](0,0){$a$}
    \uput[180](0,2){$e$}
    \uput[180](0,4){$d$}
    \uput[90](0,6){$a$}
    \uput[-90](3,0){$b$}
    \uput[-60](3,2){$i$}
    \uput[60](3,4){$f$}
    \uput[90](3,6){$b$}
    \uput[-90](6,0){$c$}
    \uput[-60](6,2){$j$}
    \uput[60](6,4){$g$}
    \uput[90](6,6){$c$}
    \uput[-90](9,0){$a$}
    \uput[0](9,2){$d$}
    \uput[0](9,4){$e$}
    \uput[90](9,6){$a$}
    \end{pspicture}
  \end{center}
  \caption{A triangulation of the Klein bottle}
\label{Kleinbot}
\end{figure}
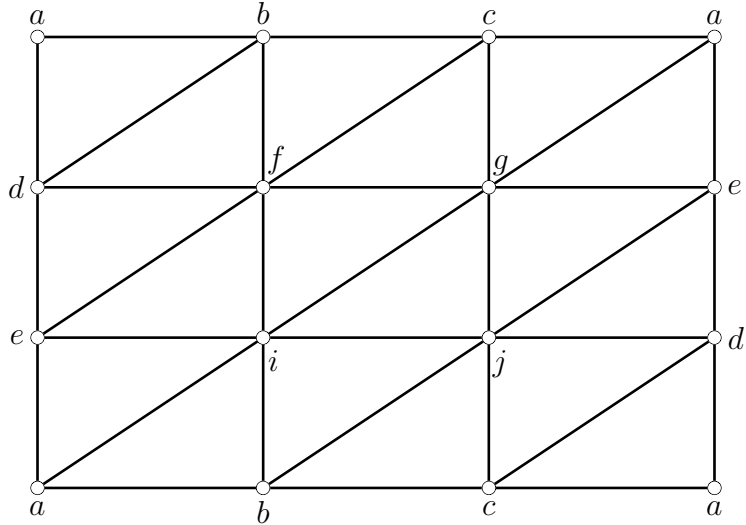

\medskip
Figure \ref{Kleinbot} shows a triangulation of  a surface
called the {\it Klein bottle\/}.

The geometric realization of the above triangulation is
obtained by pasting together the pairs of 
edges labeled $(a,d)$, $(d,e)$, $(e,a)$, and the pairs of edges
labeled  $(a,b)$, $(b,c)$, $(c,a)$.
Again, some of the gluing requires a ``twist'', since
some paired edges have opposite orientation.
Visualizing this surface in $\affreal^3$ not too difficult,
but self-intersection cannnot be avoided.

\medskip
We are now going to state a proposition characterizing the
complexes $K$ that correspond to triangulations of surfaces.
The following notational conventions will be used:
vertices (or nodes, i.e., $0$-simplices) will be denoted as $\alpha$,
edges ($1$-simplices) will be denoted as $a$, and triangles
($2$-simplices) will be denoted as $A$.
We will also denote an edge as $a = (\alpha_1 \alpha_2)$,
and a triangle as $A = (a_1 a_2 a_3)$, or as
$A = (\alpha_1 \alpha_2 \alpha_3)$, when we are interested in 
its vertices. For the moment, we do not care about the order.

\begin{prop}
\label{trianglem2}
A $2$-complex $K = (V, \s{S})$ is a triangulation 
$\mapdef{\sigma}{\s{S}}{2^{M}}$
of a surface $M$ such that $\sigma(s) = s_g$ for all $s\in \s{S}$ iff
the following properties hold:
\begin{enumerate}
\item[(D1)] 
Every edge $a$ is contained in exactly two triangles $A$.
\item[(D2)] 
For every vertex $\alpha$, the edges $a$ and triangles $A$ 
containing $\alpha$ can be arranged  as a cyclic sequence
$a_1,A_1,a_2,A_2,\ldots,A_{m-1}, a_m, A_m$, in the sense
that $a_i = A_{i-1}\cap A_i$ for all $i$, $2\leq i\leq m$,
and $a_1 = A_m\cap A_1$, with $m \geq 3$.
\item[(D3)]
$K$ is connected, in the sense that it cannot be written
as the union of two disjoint nonempty complexes.
\end{enumerate}
\end{prop}

\proof A proof can be found in
Ahlfors and Sario \cite{Ahlfors}.
The proof requires the notion of the winding number
of a closed curve in the plane with respect to a point,
and the concept of homotopy.
$\square$

\medskip
A $2$-complex $K$ which satisfies the conditions of Proposition
\ref{trianglem2} will be called a {\it triangulated complex\/},
and its geometric realization is called a {\it polyhedron\/}.
Thus, triangulated complexes are the complexes that correspond
to triangulated surfaces. Actually, it can be shown
that every surface admits some triangulation, and thus
the class of geometric realizations of the
triangulated complexes is the class of all  surfaces. 
We now give a quick presentation of homotopy, the fundamental
group, and homology groups.

\chapter{The Fundamental Group, Orientability}
\label{chap4b}
\section{The Fundamental Group}
\label{fundgroup}
If we want to somehow classify surfaces, we have to deal with
the issue of deciding when we consider two surfaces to be equivalent.
It seems reasonable to treat homeomorphic surfaces as equivalent,
but this leads to the problem of deciding when 
two surfaces are not homeomorphic, which is a very difficult problem.
One way to approach this problem is to forget some of 
the topological structure of a surface, and look for more
algebraic objects that can be associated with a surface.
For example, we can consider closed curves on a surface, and see
how they can be deformed. It is also fruitful to give
an algebraic structure to appropriate sets of  closed curves
on a surface, for example, a group structure.
Two important tools for studying  surfaces
were invented by Poincar\'e, the fundamental group, and the
homology groups. In this section, we take a look at the
fundamental group.
Roughly speaking, given a topological space
$E$ and some chosen point $a\in E$, a group $\pi(E, a)$
called the fundamental group of $E$ based at $a$ 
is associated with $(E, a)$, and to every
continuous map  $\mapdef{f}{(X, x)}{(Y, y)}$
such that $f(x) = y$, is associated a group homomorphism
$\mapdef{f_{*}}{\pi(X, x)}{\pi(Y, y)}$.
Thus, certain topological questions about the space $E$ can translated
into algebraic questions about the group $\pi(E,a)$.
This is the paradigm of algebraic topology.
In this section, we will focus on the concepts rather
than dwelve into technical details. For a thorough presentation
of the fundamental group and related concepts, the reader is referred to
Massey \cite{Massey87,Massey}, Munkres \cite{Munkrestop},
Bredon \cite{Bredon}, Dold \cite{Dold}, Fulton \cite{Fulton95},
Rotman \cite{Rotman}. We also recommend Sato \cite{Sato} for
an informal and yet very clear presentation.

\medskip
The intuitive idea behind the fundamental group is that closed paths
on a surface reflect some of the main topological properties of the surface.
Actually, the idea applies to any topological space
$E$. Let us choose some point $a$ in $E$ (a {\it base point\/}),
and consider all closed curves $\mapdef{\gamma}{[0, 1]}{E}$
based at $a$, that is, such that $\gamma(0) = \gamma(1) = a$.
We can compose closed curves $\gamma_1, \gamma_2$ based at $a$,
and consider the inverse $\gamma^{-1}$ of a closed curve,
but unfortunately, the operation of composition of closed curves
is not associative, and $\gamma\gamma^{-1}$ is not the identity
in general. In order to obtain a group structure, we define
a notion of equivalence of closed curves under continuous
deformations. Actually, such a notion can be defined
for any two paths with the same origin and extremity, and even for
continuous maps.

\begin{defin}
\label{homotop}
{\em
Given any two paths $\mapdef{\gamma_1}{[0, 1]}{E}$
and  $\mapdef{\gamma_2}{[0, 1]}{E}$ with the same intial point $a$
and the same terminal point $b$, i.e., 
such that $\gamma_1(0) = \gamma_2(0) = a$,
and  $\gamma_1(1) = \gamma_2(1) = b$, a map
$\mapdef{F}{[0, 1]\times [0, 1]}{E}$ is a {\it (path) homotopy\/}
between $\gamma_1$ and $\gamma_2$ if
$F$ is continuous, and if
\[\displaylignes{
\hfill F(t, 0) = \gamma_1(t),\hfill\cr
\hfill F(t, 1) = \gamma_2(t),\hfill\cr
}\]
for all $t\in [0, 1]$,
and
\[\displaylignes{
\hfill F(0, u) = a,\hfill\cr 
\hfill F(1, u) = b,\hfill\cr 
}\]
for all $u\in [0, 1]$.
In this case, we say that $\gamma_1 $ and $\gamma_2$ are {\it homotopic\/},
and this is denoted as $\gamma_1\approx \gamma_2$.

\medskip
Given any two continuous maps $\mapdef{f_1}{X}{Y}$ and
$\mapdef{f_2}{X}{Y}$ between two topological spaces $X$ and $Y$,
a  map $\mapdef{F}{X\times [0, 1]}{Y}$ is a {\it homotopy\/}
between $f_1$ and $f_2$ iff $F$ is continuous and
\[\displaylignes{
\hfill F(t, 0) = f_1(t),\hfill\cr
\hfill F(t, 1) = f_2(t),\hfill\cr
}\]
for all $t\in X$. We say that  $f_1$ and $f_2$ are homotopic,
and this is denoted as $f_1\approx f_2$.
}
\end{defin}

\medskip
Intuitively, a (path) homotopy $F$ between
two paths $\gamma_1$ and $\gamma_2$ from $a$ to $b$
is a continuous family of paths $F(t,u)$ from $a$ to $b$, giving a
deformation of  the path $\gamma_1$ into the path $\gamma_2$.
It is easily shown that homotopy is an equivalence relation on the
set of paths from $a$ to $b$.
A simple example of homotopy is given by reparameterizations.
A continuous nondecreasing function
$\mapdef{\tau}{[0, 1]}{[0, 1]}$ such that $\tau(0) = 0$ and
$\tau(1) = 1$ is called a {\it reparameterization\/}.
Then, given a path $\mapdef{\gamma}{[0, 1]}{E}$, the path
$\mapdef{\gamma\circ\tau}{[0, 1]}{E}$ is homotopic to 
$\mapdef{\gamma}{[0, 1]}{E}$, under the homotopy
\[(t, u) \mapsto \gamma((1 - u)t + u\tau(t)).\]
As another example, any two continuous maps
$\mapdef{f_1}{X}{\affreal^2}$ and  $\mapdef{f_2}{X}{\affreal^2}$
with range the affine plane $\affreal^2$
are homotopic under the homotopy defined such that
\[F(t, u) = (1 - u)f_1(t) + uf_2(t).\]
However, if we remove the origin from the plane $\affreal^2$,
we can find paths $\gamma_1$ and $\gamma_2$ from $(-1, 0)$ to
$(1, 0)$ that are not homotopic. For example, we can 
consider the upper half unit circle, and the lower half unit circle.
The problem is that the ``hole'' created by the missing origin
prevents continuous deformation of one path into the other.
Thus, we should expect that homotopy classes of closed
curves on a surface contain information about the presence
or absence of ``holes'' in a surface.

\medskip
It is easily verified that if $\gamma_1\approx \gamma_1'$
and $\gamma_2\approx \gamma_2'$, then
$\gamma_1\gamma_2\approx \gamma_1'\gamma_2'$, and that
$\gamma_1^{-1}\approx \gamma_1'^{-1}$. Thus, it makes sense to
define the composition and the inverse of homotopy classes.

\begin{defin}
\label{fundgrp}
{\em
Given any topological space $E$, for any choice of a point
$a\in E$ (a {\it base point\/}), 
the {\it fundamental group (or Poincar\'e group) 
$\pi(E, a)$ at the base point $a$\/}
is the set of homotopy classes of 
closed curves $\mapdef{\gamma}{[0, 1]}{E}$ such that
$\gamma(0) = \gamma(1) = a$, under the
multiplication operation  $[\gamma_1][\gamma_2] = [\gamma_1\gamma_2]$,
induced by the composition of closed paths based at $a$.
}
\end{defin}

\medskip
One actually needs to prove that the above multiplication operation
is associative, has the homotopy class of the constant path equal to $a$ 
as an identity, and that the inverse of the homotopy class
$[\gamma]$ is the class $[\gamma^{-1}]$.
The first two properties are left as an exercise, and 
the third property uses the homotopy
\[F(t, u) = \cases{\gamma(2t) & if $0\leq t\leq (1 - u)/2$;\cr
                   \gamma(1 - u) & if $(1 - u)/2 \leq t \leq (1 + u)/2$;\cr
                   \gamma(2 - 2t)& if $(1 + u)/2 \leq t \leq 1$.\cr
}\]

\medskip
As defined, the fundamental group depends on the choice of a base point.
Let us now assume that $E$ is arcwise connected (which is the case
for surfaces). Let $a$ and $b$ be any two distinct base points.
Since $E$ is arcwise connected, there is some path
$\alpha$ from $a$ to $b$. Then, to every closed curve
$\gamma$ based at $a$ corresponds a close curve 
$\gamma'= \alpha^{-1}\gamma\alpha$
based at $b$. It is easily verified that this map induces
a homomorphism $\mapdef{\varphi}{\pi(E, a)}{\pi(E, b)}$
between the groups $\pi(E, a)$ and $\pi(E, b)$.
The path $\alpha^{-1}$ from $b$ to $a$ induces a homomorphism
$\mapdef{\psi}{\pi(E, b)}{\pi(E, a)}$
between the groups $\pi(E, b)$ and $\pi(E, a)$.
Now, it is immediately verified that
$\varphi\circ\psi$ and $\psi\circ\varphi$ are both the identity,
which shows that the groups $\pi(E, a)$ and $\pi(E, b)$
are isomorphic. 

\medskip
Thus, when the space $E$ is arcwise connected, the fundamental
groups  $\pi(E, a)$ and $\pi(E, b)$ are isomorphic for any two 
points $a, b\in E$. 

\remarks
\begin{enumerate}
\item[(1)] 
The isomorphism  $\mapdef{\varphi}{\pi(E, a)}{\pi(E, b)}$
is not canonical, that is, it depends on the chosen path $\alpha$
from $a$ to $b$.
\item[(2)] 
In general, the fundamental group  $\pi(E, a)$ is not commutative.
\end{enumerate}

When $E$ is arcwise connected, we allow ourselves to refer to 
any of the isomorphic groups $\pi(E, a)$
as {\it the\/} fundamental group of $E$, and we denote any of
these groups as $\pi(E)$.

\medskip
The fundamental group  $\pi(E, a)$ is in fact one of several
homotopy groups  $\pi_n(E, a)$ associated with a space $E$,
and  $\pi(E, a)$ is often denoted as  $\pi_1(E, a)$. However,
we won't have any use for the more general homotopy groups.

\medskip
If $E$ is an arcwise connected topological space,
it may happen that some fundamental groups $\pi(E, a)$ is 
reduced to the trivial group $\{1\}$ consisting of the identity element.
It is easy to see that
this is equivalent to the fact that for any two points
$a, b\in E$, any two paths from $a$ to $b$ are homotopic, and thus,
the fundamental groups $\pi(E, a)$ are trivial for all $a\in E$.
This is an important case, which motivates the following definition.

\begin{defin}
\label{simplycon}
{\em
A topological space $E$ is {\it simply-connected\/} if
it is arcwise connected and for every $a\in E$,
the fundamental group $\pi(E, a)$ 
is the trivial one-element group.
}
\end{defin}

\medskip
For example, the plane and the sphere are simply connected,
but the torus is not simply connected (due to its hole).

\medskip
We now show that a continuous map between topological spaces
(with base points) induces a homomorphism of fundamental groups.
Given two topological spaces $X$ and $Y$, given a base
point $x$ in $X$ and a base point $y$ in $Y$,
for any continuous map $\mapdef{f}{(X, x)}{(Y, y)}$
such that $f(x) = y$, we can define a map
$\mapdef{f_{*}}{\pi(X, x)}{\pi(Y, y)}$ as follows:
\[f_{*}([\gamma]) = [f\circ \gamma],\]
for every homotopy class $[\gamma]\in \pi(X, x)$, where
$\mapdef{\gamma}{[0, 1]}{X}$ is a closed path based at $x$.

\medskip
It is easily verified that $f_{*}$ is well defined, that is,
does not depend on the choice of the closed curve $\gamma$
in the homotopy class $[\gamma]$.
It is also easily verified that $\mapdef{f_{*}}{\pi(X, x)}{\pi(Y, y)}$ is
a homomorphism of groups.
The map $f\mapsto f_{*}$ also has the following
important two properties. For any two continuous
maps  $\mapdef{f}{(X, x)}{(Y, y)}$ and  $\mapdef{g}{(Y, y)}{(Z, z)}$,
such that $f(x) = y$ and $g(y) = z$, we have
\[(g\circ f)_{*} = g_{*}\circ f_{*},\]
and if $\mapdef{Id}{(X, x)}{(X, x)}$ is the identity map,
then  $\mapdef{Id_{*}}{\pi(X, x)}{\pi(X, x)}$  is the identity
homomorphism.

\medskip
As a consequence, if $\mapdef{f}{(X, x)}{(Y, y)}$ is
a homeomorphism such that $f(x) = y$, then
$\mapdef{f_{*}}{\pi(X, x)}{\pi(Y, y)}$ is a group isomorphism.
This gives us a way of proving that two spaces are not
homeomorphic: show that for some appropriate base points
$x\in X$ and $y\in Y$, the fundamental groups
$\pi(X, x)$ and $\pi(Y, y)$ are not isomorphic.

\medskip
In general, it is difficult to determine the fundamental group
of a space. We will determine the fundamental group of $\affreal^n$
and of the punctured plane. For this, we need the concept
of the winding number of a closed curve in the plane.

\section{The Winding Number of a Closed Plane Curve}
\label{windnum}
Consider a closed curve $\mapdef{\gamma}{[0, 1]}{\affreal^2}$
in the plane, and let $z_0$ be a point not on $\gamma$.
In what follows, it is convenient to identify the plane
$\affreal^2$ with the set $\complex$ of complex numbers.
We wish to define a number $n(\gamma, z_0)$ which
counts how many times the closed curve $\gamma$ winds around $z_0$.

\medskip
We claim that there is some real number $\rho > 0$ such that
$|\gamma(t) - z_0| > \rho$ for all $t\in [0, 1]$. If not,
then for every integer $n\geq 0$, there is some $t_n\in [0, 1]$
such that $|\gamma(t_n) - z_0| \leq 1/n$.
Since $[0, 1]$ is compact, the sequence $(t_n)$ has
some convergent subsequence $(t_{n_{p}})$ having some limit $l\in [0, 1]$.
But then, by continuity of $\gamma$, we have $\gamma(l) = z_0$,
contradicting the fact that $z_0$ is not on $\gamma$.
Now, again since $[0, 1]$ is compact and $\gamma$ is continuous,
$\gamma$ is actually uniformly continuous.
Thus, there is some $\epsilon > 0$ such that
$|\gamma(t) - \gamma(u)| \leq \rho$ for all $t, u\in [0, 1]$,
with $|u - t| \leq \epsilon$. Letting $n$ be the smallest
integer such that $n\epsilon > 1$, and letting $t_i = i/n$,
for $0\leq i \leq n$, we get a subdivision of $[0, 1]$
into subintervals $[t_i, t_{i+1}]$, such that
$|\gamma(t) - \gamma(t_i)| \leq \rho$ for all $t \in [t_i, t_{i+1}]$,
with $0\leq i\leq n -1$.

\medskip
For every $i, 0\leq i \leq n-1$, if we let 
\[w_i = \frac{\gamma(t_{i+1}) - z_0}{\gamma(t_{i}) - z_0},\]
it is immediately verified that $|w_i - 1| < 1$, and thus,
$w_i$ has a positive real part. Thus, there is a unique
angle $\theta_i$ with  $-\frac{\pi}{2} < \theta_i < \frac{\pi}{2}$,
such that $w_i = \lambda_i(\cos\theta_i + i\sin\theta_i)$,
where $\lambda_i > 0$. Furthermore, 
because $\gamma$ is a closed curve,
\[\prod_{i = 0}^{n - 1}w_i = 
\prod_{i = 0}^{n - 1}\frac{\gamma(t_{i+1}) - z_0}{\gamma(t_{i}) - z_0}
= \frac{\gamma(t_n) - z_0}{\gamma(t_0) - z_0} 
= \frac{\gamma(1) - z_0}{\gamma(0) - z_0} = 1,\]
and the angle $\sum\theta_i$ is an integral multiple of $2\pi$.
Thus, for every subdivision of $[0, 1]$
into intervals $[t_i, t_{i+1}]$ such that  $|w_i - 1| < 1$, 
with $0\leq i\leq n-1$,
we define the {\it winding number $n(\gamma, z_0)$, or index, 
of $\gamma$ with respect to $z_0$\/}, as
\[n(\gamma, z_0) = \frac{1}{2\pi}\sum_{i = 0}^{i = n-1} \theta_i.\]

\medskip
Actually, in order for $n(\gamma, z_0)$ to be well defined,
we need to show that it does not depend on the subdivision of $[0, 1]$
into intervals $[t_i, t_{i+1}]$ (such that  $|w_i - 1| < 1$).
Since any two subdivisions of $[0, 1]$
into intervals $[t_i, t_{i+1}]$ can be refined into a common subdivision,
it is enough to show that nothing is changed is we replace any interval
$[t_i, t_{i+1}]$ by the two intervals $[t_i, \tau]$ and
$[\tau, t_{i+1}]$. 
Now, if $\theta_i'$ and $\theta_i''$ are the angles
associated with
\[\frac{\gamma(t_{i+1}) - z_0}{\gamma(\tau) - z_0},\]
and
\[\frac{\gamma(\tau) - z_0}{\gamma(t_i) - z_0},\]
we have 
\[\theta_i = \theta_i' + \theta_i'' + k2\pi,\]
where $k$ is some integer. However, since 
$-\frac{\pi}{2} < \theta_i < \frac{\pi}{2}$,
$-\frac{\pi}{2} < \theta_i' < \frac{\pi}{2}$, and
$-\frac{\pi}{2} < \theta_i'' < \frac{\pi}{2}$,
we must have $|k| < \frac{3}{4}$, which implies that $k = 0$,
since $k$ is an integer. This shows that 
$n(\gamma, z_0)$ is well defined.

\medskip
The next two propositions are easily shown using the above technique.
Proofs can be found in Ahlfors and Sario \cite{Ahlfors}.

\begin{prop}
\label{windlem1}
For every plane closed curve $\mapdef{\gamma}{[0, 1]}{\affreal^2}$,
for every $z_0$ not on $\gamma$, the index $n(\gamma, z_0)$ is
continuous on the complement of $\gamma$ in $\affreal^2$,
and in fact constant in each connected component of the complement
of $\gamma$. We have  $n(\gamma, z_0) = 0$ in the
unbounded component of the complement of $\gamma$.
\end{prop}

\begin{prop}
\label{windlem2}
For any two plane closed curve $\mapdef{\gamma_1}{[0, 1]}{\affreal^2}$
and $\mapdef{\gamma_2}{[0, 1]}{\affreal^2}$,
for every homotopy $\mapdef{F}{[0, 1]\times [0, 1]}{\affreal^2}$
between $\gamma_1$ and $\gamma_2$,
for every $z_0$ not on any $F(t, u)$, for all  $t, u\in [0, 1]$,
we have $n(\gamma_1,z_0)  = n(\gamma_2,z_0)$.
\end{prop}

\medskip
Proposition \ref{windlem2} shows that the index of a closed plane curve
is not changed under homotopy (provided that none the curves
involved go through $z_0$).
We can now compute the fundamental group of the punctured
plane, i.e., the plane from which a point is deleted.

\section{The Fundamental Group of the Punctured Plane}
\label{punctplane}
First, we note that the fundamental group of $\affreal^n$ is
the trivial group. Indeed, consider any closed curve
$\mapdef{\gamma}{[0, 1]}{\affreal^n}$ through $a = \gamma(0) = \gamma(1)$,
take $a$ as base point,
and let $a$ be the constant closed curve reduced to $a$.
Note that the map
\[(t, u) \mapsto (1 - u)\gamma(t)\]
is a homotopy between $\gamma$ and $a$. Thus, there is a single homotopy
class $[a]$, and $\pi(\affreal^n, a) = \{1\}$.

\medskip
The above reasoning also shows that the fundamental group of
an open ball, or a closed ball, is trivial. However, the next proposition
shows that the fundamental group of the punctured  plane
is the infinite cyclic group $\integs$.

\begin{prop}
\label{punctlem}
The fundamental group of the punctured  plane
is the infinite cyclic group $\integs$.
\end{prop}

\proof
Assume that the origin $z = 0$ is deleted from $\affreal^2 = \complex$,
and take $z = 1$ as base point.
The unit circle can be parameterized as $t\mapsto\cos t + i\sin t$,
and let $\alpha$ be the corresponding closed curve.
First of all, note that for every closed curve
$\mapdef{\gamma}{[0, 1]}{\affreal^2}$ based at $1$, there
is a homotopy (central projection) 
$\mapdef{F}{[0, 1]\times [0, 1]}{\affreal^2}$
deforming $\gamma$ into a curve $\beta$ lying on the unit circle.
By uniform continuity, any such curve $\beta$ can be decomposed
as $\beta = \beta_1\beta_2\cdots\beta_n$, where
each $\beta_k$ either does not pass through $z = 1$, or does not
pass through $z = -1$. It is also easy to see that $\beta_k$
can deformed  into one of the circular arcs $\delta_k$ between its endpoints.
For all k$, 2\leq k\leq n$, let $\sigma_k$ be one of the circular arcs
from $z = 1$ to the initial point of $\delta_k$, 
and let $\sigma_1 = \sigma_{n+1} = 1$. We have
\[\gamma\approx (\sigma_1\delta_1\sigma_2^{-1})\cdots
(\sigma_n\delta_n\sigma_{n+1}^{-1}),\]
and it is easily seen that each arc $\sigma_k\delta_k\sigma_{k+1}^{-1}$
is homotopic either to $\alpha$, or $\alpha^{-1}$, or $1$.
Thus, $\gamma\approx\alpha^m$, for some integer $m\in\integs$.

\medskip
It remains to prove that $\alpha^m$ is not homotopic to $1$ for $m\not= 0$.
This is where we use Proposition \ref{windlem2}. Indeed, it is immediate that
$n(\alpha^m, 0) = m$, and $n(1, 0) = 0$, and thus
$\alpha^m$ and $1$ are not homotopic when $m\not= 0$.
But then, we have shown that the homotopy classes are in bijection
with the set of integers.
$\square$

\medskip
The above proof also applies to a cicular annulus, closed or open,
and to a circle. In particular, the circle is not simply connected.

\medskip
We will need to define what it means for a surface
to be  orientable. Perhaps surprisingly,
a rigorous definition is not so easy to obtain,
but can be given using the notion of degree of a homeomorphism
from a plane region. First, we need to define the degree of a map
in the plane.

\section{The Degree of a Map in the Plane}
\label{degreemap}
Let $\mapdef{\varphi}{D}{\complex}$ be a continuous function
to the plane, where the plane is
viewed as the set $\complex$ of complex numbers,
and with domain some open set $D$ in $\complex$.
We say that $\varphi$ is {\it regular at $z_0\in D$\/} if
there is some open set $V\subseteq D$ containing $z_0$ such
that $\varphi(z) \not= \varphi(z_0)$, for all $z\in V$.
Assuming that $\varphi$ is regular at $z_0$, we will define
the {\it degree of $\varphi$ at $z_0$\/}.

\medskip
Let $\Omega$ be a punctured open disk $\{z\in V \ |\  |z - z_0| < r\}$
contained in $V$. Since $\varphi$ is regular at $z_0$,
it maps $\Omega$ into the punctured plane $\Omega'$
obtained by deleting $w_0 = \varphi(z_0)$.
Now, $\varphi$ induces a homomorphism 
$\mapdef{\varphi_{*}}{\pi(\Omega)}{\pi(\Omega')}$.
From Proposition \ref{punctlem}, both groups $\pi(\Omega)$
and $\pi(\Omega')$ are isomorphic
to $\integs$. Thus, it is easy to determine exactly what the
homorphism $\varphi_{*}$ is. We know that $\pi(\Omega)$ is generated by
the homotopy class of some circle $\alpha$ in $\Omega$ with center $a$,
and that  $\pi(\Omega')$ is generated by
the homotopy class of some circle $\beta$ in $\Omega'$ with center 
$\varphi(a)$. If $\varphi_{*}([\alpha]) = [\beta^d]$,
then the homomorphism $\varphi_{*}$ is completely determined.
If $d = 0$, then $\pi(\Omega') = 1$, and if $d\not= 0$, then
$\pi(\Omega')$ is the infinite cyclic subgroup generated by
the class of $\beta^d$. We let $d$ be the {\it degree of $\varphi$ at
$z_0$\/}, and we denote it as $d(\varphi)_{z_{0}}$.
It is easy to see that this definition does not depend on the choice
of $a$ (the center of the circle $\alpha$)
in $\Omega$, and thus, does not depend on $\Omega$.

\medskip
Next, if we have a second mapping $\psi$ regular at $w_0 = \varphi(z_0)$,
then $\psi\circ\varphi$ is regular at $z_0$, and it is immediately
verified that 
\[d(\psi\circ\varphi)_{z_{0}} = d(\psi)_{w_{0}}d(\varphi)_{z_{0}}.\] 

\medskip
Let us now assume that $D$ is a region (a connected open set),
and that $\varphi$ is a homeomorphism between $D$ and $\varphi(D)$.
By a theorem of Brouwer (the invariance of domain), 
it turns out that $\varphi(D)$ is also open, and thus, we can
define the degree of the inverse mapping $\varphi^{-1}$, and
since the identity clearly has degree $1$, we get that
$d(\varphi)d(\varphi^{-1}) = 1$, which shows that
$d(\varphi)_{z_{0}} = \pm 1$.

\medskip
In fact, Ahlfors and Sario \cite{Ahlfors} prove that
if $\varphi(D)$ has a nonempty interior, then the degree
of $\varphi$ is constant on $D$.
The proof is not difficult, but not very instructive.

\begin{prop}
\label{deglem}
Given a region $D$ in the plane, for every 
homeomorphism $\varphi$ between $D$ and $\varphi(D)$, if
$\varphi(D)$ has a nonempty interior, then the degree 
$d(\varphi)_z$ is constant for all $z\in D$, and in fact, 
$d(\varphi) = \pm 1$.
\end{prop}

\medskip
When $d(\varphi) = 1$ in Proposition \ref{deglem}, we say that
$\varphi$ is {\it sense-preserving\/}, and when 
$d(\varphi) = -1$, we say that
$\varphi$ is {\it sense-reversing\/}.
We can now define the notion of orientability.

\section{Orientability of a Surface}
\label{orient}
Given a surface $F$, we will call a region $V$ on $F$
a {\it planar region\/} if there is a homeomorphism
$\mapdef{h}{V}{U}$ from $V$ onto an open set in the plane.
From Proposition \ref{deglem}, the homeomorphisms
$\mapdef{h}{V}{U}$ can be divided into two classes, by
defining two such homeomorphisms $h_1, h_2$ as equivalent
iff $h_1\circ h_2^{-1}$ has degree $1$, i.e., is sense-preserving.
Observe that for any $h$ as above, if $\overline{h}$
is obtained from $h$ by conjugation (i.e., for every
$z\in V$, $\overline{h}(z) = \overline{h(z)}$, the
complex conjugate of $h(z)$), then $d(h\circ \overline{h^{-1}}) = -1$,
and thus $h$ and $\overline{h}$ are in different classes.
For any other such map $g$, either $h\circ g^{-1}$ or
$\overline{h}\circ g^{-1}$ is sense-preserving, and thus, there
are exactly two equivalence classes.

\medskip
The choice of one of the two classes of homeomorphims $h$ as above,
constitutes an {\it orientation\/} of $V$.
An orientation of $V$ induces an orientation on any subregion $W$
of $V$, by restriction.
If $V_1$ and $V_2$ are two planar regions, and these regions
have received an orientation, we say that these orientations
are {\it compatible\/} if they induce the same orientation
on all common subregions of $V_1$ and $V_2$.

\begin{defin}
\label{orientedf}
{\em
A surface $F$ is {\it orientable\/} if it is possible to assign
an orientation to all planar regions in such a way that
the orientations of any two overlapping planar regions
are compatible.
}
\end{defin}

\medskip
Clearly, orientability is preserved by homeomorphisms.
Thus, there are two classes of surfaces, the orientable surfaces,
and the nonorientable surfaces. An example of a nonorientable
surface is the Klein bottle. Because we defined a surface as
being connected, note that an orientable surface has exactly
two orientations.
It is also easy to see that to orient a surface, it is enough
to orient all planar regions in some open covering 
of the surface by planar regions.

\medskip
We will also need to consider bordered surfaces.

\section{Bordered Surfaces}
\label{bordsuf}
Consider a torus, and cut out a finite number of small disks
from its surface. The resulting space is no longer a surface,
but certainly of geometric interest. It is a surface with
boundary, or bordered surface. In this section, we extend our
concept of surface to handle this more general class of
bordered surfaces. In order to do so, we need to allow coverings
of surfaces using a richer class of open sets. This is achieved by
considering the open subsets of the half-space, in the subset
topology.

\begin{defin}
\label{boundman}
{\em
The {\it half-space\/} $\mathbb{H}^m$ is the subset of $\reals^m$
defined as the set
\[\{(x_1,\ldots,x_m) \ |\ x_i\in\reals,\, x_m \geq 0\}.\]
For any $m\geq 1$, a {\it (topological) $m$-manifold
with boundary\/} is
a second-countable, topological Hausdorff space $M$, together
with  an open cover $(U_i)_{i\in I}$ of open sets and
a family $(\varphi_i)_{i\in I}$ of homeomorphisms
$\mapdef{\varphi_i}{U_i}{\Omega_i}$, where each $\Omega_i$
is some open subset of $\mathbb{H}^m$ in the subset topology. 
Each pair $(U, \varphi)$
is called a {\it  coordinate system\/}, or
{\it chart\/}, of $M$, 
each homeomorphism $\mapdef{\varphi_i}{U_i}{\Omega_i}$ is called
a {\it coordinate map\/},
and its inverse  $\mapdef{\varphi^{-1}_{i}}{\Omega_i}{U_i}$ is called a 
{\it  parameterization\/} of $U_i$.
The family  $(U_i, \varphi_i)_{i\in I}$ is often called an
{\it atlas\/} for $M$.
A {\it (topological) bordered surface\/} is a connected $2$-manifold
with boundary.
}
\end{defin}

\medskip
Note that an $m$-manifold is also an $m$-manifold with boundary.

\medskip
If $\mapdef{\varphi_i}{U_i}{\Omega_i}$ is some
homeomorphism onto some open set $\Omega_i$ of $\mathbb{H}^m$ 
in the subset topology, some $p\in U_i$ may be mapped into
$\reals^{m-1}\times\reals_+$, or into the ``boundary''
$\reals^{m-1}\times\{0\}$ of $\mathbb{H}^m$.
Letting $\dBd\mathbb{H}^m = \reals^{m-1}\times\{0\}$, it can be
shown using homology, that if
some   coordinate map $\varphi$ defined on $p$
maps $p$ into $\dBd\mathbb{H}^m$, then
every coordinate map $\psi$ defined on $p$
maps $p$ into $\dBd\mathbb{H}^m$.
For $m = 2$, Ahlfors and Sario prove it using Proposition \ref{deglem}.

\medskip
Thus, $M$ is the disjoint union of two sets
$\dBd M$ and $\dInt M$, where $\dBd M$ is the
subset consisting of all points $p\in M$ that are mapped
by some (in fact, all)  coordinate map $\varphi$ defined on $p$
into $\dBd\mathbb{H}^m$, and where 
$\dInt M = M - \dBd M$. The set $\dBd M$ is called
the {\it boundary\/} of $M$, and the set $\dInt M$ is called
the {\it interior\/} of $M$, even though this terminology
clashes with some prior topological definitions.
A good example of a bordered surface is the M\"obius strip.
The boundary of the M\"obius strip is a circle.

\medskip
The boundary $\dBd M$ of $M$
may be empty, but $\dInt M$ is nonempty. Also, it can be shown
using homology, that the integer $m$ is unique.
It is clear that $\dInt M$ is open, and an $m$-manifold, 
and that $\dBd M$ is closed. If $p\in \dBd M$, and 
$\varphi$ is some coordinate map  defined on $p$,
since $\Omega=\varphi(U)$ is an open subset of $\dBd\mathbb{H}^m$, 
there is some open half ball $B_{o+}^{m}$ centered
at $\varphi(p)$ and contained in  $\Omega$ which intersects
$\dBd\mathbb{H}^m$ along an open ball $B_{o}^{m-1}$,
and if we consider $W = \varphi^{-1}(B_{o+}^{m})$,
we have an open subset of $M$ containing $p$ which is mapped
homeomorphically onto $B_{o+}^{m}$ in such that way
that every point in $W\cap \dBd M$ is mapped
onto the open ball  $B_{o}^{m-1}$. Thus, it is easy to see that
$\dBd M$ is an $(m-1)$-manifold.

\medskip
In particular, in the case $m = 2$,
the boundary $\dBd M$ is a union of curves
homeomorphic either to circles of to open line segments.
In this case, if $M$ is connected but not a surface, 
it is easy to see that $M$ is the topological closure of $\dInt M$.
We also claim that $\dInt M$ is connected, i.e. a surface.
Indeed, if this was not so, we could write 
$\dInt M = M_1\cup M_2$, for two nonempty disjoint sets
$M_1$ and $M_2$. But then, we have
$M = \overline{M_1}\cup \overline{M_2}$,
and since $M$ is connected, there is some $a\in \dBd M$ also in
$\overline{M_1}\cap \overline{M_2}\not= \emptyset$.
However, there is some open set $V$ containing $a$
whose intersection with $M$
is homeomorphic with an open half-disk, and thus connected.
Then, we have
\[V\cap M = (V\cap M_1) \cup (V\cap M_2),\] 
with $V\cap M_1$ and $V\cap M_2$ open in $V$, contradicting
the fact that $M\cap V$ is connected.
Thus, $\dInt M$ is a surface. 

\medskip
When the boundary $\dBd M$ 
of a bordered surface $M$ is empty, $M$ is just a surface.
Typically, when we refer to a bordered surface, we mean
a bordered surface with a nonempty border, and otherwise,
we just say surface.

\medskip
A bordered surface $M$ is orientable iff its interior 
$\dInt M$ is orientable. It is not difficult to show that
an orientation of $\dInt M$ induces an orientation of
the boundary $\dBd M$.
The components of the boundary $\dBd M$ are called
{\it contours\/}.

\medskip
The concept of triangulation of a bordered surface is
identical to the concept defined for a surface in
Definition \ref{triangul}, and Proposition \ref{trianglem1}
also holds.
However, a small change needs to made to Proposition
\ref{trianglem2}, see Ahlfors and Sario \cite{Ahlfors}.

\begin{prop}
\label{trianglem3}
A $2$-complex $K = (V, \s{S})$ is a triangulation 
$\mapdef{\sigma}{\s{S}}{2^{M}}$
of a bordered surface $M$ such that $\sigma(s) = s_g$ for all $s\in \s{S}$ iff
the following properties hold:
\begin{enumerate}
\item[(D1)] 
Every edge $a$ such that $a_g$ contains some point in the interior
$\dInt M$ of $M$ is contained in exactly two triangles $A$.
Every edge $a$ such that $a_g$ is inside the border $\dBd M$ of $M$
is contained in exactly one triangle $A$.
The border $\dBd M$ of $M$ consists of those
$a_g$ which belong to only one $A_g$.
A {\it border vertex or border edge\/} is a simplex $\sigma$ such that
$\sigma_g\subseteq \dBd M$.
\item[(D2)] 
For every non-border vertex $\alpha$,
the edges $a$ and triangles $A$ 
containing $\alpha$ can be arranged  as a cyclic sequence
$a_1,A_1,a_2,A_2,\ldots,A_{m-1}, a_m, A_m$, in the sense
that $a_i = A_{i-1}\cap A_i$ for all $i$, with $2\leq i\leq m$,
and $a_1 = A_m\cap A_1$, with $m \geq 3$.
\item[(D3)] 
For every border vertex $\alpha$,
the edges $a$ and triangles $A$ 
containing $\alpha$ can be arranged  in a sequence
$a_1,A_1,a_2,A_2,\ldots,A_{m-1},a_{m},A_{m},a_{m+1}$, with
$a_i = A_i\cap A_{i-1}$ for of all $i$,
with $2\leq i\leq m$, where $a_1$ and $a_{m+1}$
are border vertices only contained in $A_1$ and $A_{m}$
respectively. 
\item[(D4)] 
$K$ is connected, in the sense that it cannot be written
as the union of two disjoint nonempty complexes.
\end{enumerate}
\end{prop}

A $2$-complex $K$ which satisfies the conditions of Proposition
\ref{trianglem3} will also be called a  
{\it bordered triangulated $2$-complex\/},
and its geometric realization a {\it bordered polyhedron\/}.
Thus, bordered triangulated $2$-complexes are the complexes that correspond
to triangulated bordered surfaces. Actually, it can be shown
that every bordered surface admits some triangulation, and thus
the class of geometric realizations of the
bordered triangulated $2$-complexes is the class of all bordered surfaces. 

\medskip
We will now give a brief presentation of 
simplicial and singular homology, but first,
we need to review some facts about finitely generated
abelian groups.

\chapter{Homology Groups}
\label{chap5}
\section{Finitely Generated Abelian Groups}
\label{abgroup}
An abelian group is a commutative group. We will denote
the identity element of an abelian group as $0$, and the
inverse of an element $a$ as $-a$.
Given any natural number $n\in\natnums$, we denote
\[\underbrace{a + \cdots + a}_n\]
as $na$, and let $(-n)a$ be defined as $n(-a)$ (with $0a = 0$).
Thus, we can make sense of finite sums of the form
$\sum n_ia_i$, where $n_i\in\integs$.
Given an abelian group $G$ and a family $A = (a_j)_{j\in J}$ of
elements $a_j\in G$, we say that $G$ is {\it  generated by A\/} if 
every $a\in G$ can be written (in possibly more than one way) as
\[a = \sum_{i\in I} n_ia_i,\]
for some finite subset $I$ of $J$, and  some $n_i\in\integs$.
If $J$ is finite, we say that $G$ is {\it finitely generated by A\/}.
If every  $a\in G$ can be written in a {\it unique manner\/} as
\[a = \sum_{i\in I} n_ia_i\]
as above, we say that $G$ is {\it freely generated by A\/}, and we call
$A$ a {\it basis of $G$\/}.
In this case, it is clear that the  $a_j$ are all distinct.
We also have the following familiar property.

\medskip
If $G$ is a free abelian group generated by $A = (a_j)_{j\in J}$,
for every abelian group $H$, for every function
$\mapdef{f}{A}{H}$, there is a unique homomorphism
$\mapdef{\widehat{f}}{G}{H}$, such that $\widehat{f}(a_j) = f(a_j)$,
for all $j\in J$.

\remark
If $G$ is a free  abelian group,
one can show that the cardinality of all bases is the same.
When $G$ is free and finitely generated by $(a_1, \ldots, a_n)$,
this can be proved as follows.
Consider the quotient of the group $G$ modulo the subgroup
$2G$ consisting of all elements of the form $g + g$, where $g\in G$.
It is immediately verified that each coset of $G/2G$ is of the form
\[\epsilon_1a_1 + \cdots + \epsilon_na_n + 2G,\]
where $\epsilon_i = 0$ or $\epsilon_i = 1$,
and thus, $G/2G$ has $2^n$ elements. Thus, $n$ only depends
on $G$. The number $n$ is called the {\it dimension\/} of $G$.
\endremark

Given a family $A = (a_j)_{j\in J}$, we will need to construct a free abelian
group generated by $A$. This can be done easily as follows.
Consider the set $F(A)$ of all functions
$\mapdef{\varphi}{A}{\integs}$, such that $\varphi(a) \not= 0$
for only finitely many $a\in A$. We define addition on $F(A)$
pointwise, that is,  $\varphi + \psi$ is the function
such that  $(\varphi + \psi)(a) = \varphi(a) + \psi(a)$, for
all $a\in A$. 

\medskip
It is immediately verified that $F(A)$ is an abelian group, and
if we identify each $a_j$ with the function
$\mapdef{\varphi_j}{A}{\integs}$, such that $\varphi_j(a_j) = 1$,
and $\varphi_j(a_i) = 0$ for all $i\not= j$, it is clear
that $F(A)$ is freely generated by $A$.
It is also clear that every $\varphi\in F(A)$ can be uniquely
written as
\[\varphi = \sum_{i\in I} n_i\varphi_{i},\]
for some finite subset $I$ of $J$ such that  $n_i = \varphi(a_i)\not= 0$. 
For notational simplicity, we write $\varphi$ as
\[\varphi = \sum_{i\in I} n_ia_{i}.\]

\medskip
Given an abelian group $G$, for any $a\in G$, we say that
$a$ has {\it finite order\/} if there is some $n\not= 0$
in $\natnums$ such that $na = 0$. If $a\in G$ has finite order, there
is a least $n \not=0$ in $\natnums$ 
such that $na = 0$, called the {\it order of $a$\/}.
It is immediately verified that the subset $T$ of $G$ consisting of
all elements of finite order is a subroup of $G$, called the
{\it torsion subgroup of $G$\/}. When $T = \{0\}$, we
say that $G$ is {\it torsion-free\/}. One should be careful
that a torsion-free abelian group is not necessarily free.
For example, the field $\rats$ of rationals is torsion-free,
but not a free abelian group. 

\medskip
Clearly, the map $(n, a) \mapsto na$ from $\integs\times G$ to $G$
satisfies the properties
\[\eqaligneno{
(m + n)a &= ma + na,\cr
m(a + b) &= ma + nb,\cr
(mn)a &= m(na),\cr
1a &= a,\cr
}\]
which hold in vector spaces. However, $\integs$ is not a field.
The abelian group $G$ is just what is called a {\it $\integs$-module\/}.
Nevertheless, many concepts defined for vector spaces transfer
to $\integs$-modules. For example, given an abelian group $G$
and some subgroups $H_1, \ldots, H_n$, we can define the {\it (internal) sum\/}
\[H_1 + \cdots + H_n\]
of the $H_i$ as the abelian group consisting of all
sums of the form $a_1 + \cdots + a_n$, where $a_i\in H_i$.
If in addition, $G = H_1 + \cdots + H_n$ and 
$H_i\cap H_j = \{0\}$ for all $i, j$, with $i\not= j$, we say that
$G$ is the {\it direct sum of the $H_i$\/}, and this is denoted as
\[G = H_1 \oplus \cdots \oplus H_n.\]
When $H_1 = \ldots = H_n = H$, we abbreviate $H\oplus\cdots\oplus H$ as $H^n$. 
Homomorphims between abelian groups are
$\integs$-linear maps. We can also talk about linearly
independent families in $G$, except that the scalars are
in $\integs$. The {\it rank\/} of an abelian group is the 
maximum of the sizes of linearly independent families in $G$. 
We can also define (external) direct sums.

\medskip
Given a family $(G_i)_{i\in I}$ of abelian groups, the
{\it (external) direct sum\/} $\bigoplus_{i\in I} G_i$ is the
set of all function $\mapdef{f}{I}{\bigcup_{i\in I} G_i}$
such that $f(i) \in G_i$, for all $\in I$, and $f(i) = 0$
for all but finitely many $i\in I$. An element
$f\in \bigoplus_{i\in I} G_i$ is usually denoted as $(f_i)_{i\in I}$.
Addition is defined component-wise, that is, given two
functions $f = (f_i)_{i\in I}$ and $g = (g_i)_{i\in I}$
in $\bigoplus_{i\in I} G_i$,
we define $(f + g)$ such that 
\[(f + g)_i = f_i + g_i,\]
for all $i\in I$. It is immediately verified that
$\bigoplus_{i\in I} G_i$ is an abelian group.
For every $i\in I$, there is an injective homomorphism
$\mapdef{in_i}{G_i}{\bigoplus_{i\in I} G_i}$,
defined such that for every $x\in G_i$, $in_i(x)(i) = x$, and 
$in_i(x)(j) = 0$ iff $j \not= i$.
If $G = \bigoplus_{i\in I} G_i$ is an external direct sum,
it is immediately verified that
$G = \bigoplus_{i\in I} in_i(G_i)$, as an internal direct sum.
The difference is that $G$ must have been already defined
for an internal direct sum to make sense. 
For notational simplicity, we will usually identify $in_i(G_i)$ with $G_i$.

\medskip
The structure of finitely generated abelian groups can be 
completely described. Actually, the following result is a special case
of the structure theorem for finitely generated modules over a principal
ring. Recall that $\integs$ is a principal ring, which means that
every ideal $\s{I}$ in $\integs$ is of the form $d\integs$, for
some $d\in \natnums$.
For the sake of completeness, we present the following result,
whose neat proof is due to Pierre Samuel.

\begin{prop}
\label{finab}
Let $G$ be a free abelian group
finitely generated by $(a_1, \ldots, a_n)$, and let $H$ be
any subroup of $G$. Then, $H$ is a free abelian group, and there is
a basis $(e_1, ..., e_n)$ of $G$, some $q\leq n$, and some
positive natural numbers $n_1, \ldots, n_q$, such that
$(n_1e_1,\ldots, n_qe_q)$ is a basis of $H$, and
$n_i$ divides $n_{i+1}$ for all $i$, with $1\leq i\leq q-1$.
\end{prop}

\proof The proposition is trivial when $H = \{0\}$, and thus, we
assume that $H$ is nontrivial. Let $L(G,\integs)$ we the
set of homomorphisms from $G$ to $\integs$. For any $f\in L(G,\integs)$,
it is immediately verified that $f(H)$ is an ideal in $\integs$.
Thus, $f(H) = n_h\integs$, for some $n_h\in \natnums$, since
every ideal in $\integs$ is a principal ideal.
Since $\integs$ is finitely generated, any nonempty family
of ideals has a maximal element, and
let $f$ be a homomorphism such that $n_h\integs$ is a maximal ideal
in $\integs$. Let $\mapdef{\pi}{G}{\integs}$
be the $i$-th projection, i.e., $\pi_i$ is defined such that
$\pi_i(m_1a_1 + \cdots + m_na_n) = m_i$. It is clear that $\pi_i$
is a homomorphism, and since $H$ is nontrivial, one of the
$\pi_i(H)$ is nontrivial, and $n_h\not= 0$. There is some
$b\in H$ such that $f(b) = n_h$.

\medskip
We claim that for every $g\in  L(G,\integs)$, the number $n_h$
divides $g(b)$. Indeed, if $d$ is the gcd of $n_h$ and $g(b)$,
by the Bezout identity, we can write
\[d = rn_h + sg(b),\]
for some $r, s\in\integs$, and thus
\[d = rf(b) + sg(b) = (rf + sg)(b).\]
However, $rf + sg\in L(G,\integs)$, and thus,
\[n_h\integs \subseteq d\integs \subseteq (rf + sg)(H),\]
since $d$ divides $n_h$, and by maximality of $n_h\integs$,
we must have $n_h\integs = d\integs$, which implies that
$d = n_h$, and thus, $n_h$ divides $g(b)$.
In particular, $n_h$ divides each $\pi_i(b)$, and let
$\pi_i(b) = n_hp_i$, with $p_i\in\integs$.

\medskip
Let $a = p_1a_1 + \cdots + p_na_n$. Note that
\[b = \pi_1(b)a_1 + \cdots + \pi_n(b)a_n = n_hp_1a_1 + \cdots + n_hp_na_n,\]
and thus, $b = n_ha$.
Since $n_h = f(b) = f(n_ha) = n_hf(a)$, and since
$n_h\not= 0$, we must have $f(a) = 1$.

\medskip
Next, we claim that
\[G = a\integs \oplus f^{-1}(0),\]
and
\[H = b\integs \oplus (H\cap f^{-1}(0)),\]
with $b = n_ha$.

\medskip
Indeed, every $x\in G$ can be written as 
\[x = f(x)a + (x - f(x)a),\]
and since $f(a) = 1$, we have
$f(x - f(x)a) = f(x) - f(x)f(a) = f(x) - f(x) = 0$.
Thus, $G = a\integs + f^{-1}(0)$.
Similarly, for any $x\in H$, we have $f(x) = rn_h$, for some
$r\in\integs$, and thus,
\[x = f(x)a + (x - f(x)a) = rn_ha + (x - f(x)a) = 
rb + (x - f(x)a),\]
we still have $x - f(x)a\in f^{-1}(0)$, and clearly,
$x - f(x)a = x - rn_ha = x - rb \in H$, since $b\in H$.
Thus, $H = b\integs + (H\cap f^{-1}(0))$.

\medskip
To prove that we have a direct sum, it is enough to prove that
$a\integs \cap f^{-1}(0) = \{0\}$.
For any $x = ra \in a\integs$, if $f(x) = 0$,
then $f(ra) = rf(a) = r = 0$, since $f(a) = 1$, and thus,
$x = 0$. Therefore, the sums are direct sums.

\medskip
We can now prove that $H$ is a free abelian group by induction
on the size $q$ of a maximal linearly independent family for $H$.
If $q = 0$, the result is trivial. Otherwise, since
\[H = b\integs \oplus (H\cap f^{-1}(0)),\]
it is clear that $H\cap f^{-1}(0)$ is a subgroup of $G$
and that every maximal linearly independent family in 
$H\cap f^{-1}(0)$ has at most $q-1$ elements. By the induction
hypothesis, $H\cap f^{-1}(0)$ is a free abelian group, and by adding
$b$ to a basis of $H\cap f^{-1}(0)$, we obtain a basis for $H$,
since the sum is direct.

\medskip
The second part is shown by induction on the dimension $n$ of $G$.
The case $n = 0$ is trivial. Otherwise, since
\[G = a\integs \oplus f^{-1}(0),\]
and since by the previous argument,
$f^{-1}(0)$ is also free, it is easy to see that
$f^{-1}(0)$ has dimension $n-1$.
By the induction hypothesis applied to its subgroup
$H\cap f^{-1}(0)$, there is  a basis
$(e_2, \ldots, e_n)$ of $f^{-1}(0)$, some $q\leq n$,
and some positive natural numbers
$n_2, \ldots, n_q$, such that,
$(n_2e_2,\ldots, n_qe_q)$ is a basis of $H\cap f^{-1}(0)$, and
$n_i$ divides $n_{i+1}$ for all $i$, with $2\leq i\leq q-1$.
Let $e_1 = a$, and $n_1 = n_h$, as above.
It is clear that $(e_1, \ldots, e_n)$ is a basis of $G$,
and that that $(n_1e_1,\ldots, n_qe_q)$ is a basis of $H$, since
the sums are direct, and $b = n_1e_1 = n_ha$.
It remains to show that $n_1$ divides $n_2$.
Consider the homomorphism $\mapdef{g}{G}{\integs}$
such that $g(e_1) = g(e_2) = 1$, and $g(e_i) = 0$,
for all $i$, with $3\leq i\leq n$.
We have $n_h = n_1 = g(n_1e_1) = g(b) \in g(H)$, and thus,
$n_h\integs \subseteq g(H)$. Since $n_h\integs$ is maximal,
we must have $g(H) = n_h\integs = n_1\integs$.
Since $n_2 = g(n_2e_2) \in g(H)$, we have $n_2\in n_1\integs$,
which shows that $n_1$ divides $n_2$. 
$\square$

\medskip
Using Proposition \ref{finab}, we can also show the following
useful result.

\begin{prop}
\label{finab2}
Let $G$ be a finitely generated abelian group.
There is some natural number $m\geq 0$ 
and some positive natural numbers $n_1, \ldots, n_q$, such that $H$ is
isomorphic to the direct sum
\[\integs^m \oplus \integs/n_1\integs\oplus\cdots\oplus\integs/n_q\integs,\]
and where $n_i$ divides $n_{i+1}$ for all $i$, with $1\leq i\leq q-1$.
\end{prop}

\proof Assume that $G$ is generated by $A = (a_1,\ldots,a_n)$,
and let $F(A)$ be the free abelian group generated by $A$.
The inclusion map $\mapdef{i}{A}{G}$ can be extended
to a unique homomorphism
$\mapdef{f}{F(A)}{G}$ which is surjective since $A$ generates $G$,
and thus, $G$ is isomorphic to $F(A)/f^{-1}(0)$.
By Proposition \ref{finab}, $H = f^{-1}(0)$ is a free abelian group, and there is
a basis $(e_1, ..., e_n)$ of $G$, some $p\leq n$, and some
positive natural numbers $k_1, \ldots, k_p$, such that
$(k_1e_1,\ldots, k_pe_p)$ is a basis of $H$, and
$k_i$ divides $k_{i+1}$ for all $i$, with $1\leq i\leq p-1$.
Let $r$, $0\leq r\leq p$, be the largest natural number
such that $k_1 = \ldots = k_r = 1$, rename $k_{r+i}$
as $n_i$, where $1\leq i\leq p - r$, and let $q = p - r$. 
Then, we can write
\[H = \integs^{p-q}\oplus n_{1}\integs\oplus\cdots\oplus n_{q}\integs,\]
and since $F(A)$ is isomorphic to $\integs^n$, it is easy to verify that
$F(A)/H$ is isomorphic to
\[Z^{n - p}\oplus 
\integs/n_{1}\integs\oplus\cdots\oplus\integs/n_{q}\integs,\]
which proves the proposition.
$\square$

\medskip
Observe that $\integs/n_1\integs\oplus\cdots\oplus\integs/n_q\integs$
is the torsion subgroup of $G$. Thus, as a corollary of Proposition \ref{finab2},
we obtain the fact that every finitely generated abelian group $G$ is
a direct sum $G = Z^m\oplus T$, where $T$ is the torsion subroup
of $G$, and $Z^m$ is the free abelian group of dimension $m$.
It is  easy to verify that $m$ is the rank (the maximal dimension
of linearly independent sets in $G$) of $G$, and it is called
the {\it Betti number\/} of $G$. 
It can also be shown that  $q$, and the $n_i$, only depend on $G$.

\medskip
One more result will be needed to compute the homology groups
of (two-dimensional) polyhedras. 
The proof is not difficult and can be found 
in most books (a version is given in Ahlfors and Sario \cite{Ahlfors}).
Let us denote the rank of an abelian group $G$ as $r(G)$.

\begin{prop}
\label{finab3}
If 
\[\shortexact{E}{f}{F}{g}{G}\]
is a short exact sequence of homomorphisms of abelian groups
and $F$ has finite rank, then
$r(F) = r(E) + r(G)$. In particular, if
$G$ is an abelian group of finite rank and  $H$ is a subroup of $G$,
then  $r(G) = r(H) + r(G/H)$.
\end{prop}

\medskip
We are now ready to define the simplicial and the singular homology groups.

\section{Simplicial and Singular Homology}
\label{simphomolog}
There are several kinds of homology theories.
In this section, we take a quick look at two such theories,
simplicial homology,
one of the most computational theories, and singular homology theory,
one of the most general and yet fairly intuitive.
For a comprehensive treatment of homology and algebraic
topology in general, we refer the reader to
Massey \cite{Massey}, Munkres \cite{Munkresalg},
Bredon \cite{Bredon}, Fulton \cite{Fulton95}, Dold \cite{Dold},
Rotman \cite{Rotman}, Amstrong \cite{Amstrong}, and
Kinsey \cite{Kinsey}.
An excellent overview of algebraic topology, 
following a more intuitive approach, is presented in Sato \cite{Sato}.

\medskip
Let $K = (V,\s{S})$ be a complex. 
The essence of simplicial homology
is to associate some abelian groups $H_p(K)$ with $K$.
This is done by first defining some free abelian groups
$C_p(K)$ made out of oriented $p$-simplices.
One of the main new ingredients is that every oriented
$p$-simplex $\sigma$ is assigned a {\it boundary\/} $\partial_{p}\sigma$.
Technically, this is achieved by defining homomorphisms
\[\mapdef{\partial_{p}}{C_{p}(K)}{C_{p-1}(K)},\]
with the property that $\partial_{p-1}\circ\partial_{p} = 0$.
Letting $Z_p(K)$ be the kernel of $\partial_p$, and
\[B_p(K) = \partial_{p+1}(C_{p+1}(K))\] 
be the image of $\partial_{p+1}$
in $C_p(K)$, 
since $\partial_{p}\circ\partial_{p+1} = 0$, the group $B_p(K)$ is
a subgroup of the group $Z_p(K)$, and we define the homology group
$H_p(K)$ as the quotient group 
\[H_p(K) = Z_p(K)/B_p(K).\]
What makes the homology groups of a complex interesting, is that
they only depend on the geometric realization $K_g$ of the complex $K$,
and not on the various complexes representing $K_g$.
Proving this fact requires relatively hard work, and we refer
the reader to Munkres \cite{Munkresalg} 
or Rotman \cite{Rotman}, for a proof.

\medskip
The first step in defining simplicial homology groups is to
define oriented simplices. Given a complex $K = (V, \s{S})$, recall that
an $n$-simplex is a subset $\sigma = \{\alpha_0,\ldots,\alpha_{n}\}$
of $V$ that belongs to the family $\s{S}$. 
Thus, the set $\sigma$ corresponds to $(n+1)!$ linearly ordered sequences
$\mapdef{s}{\{1, 2, \ldots, n+1\}}{\sigma}$, where each $s$
is a bijection.
We define an equivalence relation on these sequences
by saying that two sequences 
$\mapdef{s_1}{\{1, 2, \ldots, n+1\}}{\sigma}$ and
$\mapdef{s_2}{\{1, 2, \ldots, n+1\}}{\sigma}$ are equivalent iff
$\pi = s_2^{-1}\circ s_1$ is a permutation  of even signature
($\pi$ is the product of an even number of transpositions)

\medskip
The two equivalence classes associated with
$\sigma$ are called {\it oriented simplices\/}, and if
$\sigma = \{\alpha_0,\ldots,\alpha_{n}\}$, we denote
the equivalence class of $s$ as
$[s(1),\ldots,s(n+1)]$,
where $s$ is one of the sequences
$\mapdef{s}{\{1, 2, \ldots, n+1\}}{\sigma}$.
We also say that the two classes associated with $\sigma$
are the {\it orientations of $\sigma$\/}.
Two oriented simplices $\sigma_1$ and $\sigma_2$
are said to have {\it opposite orientation\/} if
they are the two classes associated with some simplex $\sigma$.
Given an oriented simplex $\sigma$, we denote the
oriented simplex having the opposite orientation as $-\sigma$,
with the convention that $-(-\sigma) = \sigma$.

\medskip
For example, if $\sigma = \{a_1, a_2, a_3\}$ is
a $3$-simplex (a triangle), there are six ordered sequences,
the sequences $\lag a_3, a_2, a_1\rag$, $\lag a_2, a_1, a_3\rag$, 
and $\lag a_1, a_3, a_2\rag$, are equivalent, and
the sequences $\lag a_1, a_2, a_3\rag$, $\lag a_2, a_3, a_1\rag$, 
and  $\lag a_3, a_1, a_2\rag$,  are also equivalent.
Thus, we have the two oriented simplices, $[a_1, a_2, a_3]$ and
$[a_3, a_2, a_1]$.
We now define $p$-chains.

\begin{defin}
\label{pchain}
{\em
Given a complex $K = (V, \s{S})$, a {\it $p$-chain\/} on $K$
is a function $c$ from the set of oriented $p$-simplices to $\integs$,
such that,
\begin{enumerate}
\item[(1)] 
$c(-\sigma) = -c(\sigma)$, iff $\sigma$ and $-\sigma$ 
have opposite orientation;
\item[(2)] 
$c(\sigma) = 0$, for all but finitely many simplices $\sigma$.
\end{enumerate}
We define addition of $p$-chains pointwise, i.e., 
$c_1 + c_2$ is the $p$-chain such that
$(c_1 + c_2)(\sigma) = c_1(\sigma) + c_2(\sigma)$, for every
oriented $p$-simplex $\sigma$.
The group of $p$-chains is denoted as $C_p(K)$. 
If $p < 0$ or $p > \dimm{(K)}$,
we set $C_p(K) = \{0\}$. 
}
\end{defin}

\medskip
To every oriented $p$-simplex $\sigma$ is associated 
an {\it elementary $p$-chain\/} $c$, defined such that,

\medskip
$c(\sigma) = 1$,

\medskip
$c(-\sigma) = -1$, where $-\sigma$ is the opposite orientation of $\sigma$,
and

\medskip
$c(\sigma') = 0$, for all other oriented simplices $\sigma'$.

\medskip
We will often denote the elementary $p$-chain associated with
the oriented $p$-simplex $\sigma$ also as $\sigma$.

\medskip
The following proposition is obvious, and simply confirms the fact
that $C_p(K)$ is indeed a free abelian group.

\begin{prop}
\label{freehom}
For every complex $K = (V, \s{S})$, for every $p$,
the group $C_p(K)$ is a free abelian group. For every
choice of an orientation for every $p$-simplex,
the corresponding elementary chains form a basis for $C_p(K)$.
\end{prop}

\medskip
The only point worth elaborating is that except for
$C_0(K)$, where no choice is involved, there is no canonical basis
for $C_p(K)$ for $p \geq 1$, since different choices for the orientations
of the simplices yield different bases.

\medskip
If there are $m_p$ $p$-simplices in $K$, the above proposition shows
that $C_p(K) = \integs^{m_{p}}$.

\medskip
As an immediate consequence of Proposition \ref{freehom}, for any 
abelian group $G$ and any function $f$ mapping the oriented
$p$-simplices of a complex $K$ to $G$, and such that
$f(-\sigma) = -f(\sigma)$ for every oriented $p$-simplex $\sigma$,
there is a unique homomorphism
$\mapdef{\widehat{f}}{C_p(K)}{G}$ extending $f$.

\medskip
We now define the boundary maps 
$\mapdef{\partial_{p}}{C_{p}(K)}{C_{p-1}(K)}$.

\begin{defin}
\label{boundmap1}
{\em
Given a complex $K = (V, \s{S})$, for every
oriented $p$-simplex 
\[\sigma = [\alpha_0,\ldots,\alpha_{p}],\]
we define the {\it bounday $\partial_{p}\sigma$ of $\sigma$\/} as
\[\partial_{p}\sigma = \sum_{i = 0}^{p}
(-1)^{i}[\alpha_0,\ldots,\widehat{\alpha_i},\ldots,\alpha_{p}],\]
where $[\alpha_0,\ldots,\widehat{\alpha_i},\ldots,\alpha_{p}]$
denotes the oriented $p-1$-simplex obtained by deleting vertex $\alpha_i$.
The {\it boundary map
$\mapdef{\partial_{p}}{C_{p}(K)}{C_{p-1}(K)}$\/} is the unique
homomorphism extending $\partial_p$ on oriented $p$-simplices.
For $p\leq 0$, $\partial_p$ is the null homomorphism.
}
\end{defin}

\medskip
One must verify that $\partial_{p}(-\sigma) = -\partial_{p}\sigma$,
but this is immediate. If $\sigma = [\alpha_0, \alpha_1]$,
then 
\[\partial_1\sigma = \alpha_1 - \alpha_0.\]
If  $\sigma = [\alpha_0, \alpha_1, \alpha_2]$, then
\[\partial_2\sigma = [\alpha_1, \alpha_2] - [\alpha_0, \alpha_2]
+ [\alpha_0, \alpha_1]
= [\alpha_1, \alpha_2] + [\alpha_2, \alpha_0]
+ [\alpha_0, \alpha_1].
\]
If  $\sigma = [\alpha_0, \alpha_1, \alpha_2, \alpha_3]$, then
\[\partial_3\sigma = [\alpha_1, \alpha_2, \alpha_3] 
- [\alpha_0, \alpha_2, \alpha_3]
+ [\alpha_0, \alpha_1, \alpha_3] - [\alpha_0, \alpha_1, \alpha_2].\]

\medskip
We have the following fundamental property.

\begin{prop}
\label{boundlem1}
For every complex $K = (V, \s{S})$, for every $p$, we have
$\partial_{p-1}\circ\partial_{p} = 0$.
\end{prop}

\proof For any oriented $p$-simplex
$\sigma =  [\alpha_0,\ldots,\alpha_{p}]$, we have
\[\eqaligneno{
\partial_{p-1}\circ\partial_{p}\sigma &= \sum_{i = 0}^{p}
(-1)^{i}\partial_{p-1}[\alpha_0,\ldots,\widehat{\alpha_i},
\ldots,\alpha_{p}],\cr
&= \quad\sum_{i = 0}^{p}\sum_{j = 0}^{i - 1} (-1)^{i}(-1)^{j}
[\alpha_0,\ldots,\widehat{\alpha_j},\ldots,\widehat{\alpha_i},
\ldots,\alpha_{p}]\cr
&\quad + 
\sum_{i = 0}^{p}\sum_{j = i + 1}^{p} (-1)^{i}(-1)^{j-1}
[\alpha_0,\ldots,\widehat{\alpha_i},\ldots,\widehat{\alpha_j},
\ldots,\alpha_{p}]\cr
&= 0.\cr
}\]
The rest of the proof follows from the fact that
$\mapdef{\partial_{p}}{C_{p}(K)}{C_{p-1}(K)}$ is the unique
homomorphism extending $\partial_p$ on oriented $p$-simplices.
$\square$

\medskip
In view of Proposition \ref{boundlem1}, the image $\partial_{p+1}(C_{p+1}(K))$
of $\mapdef{\partial_{p+1}}{C_{p+1}(K)}{C_{p}(K)}$ is a subgroup of the
kernel $\partial_p^{-1}(0)$ of $\mapdef{\partial_{p}}{C_{p}(K)}{C_{p-1}(K)}$.
This motivates the following definition.

\begin{defin}
\label{homgrps}
{\em
Given a complex $K = (V, \s{S})$,
the kernel $\partial_p^{-1}(0)$ of the homomorphism 
$\mapdef{\partial_{p}}{C_{p}(K)}{C_{p-1}(K)}$ is denoted
as $Z_p(K)$, and the elements of $Z_p(K)$ are called
{\it $p$-cycles\/}. The image $\partial_{p+1}(C_{p+1})$
of the homomorphism 
$\mapdef{\partial_{p+1}}{C_{p+1}(K)}{C_{p}(K)}$ is denoted
as $B_p(K)$, and the elements of $B_p(K)$ are called
{\it $p$-boundaries\/}. 
The {\it $p$-th homology group $H_p(K)$\/} is the quotient group
\[H_p(K) = Z_p(K)/B_p(K).\]
Two $p$-chains $c, c'$ are said to be {\it homologous\/} if
there is some $(p+1)$-chain $d$ such that 
$c = c' + \partial_{p+1}d$. 
}
\end{defin}

\medskip
We will often omit the subscript $p$ in $\partial_p$.

\medskip
At this stage, we could determine the homology groups of the
finite (two-dimensional) polyhedras. 
However, we are really interested in the homology groups
of geometric realizations of complexes, in particular, compact surfaces,
and so far, we have not defined homology groups for topological spaces.

\medskip
It is possible to define homology groups for arbitrary topological
spaces, using what is called {\it singular homology\/}.
Then, it can be shown, although this requires some hard work, 
that the homology groups of a space $X$ which is the
geometric realization of some complex $K$ are independent 
of the complex $K$ such that $X = K_g$, and equal to the homology
groups of any such complex. 

\medskip
The idea behind singular homology is to define a more general notion
of an $n$-simplex associated with a topological space $X$,
and it is natural to consider continuous maps from some
standard simplices to $X$.
Recall that given any set $I$, we defined the real vector space
$\Ispac{\reals}{I}$ freely generated by $I$ 
(just before Definition \ref{geomreal}).
In particular, for $I = \natnums$ (the natural numbers), 
we obtain an infinite dimensional vector space $\Ispac{\reals}{\natnums}$,
whose elements are the countably infinite sequences
$(\lambda_i)_{i\in\natnums}$ of reals, with $\lambda_i = 0$ for all but
finitely many $i\in \natnums$.
For any $p\in \natnums$, we let $e_i\in \Ispac{\reals}{\natnums}$ be the
sequence such that $e_i(i) = 1$ and $e_i(j) = 0$ for all $j\not= i$,
and we let $\Delta_p$ be the $p$-simplex spanned by
$(e_0,\ldots, e_p)$, that is, the subset of $\Ispac{\reals}{\natnums}$
consisting of all points of the form
\[\sum_{i = 0}^{p}\lambda_ie_i,\quad\hbox{with}\quad 
\sum_{i = 0}^{p}\lambda_i = 1,\ \hbox{and}\ \lambda_i\geq 0.\]
We call $\Delta_p$ the {\it standard $p$-simplex\/}.
Note that $\Delta_{p-1}$ is a face of $\Delta_p$.

\begin{defin}
\label{singsimp}
{\em
Given a topological space $X$, a {\it singular $p$-simplex\/} is any
continuous map $\mapdef{T}{\Delta_p}{X}$. The free abelian group
generated by the singular $p$-simplices is called the 
{\it $p$-th singular chain group\/}, and is denoted as $S_p(X)$. 
}
\end{defin}

\medskip
Given any $p+1$ points $a_0,\dots, a_p$ in $\Ispac{\reals}{\natnums}$,
there is a unique affine map
$\mapdef{f}{\Delta_p}{\Ispac{\reals}{\natnums}}$, such that
$f(e_i) = a_i$, for all $i$, $0\leq i\leq p$, namely the map
such that
\[f(\sum_{i = 0}^{p}\lambda_ie_i) = 
\sum_{i = 0}^{p}\lambda_ia_i,\]
for all $\lambda_i$ such that
$\sum_{i = 0}^{p}\lambda_i = 1$, and $\lambda_i\geq 0$.
This map is called the {\it affine singular simplex\/} determined
by $a_0,\ldots, a_p$, and it is denoted as
$l(a_0,\ldots, a_p)$. In particular, the map
\[l(e_0,\ldots, \widehat{e_i},\ldots,e_p),\]
where the hat over $e_i$ means that $e_i$ is omited, is a map
from $\Delta_{p-1}$ onto a face of $\Delta_p$.
We can consider it as a map from $\Delta_{p-1}$ to $\Delta_p$
(although it is defined as a map from $\Delta_{p-1}$ to 
$\Ispac{\reals}{\natnums}$), and call it the $i$-th face
of $\Delta_p$.

\medskip
Then, if $\mapdef{T}{\Delta_p}{X}$ is a singular $p$-simplex,
we can form the map
\[\mapdef{T\circ l(e_0,\ldots, \widehat{e_i},\ldots,e_p)}{\Delta_{p-1}}{X},\]
which is a singular $p-1$-simplex, which we think of as the $i$-th face
of $T$.
Actually, for $p = 1$, a singular $p$-simplex 
$\mapdef{T}{\Delta_p}{X}$ can be viewed as curve on $X$,
and its faces are its two endpoints.
For $p = 2$,  a singular $p$-simplex 
$\mapdef{T}{\Delta_p}{X}$ can be viewed as triangular surface
patch on $X$, and its faces are its three boundary curves.
For $p = 3$,  a singular $p$-simplex 
$\mapdef{T}{\Delta_p}{X}$ can be viewed as tetrahedral ``volume patch'' 
on $X$, and its faces are its four boundary surface patches.
We can give similar higher-order descriptions when $p > 3$.

\medskip
We can now define the boundary maps
$\mapdef{\partial_{p}}{S_{p}(X)}{S_{p-1}(X)}$.

\begin{defin}
\label{boundmap2}
{\em
Given a topological space $X$, for every
singular  $p$-simplex $\mapdef{T}{\Delta_p}{X}$,
we define the {\it bounday $\partial_{p}T$ of $T$\/} as
\[\partial_{p}T = \sum_{i = 0}^{p}
(-1)^{i}\, T\circ l(e_0,\ldots, \widehat{e_i},\ldots,e_p).\]
The {\it boundary map
$\mapdef{\partial_{p}}{S_{p}(X)}{S_{p-1}(X)}$\/} is the unique
homomorphism extending $\partial_p$ on singular  $p$-simplices.
For $p\leq 0$, $\partial_p$ is the null homomorphism.
Given a continuous map $\mapdef{f}{X}{Y}$ between
two topological spaces $X$ and $Y$, the
homomorphism $\mapdef{f_{\sharp, p}}{S_p(X)}{S_p(Y)}$
is defined such that
\[f_{\sharp, p}(T)  = f\circ T,\]
for every singular $p$-simplex $\mapdef{T}{\Delta_p}{X}$.
}
\end{defin}

\medskip
The next easy proposition gives the main properties of
$\partial$.

\begin{prop}
\label{boundlem2}
For every continuous map $\mapdef{f}{X}{Y}$ between
two topological spaces $X$ and $Y$, the maps $f_{\sharp, p}$ and
$\partial_p$ commute  for every $p$, i.e.,
\[\partial_{p}\circ f_{\sharp, p} = f_{\sharp, p-1}\circ \partial_{p}.\]
We also have
$\partial_{p-1}\circ\partial_{p} = 0$.
\end{prop}

\proof For any singular $p$-simplex $\mapdef{T}{\Delta_p}{X}$,
we have
\[\partial_p f_{\sharp, p}(T) =
\sum_{i = 0}^{p} (-1)^i\, (f\circ T)\circ 
l(e_0,\ldots, \widehat{e_i},\ldots,e_p),\]
and
\[f_{\sharp, p-1}(\partial_{p} T) =
\sum_{i = 0}^{p} (-1)^i\, f\circ (T\circ 
l(e_0,\ldots, \widehat{e_i},\ldots,e_p)),\]
and the equality follows by associativity of composition.
We also have
\[\eqaligneno{
\partial_{p} l(a_0,\ldots, a_p) 
&= \sum_{i = 0}^{p} (-1)^{i}\,
l(a_0,\ldots, a_p) \circ l(e_0,\ldots, \widehat{e_i},\ldots,e_p)\cr
&= \sum_{i = 0}^{p} (-1)^{i}\,
l(a_0,\ldots, \widehat{a_i},\ldots,a_p),
}\]
since the composition of affine maps is affine.
Then, we can compute 
$\partial_{p-1}\partial_{p} l(a_0,\ldots, a_p)$
as we did in Proposition \ref{boundlem1}, and the proof is similar,
except that we have to insert an $l$ at appropriate places. 
The rest of the proof follows from the fact that
\[\partial_{p-1}\partial_{p} T = \partial_{p-1}\partial_{p}
(T_{\sharp}(l(e_0,\ldots, e_p))),\]
since $l(e_0,\ldots, e_p)$ is simply the inclusion of
$\Delta_p$ in $\Ispac{\reals}{\natnums}$,
and that $\partial$ commutes with $T_{\sharp}$.
$\square$

\medskip
In view of Proposition \ref{boundlem2}, the image $\partial_{p+1}(S_{p+1}(X))$
of $\mapdef{\partial_{p+1}}{S_{p+1}(X)}{S_{p}(X)}$ is a subgroup of the
kernel $\partial_p^{-1}(0)$ of $\mapdef{\partial_{p}}{S_{p}(X)}{S_{p-1}(X)}$.
This motivates the following definition.

\begin{defin}
\label{homgrps2}
{\em
Given a topological space $X$,
the kernel $\partial_p^{-1}(0)$ of the homomorphism 
$\mapdef{\partial_{p}}{S_{p}(X)}{S_{p-1}(X)}$ is denoted
as $Z_p(X)$, and the elements of $Z_p(X)$ are called
{\it singular $p$-cycles\/}. The image $\partial_{p+1}(S_{p+1})$
of the homomorphism 
$\mapdef{\partial_{p+1}}{S_{p+1}(X)}{S_{p}(X)}$ is denoted
as $B_p(X)$, and the elements of $B_p(X)$ are called
{\it singular $p$-boundaries\/}. 
The {\it $p$-th singular homology group $H_p(X)$\/} is the quotient group
\[H_p(X) = Z_p(X)/B_p(X).\]
}
\end{defin}

\medskip
If $\mapdef{f}{X}{Y}$ is a continuous map, the fact that
\[\partial_{p}\circ f_{\sharp, p} = f_{\sharp, p-1}\circ \partial_{p}\]
allows us to define homomorphisms
$\mapdef{f_{*, p}}{H_{p}(X)}{H_{p}(Y)}$, and it it easily verified that
\[(g\circ f)_{*, p} = g_{*, p}\circ f_{*, p},\]
and that $\mapdef{Id_{*, p}}{H_{p}(X)}{H_{p}(Y)}$ is the identity
homomorphism, when  $\mapdef{Id}{X}{Y}$ is the identity.
As a corollary, if $\mapdef{f}{X}{Y}$ is a homeomorphism,
then each $\mapdef{f_{*, p}}{H_{p}(X)}{H_{p}(Y)}$ is a group
isomorphism. This gives us a way of showing that two spaces
are not homeomorphic, by showing that some homology groups
$H_p(X)$ and $H_p(Y)$ are not  isomorphic.

\medskip
It is fairly easy to show that $H_0(X)$ is a free abelian group,
and that if the path components of $X$ are the family
$(X_i)_{i\in I}$, then $H_0(X)$ is isomorphic to the direct
sum $\bigoplus_{i\in I} \integs$. In particular, if
$X$ is arcwise connected, $H_0(X) = \integs$.

\medskip
The following important theorem shows  the relationship
between simplicial homology and singular homology.
The proof is fairly involved, and can be found in Munkres \cite{Munkresalg},
or Rotman \cite{Rotman}.

\begin{thm}
\label{simpsing}
Given any polytope $X$, if $X = K_g = K'_g$
is the geometric realization of any two complexes $K$ and $K'$, then 
\[H_p(X) = H_p(K) = H_p(K'),\]
for all $p\geq 0$.
\end{thm}

\medskip
Theorem \ref{simpsing} implies that $H_p(X)$ is finitely
generated for all $p\geq 0$. It is immediate that if $K$ has
dimension $m$, then $H_p(X) = 0$ for $p > m$,
and it can be shown that $H_m(X)$ is a free abelian group.

\medskip
A fundamental invariant of finite complexes is the
Euler-Poincar\'e characteristic.

\begin{defin}
\label{Eulerchar}
{\em
Given a finite complex $K = (V, \s{S})$ of dimension $m$,
letting $m_p$ be the number of $p$-simplices in $K$,
we define the {\it Euler-Poincar\'e characteristic $\chi(K)$ of $K$\/} as
\[\chi(K) = \sum_{p = 0}^{m} (-1)^{p}\> m_p.\]
}
\end{defin}

\medskip
The following remarkable theorem holds.

\begin{thm}
\label{Eulercharthm}
Given a finite complex $K = (V, \s{S})$ of dimension $m$,
we have
\[\chi(K) = \sum_{p = 0}^{m} (-1)^{p}\> r(H_p(K)),\]
the alternating sum of the Betti numbers (the ranks) 
of the homology groups of $K$. 
\end{thm}

\proof We know that $C_p(K)$ is a free group of rank
$m_p$. Since $H_p(K) = Z_p(K)/B_p(K)$, by Proposition \ref{finab3}, we have
\[r(H_p(K)) = r(Z_p(K)) - r(B_p(K)).\]
Since we have a short exact sequence
\[\shortexact{Z_p(K)}{}{C_p(K)}{\partial_p}{B_{p-1}(K)},\]
again, by Proposition \ref{finab3}, we have
\[r(C_p(K)) = m_p = r(Z_p(K)) + r(B_{p-1}(K)).\]
Also, note that $B_m(K) = 0$, and $B_{-1}(K) = 0$.
Then, we have
\[\eqaligneno{
\chi(K) &= \sum_{p = 0}^{m} (-1)^{p}\> m_p\cr
&= \sum_{p = 0}^{m} (-1)^{p}\> (r(Z_p(K)) + r(B_{p-1}(K)))\cr
&= \sum_{p = 0}^{m} (-1)^{p}\> r(Z_p(K)) 
+ \sum_{p = 0}^{m} (-1)^{p}\> r(B_{p-1}(K)).\cr
}\]
Using the fact that $B_m(K) = 0$, and $B_{-1}(K) = 0$, we get
\[\eqaligneno{
\chi(K) &= \sum_{p = 0}^{m} (-1)^{p}\> r(Z_p(K)) 
+ \sum_{p = 0}^{m} (-1)^{p+1}\> r(B_{p}(K))\cr
&= \sum_{p = 0}^{m} (-1)^{p}\> (r(Z_p(K)) - r(B_{p}(K)))\cr
&= \sum_{p = 0}^{m} (-1)^{p}\> r(H_p(K)).\cr
}\]
$\square$

\medskip
A striking corollary of Theorem \ref{Eulercharthm}
(together with Theorem \ref{simpsing}), is that
the Euler-Poincar\'e characteristic $\chi(K)$ of a complex
of finite dimension $m$ only depends on the
geometric realization $K_g$ of $K$, since
it only depends on the homology groups $H_p(K) = H_p(K_g)$ of 
the polytope $K_g$.
Thus, the Euler-Poincar\'e characteristic is an invariant of
all the finite complexes corresponding to the same
polytope $X = K_g$, and we can say that it is the
Euler-Poincar\'e characteristic of the polytope $X = K_g$, and
denote it as $\chi(X)$. In particular, this is true of
surfaces that admit a triangulation, and as we shall see
shortly,  the Euler-Poincar\'e characteristic in one of the major
ingredients in the classification of the compact surfaces.
In this case, $\chi(K) = m_0 - m_1 + m_2$, where
$m_0$ is the number of vertices, $m_1$ the number of edges,
and $m_2$ the number of triangles, in $K$. We warn the reader
that Ahlfors and Sario have flipped the signs, and
define the Euler-Poincar\'e characteristic as $-m_0 + m_1 - m_2$.
\medskip

\medskip
Going back to the triangulations of the sphere, the torus,
the projective space, and the Klein bottle, it is easy to see
that their Euler-Poincar\'e characteristic is 
$2$ (sphere), $0$ (torus), $1$ (projective space), and
$0$ (Klein bottle).

\medskip
At this point, we are ready to compute the homology groups of
finite (two-dimensional) polyhedras.

\section{Homology Groups of the Finite  Polyhedras}
\label{grpsurf}
Since a polyhedron is the geometric realization of 
a triangulated $2$-complex, it is possible
to determine the homology groups of the (finite) polyhedras.
We say that a triangulated $2$-complex $K$ is orientable if
its geometric realization $K_g$ is orientable.
We will consider the finite, bordered, orientable, and nonorientable,
triangulated $2$-complexes.
First, note that $C_p(K)$ is the trivial group for $p<0$ and $p>2$,
and thus, we just have to consider the cases where $p = 0, 1, 2$.
We will use the notation $c\sim c'$, to denote that two
$p$-chains are homologous, which means that
$c = c' + \partial_{p+1} d$, for some $(p+1)$-chain $d$.

\medskip
The first proposition is very easy, and is just a special case of
the fact that $H_0(X) = \integs$ for an arcwise connected space $X$.

\begin{prop}
\label{homol1}
For every triangulated $2$-complex (finite or not) $K$, 
we have $H_0(K) = \integs$.
\end{prop}

\proof When $p = 0$, we have $Z_0(K) = C_0(K)$,
and thus, $H_0(K) = C_0(K)/B_0(K)$.
Thus, we have to figure out what the $0$-boundaries are.
If $ c = \sum x_i\partial a_i$ is a $0$-boundary,
each $a_i$ is an oriented edge $[\alpha_i, \beta_i]$, and we have
\[ c = \sum x_i\partial a_i = \sum x_i\beta_i - \sum x_i\alpha_i,\]
which shows that the sum of all the coefficients of the vertices is $0$.
Thus, it is impossible for a $0$-chain of the form
$x\alpha$, where $x\not=0$, to be homologous to $0$.
On the other hand, we claim that $\alpha \sim \beta$ for any
two vertices $\alpha, \beta$. Indeed, since we assumed that $K$ is connected,
there is a path from $\alpha$ to $\beta$ consisting of edges
\[[\alpha, \alpha_1],\ldots,[\alpha_{n}, \beta],\]
and the $1$-chain
\[ c = [\alpha, \alpha_1]+\ldots+[\alpha_{n}, \beta]\]
has boundary 
\[\partial c = \beta - \alpha,\]
which shows that 
$\alpha \sim \beta$.
But then, $H_0(K)$ is the infinite cyclic group generated
by any vertex.
$\square$

\medskip
Next, we determine the groups $H_2(K)$.

\begin{prop}
\label{homol2}
For every triangulated $2$-complex (finite or not) $K$, 
either $H_2(K) = \integs$
or $H_2(K) = 0$.  Furthermore, $H_2(K) = \integs$
iff $K$ is finite, has no border and is orientable, else $H_2(K) = 0$.
\end{prop}

\proof When $p = 2$, we have $B_2(K) = 0$, and $H_2(K) = Z_{2}(K)$.
Thus, we have to figure out what the $2$-cycles are.
Consider a $2$-chain $ c = \sum x_i A_i$, where each $A_i$
is an  oriented triangle $[\alpha_0, \alpha_1, \alpha_2]$, and assume that
$ c$ is a cycle, which means that
\[\partial c = \sum x_i\partial A_i = 0.\]
Whenever $A_i$ and $A_j$ have an edge $a$ in common, 
the contribution of $a$ to $\partial c$ is either
$x_i a + x_j a$, or $x_i a - x_j a$, or $-x_i a + x_j a$, or
$-x_i a  -x_j a$,   which implies that $x_i = \epsilon x_j$,
with $\epsilon = \pm 1$. Consequently, if $A_i$ and $A_j$
are joined by a path of pairwise adjacent triangles, $A_k$,
all in $c$, then $|x_i| = |x_j|$. 
However, Proposition \ref{trianglem2} and Proposition \ref{trianglem3} imply that
any two triangles $A_i$ and $A_j$ in $K$ are connected by
a sequence of pairwise adjacent triangles. If some triangle 
in the path does not belong to $c$, then there are two adjacent 
triangles in the path, $A_h$ and $A_k$, with $A_h$ in $c$ and $A_k$ not in $c$
such that all the triangles in the path from $A_i$ to $A_h$ belong to $c$.
But then, $A_h$  has an edge not adjacent to any other triangle in
$c$, so $x_h = 0$ and thus, $x_i = 0$. The same reasoning applied to $A_j$
shows that $x_j = 0$. If all triangles in the path from $A_i$ to $A_j$
belong to $c$, then we already know that $|x_i| = |x_j|$.
Therefore, all  $x_i$'s have the same absolute value.
If $K$ is infinite, there must be some $A_i$ in the finite sum
which is adjacent to some triangle $A_j$ not in the finite sum,
and the contribution of the edge common to $A_i$ and $A_j$ to
$\partial c$ must be zero, which implies that $x_i = 0$
for all $i$. Similarly, the coefficient of every triangle with an edge
in the border must be zero.
Thus, in these cases, $c \sim 0$, and 
$H_2(K) = 0$.

\medskip
Let us now assume that $K$ is a finite 
triangulated $2$-complex without a border.
The above reasoning showed that any nonzero $2$-cycle, $c$, can be written as
\[ c = \sum \epsilon_i x A_i,\]
where $x = |x_i| > 0$ for all $i$, and $\epsilon_i = \pm 1$.
Since $\partial c  = 0$, $\sum\epsilon_i A_i$ is also
a $2$-cycle. For any other nonzero $2$-cycle, $\sum y_i A_i$,
we can subtract $\epsilon_1 y_1(\sum \epsilon_i A_i)$
from  $\sum y_i A_i$, and we get the cycle
\[\sum_{i\not= 1} (y_i - \epsilon_1\epsilon_i y_1) A_i,\]
in which $A_1$ has coefficient $0$. But then, since all the
coefficients have the same absolute value, we must have
$y_i = \epsilon_1\epsilon_i y_1$ for all $i\not= 1$,
and thus, 
\[\sum y_i A_i = \epsilon_1 y_1 (\sum \epsilon_i A_i).\]
This shows that either $H_2(K) = 0$, or $H_2(K) = \integs$.

\medskip
It remains to prove that $K$ is orientable iff $H_2(K) = \integs$.
The idea is that in this case, we can choose an orientation
such that $\sum A_i$ is a $2$-cycle. The proof is not really
difficult, but a little involved, and  the reader
is referred to Ahlfors and Sario \cite{Ahlfors} for details.
$\square$

\medskip
Finally, we need to determine $H_1(K)$. We will only do so
for finite triangulated $2$-complexes, and refer the reader 
to Ahlfors and Sario \cite{Ahlfors} for the infinite case.

\begin{prop}
\label{homol3}
For every finite triangulated $2$-complex $K$, 
either $H_1(K) = \integs^{m_1}$,
or $H_1(K) = \integs^{m_1}\oplus \integs/2\integs$, 
the second case occurring iff $K$ has no border and is
nonorientable.
\end{prop}

\proof 
The first step is to determine the torsion subgroup of $H_1(K)$.
Let $c$ be a $1$-cycle, and assume that $mc \sim 0$ for some
$m > 0$, i.e., there is some $2$-chain $\sum x_iA_i$ such that
$mc = \sum x_i\partial A_i$.
If $A_i$ and $A_j$ have a common edge $a$, the contribution of
$a$ to $\sum x_i\partial A_i$ is either
$x_i a + x_j a$, or $x_i a - x_j a$, or $-x_i a + x_j a$, or
$-x_i a  -x_j a$,   which implies that either 
$x_i \equiv  x_j (\hbox{mod}\> m)$,
or $x_i \equiv  -x_j (\hbox{mod}\> m)$.
Because of the connectedness of $K$, the above actually
holds for all $i, j$. If $K$ is bordered, there is some
$A_i$ which contains a border edge not adjacent to any other
triangle, and thus $x_i$ must be divisible by $m$, which
implies that every $x_i$ is divisible by $m$.
Thus, $c \sim 0$. Note that a similar reasoning applies when
$K$ is infinite, but we are not considering this case.
If $K$ has no border and is orientable, by a previous remark, we can
assume that $\sum A_i$ is a cycle. Then,
$\sum \partial A_i = 0$, and we can write
\[mc = \sum(x_i - x_1)\partial A_i.\]
Due to the connectness of $K$, the above argument shows
that every $x_i - x_1$ is divisible by $m$, which
shows that $c\sim 0$. Thus, the torsion group is $0$.

\medskip
Let us now assume that $K$ has no border and is nonorientable.
Then, by a previous remark, there are no $2$-cycles except $0$.
Thus, the coefficients in $\sum \partial A_i$ must
be either $0$ or $\pm 2$. Let $\sum\partial A_i = 2z$.
Then, $2z \sim 0$, but $z$ is not homologous to $0$,
since from $z = \sum x_i\partial A_i$, we would get
$\sum(2x_i -1)\partial A_i \sim 0$, contrary to the fact
that there are no $2$-cycles except $0$. Thus, $z$ is of order $2$.

\medskip
Consider again  $mc = \sum x_i\partial A_i$.
Since $x_i \equiv  x_j (\hbox{mod}\> m)$,
or $x_i \equiv  -x_j (\hbox{mod}\> m)$, for all $i, j$, we can write
\[mc = x_1\sum\epsilon_i\partial A_i + m\sum t_i\partial A_i,\]
with $\epsilon_i = \pm 1$, and at least some coefficient
of $\sum\epsilon_i\partial A_i$ is $\pm 2$, since otherwise
$\sum\epsilon_i A_i$ would be a nonnull $2$-cycle.
But then, $2x_1$ is divisible by $m$, and this implies that
$2c\sim 0$. If $2c = \sum u_i\partial A_i$,
the $u_i$ are either all odd or all even. If they are all even,
we get $c\sim 0$, and if they are all odd, we get $c\sim z$.
Hence, $z$ is the only element of finite order,
and the torsion group if $\integs/2\integs$.

\medskip
Finally, having determined the torsion group of $H_1(K)$,
by the corollary of Proposition \ref{finab2}, 
we know that $H_{1}(K) = \integs^{m_{1}}\oplus T$,
where $m_1$ is the rank of $H_1(K)$, and the proposition follows.
$\square$

\medskip
Recalling Proposition \ref{Eulercharthm},
the Euler-Poincar\'e characteristic $\chi(K)$ is given by
\[\chi(K) = r(H_{0}(K)) - r(H_{1}(K)) + r(H_{2}(K)),\]
and we have determined that $r(H_{0}(K)) = 1$ and 
either $r(H_{2}(K)) = 0$
when $K$  has a border or has no border and is nonorientable, 
or  $r(H_{2}(K)) = 1$
when $K$  has no border and is orientable.

\medskip
Thus, the rank $m_1$ of $H_{1}(K)$ is either 
\[m_1 = 2 - \chi(K)\]
if $K$ has no border and is orientable, and 
\[m_1 = 1 - \chi(K)\]
otherwise. This implies that $\chi(K) \leq 2$.

\medskip
We will now prove the classification theorem for compact
(two-dimensional) polyhedras.

\chapter{The Classification Theorem for Compact Surfaces}
\label{chap6}
\section{Cell Complexes}
\label{cellcplx}
It is remarkable that  the compact (two-dimensional)
polyhedras can be characterized up to homeomorphism. 
This situation is exceptional, as such a result is known to be 
essentially impossible for
compact $m$-manifolds for $m\geq 4$, and still open for
compact $3$-manifolds. In fact, it is possible to characterize
the compact (two-dimensional)
polyhedras in terms of a simple extension of
the notion of a complex, called cell complex  by Ahlfors and Sario.
What happens is that it is possible to define an equivalence relation
on cell complexes, and it can be shown that every
cell complex is equivalent to some specific normal form.
Furthermore, every cell complex has  a geometric
realization which is a surface,
and equivalent cell complexes have homeomorphic
geometric realizations. Also, every cell complex is equivalent
to a triangulated $2$-complex.
Finally, we can show that
the geometric realizations of distinct normal forms are not
homeomorphic.

\medskip
The classification theorem for compact surfaces is presented
(in slightly different ways) in Massey \cite{Massey87}
Amstrong \cite{Amstrong}, and Kinsey \cite{Kinsey}.
In the above references, the presentation is sometimes quite informal.
The classification theorem is also presented in
Ahlfors and Sario \cite{Ahlfors}, and there,
the presentation is  formal and not always easy to follow. 
We tried to strike a middle ground in the degree of formality.
It should be noted that the combinatorial part of
the proof (Section \ref{cpxnormal}) is  heavily
inspired by the proof given in Seifert and Threlfall \cite{Seifert80}.
One should also take a look at Chapter 1 of 
Thurston \cite{Thurston97}, especially Problem 1.3.12.
Thurston's book is also highly recommended
as a wonderful and insighful introduction to the topology
and geometry of three-dimensional manifolds, but that's
another story.

\medskip
The first step is to define cell complexes. The intuitive idea
is to generalize a little bit the notion of a triangulation,
and consider objects made of oriented faces, each face having some 
boundary. A boundary is a cyclically ordered list of oriented edges.
We can think of each face as a circular
closed disk, and of the edges in a boundary as circular arcs on
the boundaries of these disks.
A cell complex  represents the surface obtained
by identifying identical boundary edges.

\medskip
Technically, in order to deal with the notion of orientation,
given any set $X$, it is convenient to introduce the set
$X^{-1} = \{x^{-1}\ |\ x\in X\}$ of formal inverses of elements in $X$.
We will say that the elements of $X\cup X^{-1}$ are {\it oriented\/}.
It is also convenient to assume that $(x^{-1})^{-1} = x$, for
every $x\in X$.
It turns out that cell complexes can be defined
using only faces and boundaries, and that the notion of a vertex can be
defined from the way edges occur in boundaries. 
This way of dealing with vertices is a bit counterintuitive,
but we haven't found a better way to present cell complexes.
We now give precise definitions.

\begin{defin}
\label{cellcomp}
{\em
A {\it cell complex K\/} consists of a triple
$K = (F, E, B)$, where $F$ is a finite nonempty set of {\it faces\/},
$E$ is a finite set of {\it edges\/}, and
$\mapdef{B}{(F\cup F^{-1})}{(E\cup E^{-1})^*}$
is the {\it boundary function\/}, which assigns
to each oriented face $A\in F\cup F^{-1}$ a cyclically ordered
sequence $a_1\ldots a_n$ of oriented edges in $E\cup E^{-1}$,
the {\it boundary of $A$\/},
in such a way that $B(A^{-1}) = a_n^{-1}\ldots a_1^{-1}$
(the reversal of the sequence $a_1^{-1}\ldots a_n^{-1}$).
By a cyclically ordered sequence, we mean that we do not
distinguish between the sequence $a_1\ldots a_n$ and 
any sequence obtained from it by a cyclic permutation.
In particular, the successor of $a_n$ is $a_1$.
Furthermore, the following conditions must hold:
\begin{enumerate}
\item[(1)] 
Every oriented edge $a\in E\cup E^{-1}$ occurs either once or twice
as an element of a boundary. In particular, this means that
if $a$ occurs twice in some boundary, then it does not occur in any
other boundary.
\item[(2)]
$K$ is connected. This means that $K$ is not the union of
two disjoint systems satisfying the conditions.
\end{enumerate}
}
\end{defin}

\medskip
It is possible that $F = \{A\}$ and $E = \emptyset$,
in which case $B(A) = B(A^{-1}) = \epsilon$, the empty sequence.

\medskip
For short, we will often say face and edge, rather than oriented face
or oriented edge.
As we said earlier, the notion of a vertex is defined in terms
of faces and boundaries. The intuition is that
a vertex is adjacent to pairs of incoming and outgoing edges.
Using inverses of edges, we can define a vertex as the 
sequence of incoming edges into that vertex. When the vertex is
not a boundary vertex, these edges form a cyclic sequence,
and when the vertex is a border vertex, such a sequence has
two endpoints with no successors.

\begin{defin}
\label{vertex}
{\em
Given a cell complex $K = (F, E, B)$, for any edge $a\in E\cup E^{-1}$,
a {\it successor\/} of $a$ is an edge $b$ such that $b$ is the successor
of $a$ in some boundary $B(A)$. If $a$ occurs in two places
in the set of boundaries, it has a {\it a pair of successors\/}
(possibly identical), and otherwise it has a {\it single successor\/}.
A cyclically ordered sequence $\alpha = (a_1,\ldots, a_n)$ is called
an {\it inner vertex\/} if every $a_i$ has $a_{i-1}^{-1}$
and $a_{i+1}^{-1}$ as pair of successors (note that
$a_1$ has $a_{n}^{-1}$ and $a_{2}^{-1}$ as pair of successors, and
$a_n$ has $a_{n-1}^{-1}$ and $a_{1}^{-1}$ as pair of successors).
A {\it border vertex\/} is a 
cyclically ordered sequence $\alpha = (a_1,\ldots, a_n)$ 
such that the above condition holds for all $i$, 
$2 \leq i \leq n-1$, while $a_1$ has $a_2^{-1}$ as only successor,
and $a_{n}$ has $a_{n-1}^{-1}$ as only successor.
An edge $a\in E\cup E^{-1}$ is a {\it border edge\/} if
it occurs once in a single boundary, and otherwise an
{\it inner edge\/}.
}
\end{defin}

\medskip
Given any edge $a\in E\cup E^{-1}$, we can determine a unique
vertex $\alpha$ as follows: the neighbors of $a$ in the vertex $\alpha$
are the inverses of its successor(s). Repeat this step in both directions
until either the cycle closes, or we hit sides with only one successor.
The vertex $\alpha$ in question is the list of the incoming edges into it.
For this reason, we say that {\it $a$ leads to $\alpha$\/}.
Note that when a vertex $\alpha = (a)$ contains
a single edge $a$, there must be an occurrence of the form
$aa^{-1}$ in some boundary. Also, note that if $(a, a^{-1})$ is a vertex,
then it is an inner vertex, and if $(a, b^{-1})$ is a vertex with
$a\not= b$, then it is a border vertex.

\medskip
Vertices can also characterized in another way which
will be useful later on. Intuitively, two edges $a$ and $b$
are equivalent iff they have the same terminal vertex.

\medskip
We define a relation $\lambda$ on
edges as follows: $a \lambda b$ iff $b^{-1}$ is the successor
of $a$ in some boundary. Note that this relation is symmetric.
Indeed, if $ab^{-1}$ appears in the boundary of some face $A$,
then $ba^{-1}$ appears in the boundary of $A^{-1}$.
Let $\Lambda$ be the reflexive and transitive closure of $\lambda$.
Since $\lambda$ is symmetric, $\Lambda$ is an equivalence
relation. We leave as an easy exercise to prove that
the equivalence class of an edge $a$ is the vertex $\alpha$
that $a$ leads to. Thus, vertices induce a partition
of $E\cup E^{-1}$. We say that an edge $a$ is an edge from
a vertex $\alpha$ to a vertex $\beta$ if $a^{-1}\in \alpha$
and $a\in \beta$. Then, by a familiar reasoning,
we can show that the fact that $K$ is connected implies that
there is a path between any two vertices.

\medskip
Figure \ref{cellfig1} shows a cell complex with border.
The cell complex has three faces with boundaries
$abc$, $bed^{-1}$, and $adf^{-1}$.
It has one inner vertex $b^{-1}ad^{-1}$ and three
border vertices $edf$, $c^{-1}be^{-1}$, and $ca^{-1}f^{-1}$.

\begin{figure}
  \begin{center}
    \begin{pspicture}(0,-4.5)(6,4.2)
    \psline[linewidth=1pt]{->}(0,0)(4,0)
    \psline[linewidth=1pt]{->}(0,0)(2,-4)
    \psline[linewidth=1pt]{->}(2,4)(0,0)
    \psline[linewidth=1pt]{->}(4,0)(2,4)
    \psline[linewidth=1pt]{->}(4,0)(6,4)
    \psline[linewidth=1pt]{->}(2,4)(6,4)
    \psline[linewidth=1pt]{->}(4,0)(2,-4)
    \uput[90](2,0){$a$}
    \uput[90](3,2){$b$}
    \uput[90](1,2){$c$}
    \uput[90](5,2){$d$}
    \uput[90](3,-2){$d$}
    \uput[90](1,-2){$f$}
    \uput[90](4,4){$e$}
    \end{pspicture}
  \end{center}
  \caption{A cell complex with border}
\label{cellfig1}
\end{figure}

If we fold the above cell complex by identifying the two edges labeled
$d$, we get a tetrahedron with one face omitted, the face opposite
the inner vertex, the endpoint of edge $a$.

\medskip
There is a natural way to view a triangulated complex
as a cell complex, and it is not hard to see that the following
conditions allow us to view a cell complex as a triangulated complex.
\begin{enumerate}
\item[(C1)] 
If $a, b$ are distinct edges leading to the same vertex,
then $a^{-1}$ and $b^{-1}$ lead to distinct vertices.
\item[(C2)] 
The boundary of every face is a triple $abc$.
\item[(C3)] 
Different faces have different boundaries.
\end{enumerate}

It is easy to see that $a$ and $a^{-1}$ cannot lead to the same vertex,
and that in a face $abc$, the edges $a, b, c$ are distinct.

\section{Normal Form for Cell Complexes}
\label{cpxnormal}
We now introduce a notion of elementary subdivision of
cell complexes which is crucial in obtaining the classification
theorem.

\begin{defin}
\label{subdiv}
{\em
Given any two cells complexes $K$ and $K'$,
we say that {\it $K'$ is an elementary subdivision of $K$\/} if
$K'$ is obtained from $K$ by one of the following two operations:
\begin{enumerate}
\item[(P1)]
Any two edges $a$ and $a^{-1}$ in $K$ are replaced by
$bc$ and $c^{-1}b^{-1}$ in all boundaries, where $b, c$ are 
distinct edges of $K'$ not in $K$.
\item[(P2)] 
Any face $A$ in $K$ with boundary $a_1\ldots a_{p}a_{p+1}\ldots a_n$
is replaced by two faces $A'$ and $A''$ in $K'$, with
boundaries $a_1\ldots a_{p} d$ and 
$d^{-1}a_{p+1} \ldots a_n$, where $d$ is an edge in $K'$ not in $K$.
Of course, the corresponding replacement is applied to $A^{-1}$.
\end{enumerate}
We say that a cell complex $K'$ is a {\it refinement\/}
of a cell complex $K$ if $K$ and $K'$ are related
in the reflexive and transitive closure of the elementary subdivision
relation, and we say that  $K$ and $K'$ are {\it equivalent\/} 
if they are related in the least equivalence relation
containing the elementary subdivision relation.
}
\end{defin}

\medskip
As we will see shortly, every cell complex is equivalent to
some special cell complex in normal form.
First, we show that a topological space $|K|$ can
be associated with a cell complex $K$, that this space
is the same for all cell complexes equivalent to $K$, and that
it is a surface.

\medskip
Given a cell complex $K$, we associate with $K$ a topological space 
$|K|$ as follows. Let us first assume that no face has
the empty sequence as a boundary. Then, we assign to each face
$A$ a circular disk, and if the boundary of $A$ is
$a_1\ldots a_m$, we divide the boundary of the disk into
$m$ oriented arcs. These arcs, in clockwise order are named
$a_1\ldots a_m$, while the opposite arcs are named
$a_1^{-1}\ldots a_m^{-1}$. We then form the quotient space
obtained by identifying arcs having the same name in the
various disks (this requires using homeomorphisms between
arcs named identically, etc).

\medskip
We leave as an exercise to prove that equivalent cell complexes
are mapped to homeomorphic spaces, and that if $K$ represents a
triangulated complex, then $|K|$ is homeomorphic to $ K_g$.

\medskip
When $K$ has a single face $A$ with the null boundary,
by (P2), $K$ is equivalent to the cell complex with
two faces $A', A''$, where $A'$ has boundary $d$, and
$A''$ has boundary $d^{-1}$. In this case, $|K|$ must
be homeomorphic to a sphere. 

\medskip
In order to show that the space $|K|$ associated with a cell
complex is a surface, we prove that every cell complex
can be refined to a triangulated  $2$-complex.

\begin{prop}
\label{refinlem}
Every cell complex $K$ can be refined to a triangulated
$2$-complex.
\end{prop}

\proof Details are given in Ahlfors and Sario
\cite{Ahlfors}, and we only indicate the main steps.
The idea is to subdivide the cell complex by adding new
edges. Informally, it is helpful to view the process as
adding new vertices and new edges, but since vertices are
not primitive objects, this must be done via the refinement
operations (P1) and (P2). 
The first step is to split every edge $a$ into two edges $b$ and $c$
where $b\not= c$, using (P1), introducing new border
vertices $(b, c^{-1})$. The effect is that for every edge $a$
(old or new), $a$ and $a^{-1}$ lead to distinct vertices.
Then, for every boundary $B = a_1\ldots a_n$, we have $n\geq 2$,
and intuitively, we create a ``central vertex''
$\beta = (d_1,\ldots, d_n)$, and
we join this vertex $\beta$ to every  vertex
including the newly created vertices
(except $\beta$ itself). This is done as follows:
first, using (P2), split the boundary $B = a_1\ldots a_n$
into $a_1d$ and $d^{-1}a_2\ldots a_n$, and then using (P1),
split $d$ into $d_1d_n^{-1}$, getting boundaries
$d_n^{-1}a_1d_1$ and $d_1^{-1}a_2\ldots a_nd_n$. Applying (P2) to
the boundary $d_1^{-1}a_2\ldots a_nd_n$,
we get the boundaries $d_1^{-1}a_2d_2$, $d_2^{-1}a_3d_3,\ldots$,
$d_{n-1}^{-1}a_nd_n$, and   $\beta = (d_1,\ldots, d_n)$ is indeed
an inner vertex.
At the end of this step, it is easy
to verify that (C2) and (C3) are satisfied, but (C1) may not.
Finally, we split each new triangular boundary  
$a_1a_2a_3$ into four subtriangles,
by joining the middles of its three sides. 
This is done by getting $b_1c_1b_2c_2b_3c_3$, using (P1),
and then $c_1b_2d_3$, $c_2b_3d_1$, $c_3b_1d_2$, and
$d_1^{-1}d_2^{-1}d_3^{-1}$, using (P2).
The resulting cell complex also satisfies (C1), and in fact, what we have
done is to provide a triangulation.
$\square$

\medskip
Next, we need to define cell complexes in normal form.
First, we need to define what we mean by orientability of
a cell complex, and to explain how we compute its Euler-Poincar\'e
characteristic.

\begin{defin}
\label{cellor}
{\em
Given a cell complex $K = (F, E, B)$, an {\it orientation of $K$\/}
is the choice of one of the two oriented faces $A, A^{-1}$ for
every face $A\in F$. An orientation is {\it coherent\/} if
for every edge $a$, if $a$ occurs twice in the boundaries, then
$a$ occurs in the boundary of a face $A_1$ and 
in the boundary of a face $A_2^{-1}$, where $A_1\not= A_2$.
A cell complex $K$ is {\it orientable\/} if is has some coherent
orientation. A {\it contour\/} of a cell complex
is a cyclically ordered sequence $(a_1,\ldots, a_n)$ of edges
such that $a_i$ and $a_{i+1}^{-1}$ lead to the same vertex,
and the $a_i$ belong to a single boundary.
}
\end{defin}

\medskip
It is easily seen that equivalence of cell complexes
preserves orientability. In counting contours, we do not
distinguish between  $(a_1,\ldots, a_n)$ and  
$(a_n^{-1},\ldots, a_1^{-1})$.
It is easily verified that (P1) and (P2) do not change the number
of contours.

\medskip
Given a cell complex $K = (F, E, B)$, the number of vertices is
denoted as $n_0$, the number $n_1$ of edges  is the number
of elements in $E$, and the  number $n_2$ of faces is the number
of elements in $F$. The Euler-Poincar\'e characteristic
of $K$ is $n_0 - n_1 + n_2$.
It is easily seen that (P1) increases $n_1$ by $1$,
creates one more vertex, and leaves $n_2$ unchanged.
Also, (P2) increases $n_1$ and $n_2$ by $1$ and leaves
$n_0$ unchanged. Thus, equivalence preserves the
Euler-Poincar\'e characteristic.  However, we need a small
adjustment in the case where $K$ has a single face $A$
with the null boundary. In this case, we agree that
$K$ has the ``null vertex '' $\epsilon$.
We now define the normal forms of cell complexes.
As we shall see, these normal forms have a single face and 
a single inner vertex.

\begin{defin}
\label{normcell}
{\em
A {\it cell complex in normal form, or canonical cell complex\/},
is a cell complex $K = (F, E, B)$, where
$F = \{A\}$ is a singleton set, and either
\begin{enumerate}
\item[(I)] 
$E= \{a_1,\ldots,a_p,b_1,\ldots,b_p,c_1,\ldots,c_q,h_1,\ldots,h_q\}$,
and
\[B(A) = a_1b_1a_1^{-1}b_{1}^{-1}\cdots a_pb_pa_p^{-1}b_{p}^{-1}
c_1h_1c_1^{-1}\cdots c_qh_qc_q^{-1},\]
where $p\geq 0$, $q\geq 0$, or
\item[(II)]  
$E= \{a_1,\ldots,a_p,c_1,\ldots,c_q,h_1,\ldots,h_q\}$,
and
\[B(A) = a_1a_1\cdots a_pa_p c_1h_1c_1^{-1}\cdots c_qh_qc_q^{-1},\]
where $p\geq 1$, $q\geq 0$.
\end{enumerate}
}
\end{defin}

Observe that canonical complexes of type (I) are orientable,
whereas canonical complexes of type (II) are not.
The sequences $c_ih_ic_i^{-1}$  yield $q$ border vertices
$(h_i, c_i, h_i^{-1})$, and thus $q$ contours $(h_i)$,
and in case (I), the single inner vertex
\[(a_1^{-1}, b_1, a_1, b_1^{-1}\ldots, a_p^{-1}, b_p, a_p, b_p^{-1}, 
c_1^{-1}, \ldots, c_q^{-1}),\]
and in case (II), the single inner vertex
\[(a_1^{-1}, a_1, \ldots, a_p^{-1}, a_p, 
c_1^{-1}, \ldots, c_q^{-1}).\]
Thus, in case (I), there are $q+1$ vertices, $2p + 2q$ sides, and
one face, and the Euler-Poincar\'e characteristic is
$q+1 -(2p + 2q) + 1 = 2 - 2p - q$, that is
\[\chi(K) = 2 -2p - q,\]
and in case (II), there are $q+1$ vertices, $p + 2q$ sides, and
one face, and the Euler-Poincar\'e characteristic is
$q + 1 - (p + 2q) + 1 = 2 - p - q$, that is 
\[\chi(K) = 2 - p - q.\]
Note that when $p = q = 0$, we do get $\chi(K) = 2$,
which agrees with the fact that in this case,
we assumed the existence of
a null vertex, and there is one face. This is the case
of the sphere.

\medskip
The above shows that distinct canonical complexes $K_1$ and
$K_2$ are inequivalent, since otherwise $|K_1|$ and  $|K_2|$ 
would be homeomorphic, which would imply that
$K_1$ and $K_2$ have the same number of contours, the
same kind of orientability, and the same
Euler-Poincar\'e characteristic.

\medskip
It remains to prove that every cell complex is
equivalent to a canonical cell complex, but first,
it is helpful to give more  intuition regarding the
nature of the canonical complexes.

\begin{figure}
  \begin{center}
    \begin{pspicture}(0,-0.5)(3,3.2)
    \psline[linewidth=1pt]{->}(0,0)(0,3)
    \psline[linewidth=1pt]{->}(0,3)(3,3)
    \psline[linewidth=1pt]{->}(3,0)(3,3)
    \psline[linewidth=1pt]{->}(0,0)(3,0)
    \uput[180](0,1.5){$a_1$}
    \uput[90](1.5,3){$b_1$}
    \uput[0](3,1.5){$a_1$}
    \uput[-90](1.5,0){$b_1$}
    \end{pspicture}
  \end{center}
  \caption{A cell complex corresponding to a torus}
\label{torfig2}
\end{figure}
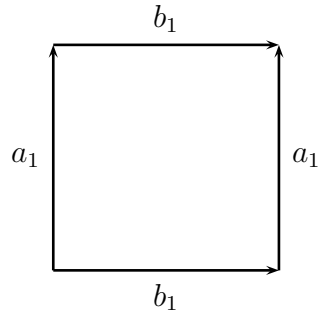

\medskip
If a canonical cell complex has the border $B(A) = a_1b_1a_1^{-1}b^{-1}_1$,
we can think of the face $A$ as a square whose opposite 
edges are oriented the same way, and labeled the same way, so that
by identification of the opposite edges labeled $a_1$ and
then of the edges labeled $b_1$, we get a surface homeomorphic to a torus.
Figure \ref{torfig2} shows such a cell complex.

If we start with a sphere and glue a torus onto the surface of the
sphere by removing some small disk from both the sphere and
the torus and gluing along the boundaries of the holes, it is as if
we had added a handle to the sphere. For this reason, 
the string $a_1b_1a_1^{-1}b^{-1}_1$ is called a {\it handle\/}.
A canonical cell complex with boundary 
$a_1b_1a_1^{-1}b_{1}^{-1}\cdots a_pb_pa_p^{-1}b_{p}^{-1}$
can be viewed as the result of attaching $p$ handles to a sphere.

\begin{figure}
  \begin{center}
    \begin{pspicture}(-1.5,-2)(1.4,1.8)
    \parametricplot[linewidth=1pt,plotstyle=curve]%
    {0}{360}{t sin 1.5 mul
             t cos 1.5 mul}
    \psline[linewidth=1pt]{->}(0,1.5)(0.01,1.5)
    \psline[linewidth=1pt]{->}(0,-1.5)(-0.01,-1.5)
    \psdots[dotstyle=o,dotscale=1.5](-1.5, 0)
    \psdots[dotstyle=o,dotscale=1.5](1.5, 0)
    \uput[90](0,1.5){$a_1$}
    \uput[-90](0,-1.5){$a_1$}
    \end{pspicture}
  \end{center}
  \caption{A cell complex corresponding to a projective plane}
\label{projplan2}
\end{figure}
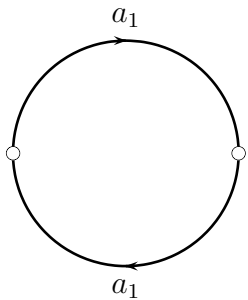

\medskip
If a canonical cell complex has the border $B(A) = a_1a_1$,
we can think of the face $A$ as a circular disk whose
boundary is divided into two semi-circles both labeled $a_1$.
The corresponding surface is obtained by identifying
diametrically opposed points on the boundary, and thus it is
homeomorphic to the projective plane. 
Figure \ref{projplan2} illustrates this situation.

There is a way of
performing such an identification resulting in a surface
with self-intersection, sometimes called a {\it cross-cap\/}.
A nice description of the process of getting a cross-cap is given in 
Hilbert and Cohn-Vossen \cite{Hilbert}.
A string of the form $aa$ is called
a {\it cross-cap\/}.
Generally, a canonical cell complex with boundary
$a_1a_1\cdots a_pa_p$ can be viewed as the result of forming
$p\geq 1$ cross-caps, starting from a  circular disk with $p-1$ circular
holes, and  performing the cross-cap identifications on
all $p$ boundaries, including the original disk itself.

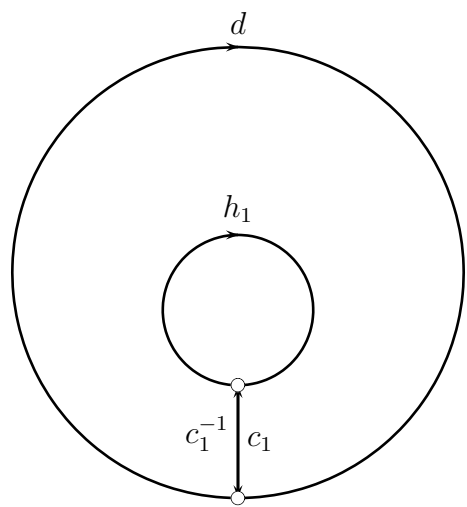
\begin{figure}
  \begin{center}
    \begin{pspicture}(-3,-3.5)(3,3.2)
    \parametricplot[linewidth=1pt,plotstyle=curve]%
    {0}{360}{t sin 3 mul
             t cos 3 mul}
    \parametricplot[linewidth=1pt,plotstyle=curve]%
    {0}{360}{t sin 
             t cos 0.5 sub}
    \psline[linewidth=1pt]{->}(0,-3)(0,-1.5)
    \psline[linewidth=1pt]{->}(0,-1.5)(0,-3)
    \psline[linewidth=1pt]{->}(0,0.5)(0.01,0.5)
    \psline[linewidth=1pt]{->}(0,3)(0.01,3)
    \psdots[dotstyle=o,dotscale=1.5](0, -3)
    \psdots[dotstyle=o,dotscale=1.5](0, -1.5)
    \uput[90](0,0.5){$h_1$}
    \uput[135](0,-2.5){$c_1^{-1}$}
    \uput[45](0,-2.5){$c_1$}
    \uput[90](0,3){$d$}
    \end{pspicture}
  \end{center}
  \caption{A disk with a hole}
\label{holefig}
\end{figure}

\medskip
A string  of the form $c_1h_1c_1^{-1}$ occurring in a border
can be interpreted as a hole with boundary $h_1$.
For instance, if the boundary of a canonical cell complex  
is $c_1h_1c_1^{-1}$, splitting the face $A$ into the two faces
$A'$ and $A''$ with boundaries $c_1h_1c_1^{-1}d$ and
$d^{-1}$, we can view the face $A'$ as a disk 
with boundary $d$ in which a small circular
disk has been removed. Choosing any point on the boundary $d$ of
$A'$, we can join this point to the boundary $h_1$ of the small circle
by an edge $c_1$, and we get a path $c_1h_1c_1^{-1}d$.
The path is a closed loop, and a string of the form
$c_1h_1c_1^{-1}$ is called a {\it loop\/}.
Figure \ref{holefig} illustrates this situation.

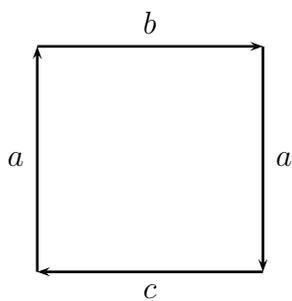
\begin{figure}
  \begin{center}
    \begin{pspicture}(0,-0.5)(3,3.2)
    \psline[linewidth=1pt]{->}(0,0)(0,3)
    \psline[linewidth=1pt]{->}(0,3)(3,3)
    \psline[linewidth=1pt]{->}(3,3)(3,0)
    \psline[linewidth=1pt]{->}(3,0)(0,0)
    \uput[180](0,1.5){$a$}
    \uput[90](1.5,3){$b$}
    \uput[0](3,1.5){$a$}
    \uput[-90](1.5,0){$c$}
    \end{pspicture}
  \end{center}
  \caption{A cell complex corresponding to a M\"obius strip}
\label{mobstrip}
\end{figure}

\medskip
We now prove a combinatorial lemma which is the key
to the classification of the compact surfaces.
First, note that the inverse of the reduction step (P1),
denoted as (P1$)^{-1}$, applies to a string of edges $bc$
provided that $b\not = c$ and $(b, c^{-1})$ is a vertex. 
The result is that such a border vertex is eliminated.
The inverse of the reduction step (P2),
denoted as (P2$)^{-1}$, applies to two faces $A_1$ and $A_2$
such that $A_1\not= A_2$, $A_1\not= A_2^{-1}$, and 
$B(A_1)$ contains some edge $d$ and $B(A_2)$ contains the edge $d^{-1}$.
The result is that $d$ (and $d^{-1}$) is eliminated.
As a preview of the proof, we show that the following cell complex,
obviously corresponding to a M\"obius strip, is equivalent to
the cell complex of type (II) with boundary
$aachc^{-1}$. The boundary of the cell complex 
shown in Figure \ref{mobstrip} is
$abac$.

First using (P2), we split $abac$ into $abd$ and $d^{-1}ac$. 
Since $abd = bda$ and the inverse face of $d^{-1}ac$
is $c^{-1}a^{-1}d = a^{-1}dc^{-1}$, by applying (P2$)^{-1}$,
we get $bddc^{-1} = ddc^{-1}b$. We can now apply (P1$)^{-1}$,
getting $ddk$. We are almost there, except that
the complex with boundary $ddk$ has no inner vertex.
We can introduce one as follows. Split $d$ into $bc$, getting
$bcbck = cbckb$. Next, apply (P2), getting
$cba$ and $a^{-1}ckb$. Since $cba = bac$ and the inverse face
of $a^{-1}ckb$ is $b^{-1}k^{-1}c^{-1}a = c^{-1}ab^{-1}k^{-1}$, 
by applying  (P2$)^{-1}$ again, we get
$baab^{-1}k^{-1} = aab^{-1}k^{-1}b$, which is of the form
$aachc^{-1}$, with $c = b^{-1}$ and $h = k^{-1}$.
Thus, the canonical cell complex with boundary $aachc^{-1}$
has the M\"obius strip as its geometric realization.
Intuitively, this corresponds to cutting out a small circular
disk in a projective plane.
This process is very nicely described in Hilbert and Cohn-Vossen
\cite{Hilbert}.

\begin{lemma}
\label{keylem}
Every cell complex $K$ is equivalent to some canonical cell complex.
\end{lemma}

\proof All the steps are given in Ahlfors and Sario
\cite{Ahlfors}, and in a slightly different and more informal 
manner in Massey \cite{Massey87}. We will only give the keys steps,
referring the reader to the above sources for details.

\medskip
The proof proceeds by steps that bring the original cell complex
closer to normal form.

\medskip
{\it Step\/} 1. Elimination of strings $aa^{-1}$ in boundaries.

\medskip
Given a boundary of the form $aa^{-1}X$, where $X$ denotes
some string of edges (possibly empty), we can use (P2) to replace
$aa^{-1}X$ by the two boundaries $ad$ and $d^{-1}a^{-1}X$, where
$d$ is new. But then, using (P1), we can contract $ad$ to a new edge
$c$ (and $d^{-1}a^{-1}$ to $c^{-1}$). But now, using (P2$)^{-1}$,
we can eliminate $c$. The net result is the elimination of $aa^{-1}$.

\medskip
{\it Step\/} 2. Vertex Reduction.

\medskip
If $p = 0, q = 0$, there is only the empty vertex, and there is
nothing to do. Otherwise,
the purpose of this step is to obtain a cell complex with a
single inner vertex, and where border vertices correspond to loops.
First, we perform step 1 until
all occurrences of the form $aa^{-1}$ have been eliminated.

\medskip
Consider an inner vertex $\alpha = (b_1,\ldots, b_m)$.
If $b_i^{-1}$ also belongs to $\alpha$ for all
$i$, $1\leq i \leq m$, and there is another inner vertex
$\beta$, since all vertices are connected, there is some
inner vertex $\delta\not=\alpha$ directly connected to $\alpha$,
which means that either some $b_i$ or $b_i^{-1}$ belongs
to $\delta$. But since the vertices form a partition of $E\cup E^{-1}$,
$\alpha = \delta$, a contradiction.

\medskip
Thus, if  $\alpha = (b_1,\ldots, b_m)$ is not the only inner vertex,
we can assume by relabeling that $b_1^{-1}$ does not
belong to $\alpha$. Also, we must have $m\geq 2$, since
otherwise there would be a string $b_1b^{-1}_1$ in some
boundary, contrary to the fact that we performed step 1 all the way.
Thus, there is a string $b_1b_2^{-1}$ in some boundary.
We claim that we can eliminate $b_2$. Indeed, since
$\alpha$ is an inner vertex, $b_2$ must occur twice in the set of
boundaries, and thus, since $b_2^{-1}$ is
a successor of $b_1$, there are boundaries
of the form $b_1b_2^{-1}X_1$ and $b_2X_2$, and using (P2),
we can split $b_1b_2^{-1}X_1$ into $b_1b_2^{-1}c$ and
$c^{-1}X_1$, where $c$ is new. Since
$b_2$ differs from $b_1, b_1^{-1}, c, c^{-1}$,
we can eliminate $b_2$ by (P2$)^{-1}$ applied to $b_2X_2 = X_2b_2$ and
$b_1b_2^{-1}c = b_2^{-1}cb_1$, getting $X_2cb_1 = cb_1X_2$.
This has the effect of shrinking $\alpha$. 
Indeed, the existence of the boundary $cb_1X_2$ implies that 
$c$ and $b_1^{-1}$ lead to the same vertex,
and the existence of the boundary $b_1b_2^{-1}c$ implies that  
$c^{-1}$ and  $b_2^{-1}$ lead to the same vertex,
and if $b_2^{-1}$ does not belong to $\alpha$, then
$b_2$ is dropped, or if  $b_2^{-1}$ belongs to $\alpha$, then
$c^{-1}$ is added to $\alpha$, but both $b_2$ and $b_2^{-1}$
are dropped.

\medskip
This process can be repeated until $\alpha = (b_1)$,
at which stage $b_1$ is eliminated using step 1.
Thus, it is possible to eliminate all inner vertices except one.
In the event that there was no inner vertex, we can always create
one using (P1) and (P2) as in the proof of Proposition \ref{refinlem}.
Thus, from now on, we will assume that there is a single inner
vertex.

\medskip
We now show that border vertices can be reduced to the form
$(h, c, h^{-1})$. The previous argument shows that we
can assume that there is a single inner vertex $\alpha$.
A border vertex is of the form $\beta = (h,b_1,\ldots, b_m,k)$,
where $h, k$ are border edges, and the $b_i$ are inner edges.
We claim that there is some border vertex $\beta = (h,b_1,\ldots, b_m,k)$
where some $b_i^{-1}$ belongs to the inner vertex $\alpha$.
Indeed, since $K$ is connected, every border vertex is
connected to $\alpha$, and thus, there is a least one
border vertex  $\beta = (h,b_1,\ldots, b_m,k)$ directly connected to
$\alpha$ by some edge.  Observe that
$h^{-1}$ and $b_1^{-1}$ lead to the same vertex, and
similarly, $b_m^{-1}$ and $k^{-1}$ lead to the same vertex.
Thus, if no  $b_i^{-1}$ belongs to $\alpha$, either
$h^{-1}$ or $k^{-1}$ belongs to $\alpha$, which
would imply that either $b_1^{-1}$ or $b_m^{-1}$ is in $\alpha$. 
Thus, such an edge from $\beta$ to $\alpha$ 
must be one of the $b_i^{-1}$.
Then by the reasoning used  in the case of an inner vertex,
we can eliminate all $b_j$ except $b_i$, and the resulting
vertex is of the form $(h, b_i, k)$.
If $h\not= k^{-1}$, we can also eliminate $b_i$
since $h^{-1}$ does not belong to $(h, b_i, k)$, and
the vertex $(h, k)$ can be eliminated using (P1$)^{-1}$.

\medskip
One can verify that reducing a border vertex to the form  
$(h, c, h^{-1})$ does not undo the reductions already performed,
and thus, at the end of step 2, we 
either obtain a cell complex with a null inner node and
loop vertices, or a single inner vertex and loop vertices.

\medskip
{\it Step\/} 3. Introduction of cross-caps.

\medskip
We may still have several faces. We claim that if there are at least
two faces, then for every face $A$, there is some face $B$ such that
$B\not= A$, $B\not= A^{-1}$, and there is some edge $a$ both
in the boundary of $A$ and in the boundary of $B$.
In this was not the case, there would be some
face $A$ such that for every face $B$ such that
$B\not= A$ and $B\not= A^{-1}$, every edge $a$ in the boundary of $B$ 
does not belong to the boundary of $A$. 
Then, every inner edge $a$ occurring in the boundary of
$A$ must have both of its occurrences in the boundary of $A$, and
of course, every border edge in the boundary of
$A$  occurs once in the boundary of $A$ alone. But then,
the cell complex consisting of the face $A$ alone and the
edges occurring in its boundary  would form a proper subsystem of $K$,
contradicting the fact that $K$ is connected.

\medskip
Thus, if there are at least two faces, from the above claim and
using (P2$)^{-1}$, we can reduce the number of faces down to one. 
It it easy to check that no new vertices are introduced,
and loops are unaffected.
Next, if some boundary contains two occurrences of the same
edge $a$, i.e., it is of the form
$aXaY$, where $X, Y$ denote strings of edges, with
$X, Y\not= \epsilon$, we show how to make the two occurrences
of $a$ adjacent. Symbolically, we show that
the following pseudo-rewrite rule is admissible:
\[aXaY \simeq bbY^{-1}X,\quad\hbox{or}\quad aaXY \simeq bYbX^{-1}.\]
Indeed, $aXaY$ can be split into $aXb$ and $b^{-1}aY$,
and since we also have the boundary 
\[(b^{-1}aY)^{-1} = Y^{-1}a^{-1}b = a^{-1}bY^{-1},\]
together with $aXb = Xba$, we can apply (P2$)^{-1}$ to
$Xba$ and $a^{-1}bY^{-1}$, obtaining $XbbY^{-1} = bbY^{-1}X$,
as claimed. Thus, we can introduce cross-caps.

\medskip
Using the formal rule $aXaY \simeq bbY^{-1}X$ again does not
alter the previous loops and cross-caps.
By repeating  step 3, we convert boundaries of the form
$aXaY$ to boundaries with cross-caps.

\medskip
{\it Step\/} 4. Introduction of handles.

\medskip
The purpose of this step is to convert boundaries of the form
$aUbVa^{-1}Xb^{-1}Y$ to boundaries $cdc^{-1}d^{-1}YXVU$
containing handles. First, we prove the pseudo-rewrite rule
\[aUVa^{-1}X \simeq bVUb^{-1}X.\]
First, we split $aUVa^{-1}X$ into $aUc = Uca$ and 
$c^{-1}Va^{-1}X = a^{-1}Xc^{-1}V$, and then we apply (P2$)^{-1}$ to
$Uca$ and $a^{-1}Xc^{-1}V$, getting $UcXc^{-1}V = c^{-1}VUcX$.
Letting $b = c^{-1}$, the rule follows.

\medskip
Now we apply the rule to $aUbVa^{-1}Xb^{-1}Y$, and we get
\[\eqaligneno{
aUbVa^{-1}Xb^{-1}Y &\simeq a_1bVUa_1^{-1}Xb^{-1}Y \cr
&\simeq a_1b_1a_1^{-1}XVUb_1^{-1}Y = a_1^{-1}XVUb_1^{-1}Ya_1b_1 \cr
&\simeq a_2^{-1}b_1^{-1} YXVU a_2b_1 = a_2b_1a_2^{-1}b_1^{-1} YXVU.\cr
}\]

\medskip
Iteration of this step preserves existing loops, cross-caps
and handles.

\medskip
At this point, one of the obstacle to the canonical form
is that we may still have a mixture of handles and cross-caps.
We now show that a handle and a cross-cap is equivalent to
three cross-caps. For this, we apply the pseudo-rewrite rule
$aaXY \simeq bYbX^{-1}$. We have

\[\eqaligneno{
aaXbcb^{-1}c^{-1}Y &\simeq a_1b^{-1}c^{-1}Ya_1c^{-1}b^{-1}X^{-1} =
b^{-1}c^{-1}Ya_1c^{-1}b^{-1}X^{-1}a_1\cr
&\simeq b_1^{-1}b_1^{-1}a_1^{-1}Xc^{-1}Ya_1c^{-1} =
c^{-1}Ya_1c^{-1}b_1^{-1}b_1^{-1}a_1^{-1}X\cr
&\simeq c_1^{-1}c_1^{-1}X^{-1}a_1b_1b_1Ya_1 = 
a_1b_1b_1Ya_1c_1^{-1}c_1^{-1}X^{-1}\cr
&\simeq a_2a_2Xc_1c_1b_1b_1Y.\cr
}\]

\medskip
At this stage, we can prove that all boundaries consist
of loops, cross-caps, or handles. The details can be found
in Ahlfors and Sario \cite{Ahlfors}.

\medskip
Finally, we have to group the loops together. This can be done using
the pseudo-rewrite rule
\[aUVa^{-1}X \simeq bVUb^{-1}X.\]
Indeed, we can write
\[chc^{-1}Xdkd^{-1}Y = c^{-1}Xdkd^{-1}Ych \simeq c_1^{-1}dkd^{-1}YXc_1h 
= c_1hc_1^{-1}dkd^{-1}YX,\]
showing that any two loops can be brought next to each other,
without altering other successions.

\medskip
When all this is done, we have obtained a canonical form, 
and the proof is complete.
$\square$

\medskip
Readers familiar with formal grammars or rewrite rules
may be intrigued by the use of the ``rewrite rules''
\[aXaY \simeq bbY^{-1}X\] 
or
\[aUVa^{-1}X \simeq bVUb^{-1}X.\]
These rules are context-sensitive, since $X$ and $Y$
stand for parts of boundaries, but they also apply to objects
not traditionally found in formal language theory or
rewrite rule theory. Indeed, the objects being rewritten are
cell complexes, which can be viewed as certain kinds of graphs.
Furthermore, since boundaries are invariant under cyclic
permutations, these rewrite rules apply modulo cyclic
permutations, something that I have never encountered in the
rewrite rule literature. Thus, it appears that a formal
treatment of such rewrite rules has not been given yet,
which poses an interesting challenge to
researchers in the field of rewrite rule theory. For example,
are such rewrite systems confluent, can normal forms
be easily found?

\medskip
We have already observed that
identification of the edges in the boundary $aba^{-1}b^{-1}$
yields a torus.  
We have also noted that identification
of the two edges in the boundary $aa$ yields the projective plane.
Lemma \ref{keylem} implies that the cell complex
consisting of a single face $A$ and the boundary
$abab^{-1}$ is equivalent to the canonical cell complex
$ccbb$. This follows immediately from the
pseudo-rewrite rule $aXaY \simeq bbY^{-1}X$.
However, it is easily seen that identification of edges
in the boundary $abab^{-1}$ yields the Klein bottle.
The lemma also showed that the cell complex with
boundary $aabbcc$ is equivalent to the cell complex with
boundary $aabcb^{-1}c^{-1}$. Thus, intuitively, it seems
that the corresponding space is a simple combination of
a projective plane and a torus, or of three projective planes. 

\medskip
We will see shortly that there is an operation on surfaces
(the connected sum) which allows us to interpret the
canonical cell complexes as combinations of elementary surfaces,
the sphere, the torus, and the projective plane.

\section{Proof of the Classification Theorem}
\label{mainthm}
Having the key Lemma \ref{keylem} at hand, we can finally prove
the fundamental theorem of the classification of
triangulated compact surfaces and compact bordered surfaces.

\begin{thm}
\label{classthm1}
Two (two-dimensional)  compact polyhedra or compact bordered 
polyhedra 
(triangulated compact surfaces or compact bordered surfaces) are
homeomorphic iff they agree in character of orientability,
number of contours, and Euler-Poincar\'e characteristic. 
\end{thm}

\proof 
If $M_1 = (K_1)_g$ and $M_2 = (K_2)_g$ are homeomorphic,
we know that $M_1$ is orientable iff $M_2$ is orientable,
and the restriction of the homeomorphism between $M_1$ and $M_2$ to
the boundaries $\partial M_1$  and $\partial M_2$,
is a homeomorphism, which implies that   $\partial M_1$  and $\partial M_2$
have the same number of arcwise components, that is, the
same number of contours. Also, we have stated that
homeomorphic spaces have isomorphic homology groups, and by
Theorem \ref{Eulercharthm}, they have the same 
Euler-Poincar\'e characteristic. 
Conversely, by Lemma \ref{keylem},
since any cell complex is equivalent to a
canonical cell complex, the triangulated $2$-complexes $K_1$ and $K_2$,
viewed as cell complexes, are equivalent to canonical cell
complexes $C_1$ and $C_2$. However, we know that equivalence
preserves orientability, the number of contours, and the
Euler-Poincar\'e characteristic, which implies that
$C_1$ and $C_2$ are identical. But then, 
$M_1 = (K_1)_g$ and $M_2 = (K_2)_g$ are both homeomorphic
to $|C_1| = |C_2|$.
$\square$

\medskip
In order to finally get a version of Theorem \ref{classthm1}
for compact surfaces or compact bordered surfaces 
(not necessarily triangulated),
we need to prove that every surface and every bordered surface
can be triangulated. This is indeed true, but the proof
is far from trivial, and it involves a strong version of the
Jordan curve theorem due to Schoenflies. At this stage,
we believe that our readers will be relieved if
we omit this proof, and refer them once again to Alhfors and Sario
\cite{Ahlfors}. It is interesting to know that
$3$-manifolds can be triangulated, but that Markov showed that
deciding whether two triangulated $4$-manifolds are homeomorphic
is undecidable (1958).
For the record, we state the following theorem
putting all the pieces of the puzzle  together.

\begin{thm}
\label{classthm2}
Two compact surfaces or compact bordered surfaces are
homeomorphic iff they agree in character of orientability,
number of contours, and Euler-Poincar\'e characteristic. 
\end{thm}

\medskip
We now explain somewhat informally what is the
connected sum operation, and how it can be used to
interpret the canonical cell complexes.
We will also indicate how the canonical cell complexes
can be used to determine the fundamental groups of the compact
surfaces and compact bordered surfaces.

\begin{defin}
\label{commsum}
{\em
Given two surfaces $S_1$ and $S_2$, their {\it connected sum\/}
$S_1\sharp S_2$ is the surface obtained by choosing two small
regions $D_1$ and $D_2$ on $S_1$ and $S_2$
both homeomorphic to some disk in the plane, 
and letting $h$ be a homeomorphism between the boundary circles
$C_1$ and $C_2$ of 
$D_1$ and $D_2$, by forming the quotient space of 
$(S_1 - \interio{D_1})\cup (S_2 - \interio{D_2})$,
by the equivalence relation defined by the relation
$\{(a, h(a))\ |\ a\in C_1\}$.
}
\end{defin}

\medskip
Intuitively, $S_1\sharp S_2$ is formed by cutting out some small
circular hole in each surface, and gluing the two surfaces along
the boundaries of these holes. It can be shown that $S_1\sharp S_2$
is a surface, and that it does not depend on the choice of
$D_1$, $D_2$, and $h$. Also, if $S_2$ is a sphere, then
$S_1\sharp S_2$ is homeomorphic to $S_1$.
It can also be shown that the Euler-Poincar\'e characteristic
of $S_1\sharp S_2$ is given by the formula
\[\chi(S_1\sharp S_2) = \chi(S_1) + \chi(S_2) - 2.\] 
Then, we can give an interpretation
of the geometric realization of a canonical cell complex.
It turns out to be the connected sum of some elementary
surfaces. Ignoring borders for the time being, assume that
we have two canonical cell complexes $S_1$ and $S_2$ 
represented by circular disks with borders
\[B_1 = a_1b_1a_1^{-1}b_{1}^{-1}\cdots 
a_{p_{1}}b_{p_{1}}a_{p_{1}}^{-1}b_{p_{1}}^{-1}\]
and
\[B_2 = c_1d_1c_1^{-1}d_{1}^{-1}\cdots 
c_{p_{2}}d_{p_{2}}c_{p_{2}}^{-1}d_{p_{2}}^{-1}.\]
Cutting a small hole with boundary $h_1$ in $S_1$
amounts to forming the new boundary
\[B_1' = a_1b_1a_1^{-1}b_{1}^{-1}\cdots 
a_{p_{1}}b_{p_{1}}a_{p_{1}}^{-1}b_{p_{1}}^{-1}h_1,\]
and similarly, 
cutting a small hole with boundary $h_2$ in $S_2$
amounts to forming the new boundary
\[B_2' = c_1d_1c_1^{-1}d_{1}^{-1}\cdots 
c_{p_{2}}d_{p_{2}}c_{p_{2}}^{-1}d_{p_{2}}^{-1}h_2^{-1}.\]
If we now glue $S_1$ and $S_2$ along $h_1$ and $h_2$
(note how we first need to reverse $B_2'$ so that
$h_1$ and $h_2$ can be glued together),
we get a figure looking like two convex polygons glued
together along one edge, and by deformation, we get a
circular disk with boundary
\[B = a_1b_1a_1^{-1}b_{1}^{-1}\cdots 
a_{p_{1}}b_{p_{1}}a_{p_{1}}^{-1}b_{p_{1}}^{-1}
c_1d_1c_1^{-1}d_{1}^{-1}\cdots 
c_{p_{2}}d_{p_{2}}c_{p_{2}}^{-1}d_{p_{2}}^{-1}.\]
A similar reasoning applies to cell complexes of type (II).

\medskip
As a consequence, the geometric realization of a cell complex
of type (I) is either a sphere, or the connected sum of $p\geq1 $ tori, 
and the geometric realization of a cell complex
of type (II) is the connected sum of $p\geq 1$ projective planes.
Furthermore, the equivalence of the cell complexes
consisting of a single face $A$ and the boundaries
$abab^{-1}$  and $aabb$, shows that the connected
sum of two projective planes is homeomorphic to the Klein bottle.
Also, the equivalence of the cell complexes with boundaries
$aabbcc$  and $aabcb^{-1}c^{-1}$ shows that
the connected sum of a projective plane and a torus is equivalent
to the connected sum of three projective planes.
Thus, we obtain another form of the classification theorem
for compact surfaces.

\begin{thm}
\label{classthm3}
Every orientable compact surface is homeomorphic either to a sphere
or to a connected sum of tori.
Every nonorientable compact surface is homeomorphic either to a 
projective plane, or a Klein bottle, or the connected sum of
a projective plane or a Klein bottle with some tori.
\end{thm}

\medskip
If  bordered compact surfaces are considered, a similar theorem
holds, but holes have to be made in the various
spaces forming the connected sum.
For more details, the reader is referred to Massey \cite{Massey87},
in which it is also shown how to build models
of bordered surfaces by gluing strips to a circular disk.

\section[Application of the Main Theorem] 
{Application of the Main Theorem: 
Determining the Fundamental Groups of Compact Surfaces} 
\label{fungrp}
We now explain  briefly how the canonical forms can be used
to determine the fundamental groups of the compact (bordered)
surfaces. This is done in two steps. The first step consists
in defining a group structure on certain closed paths in
a cell complex. The second step consists in showing
that this group is isomorphic to the fundamental group of
$|K|$.

\medskip
Given a cell complex $K = (F, E, B)$,
recall that a vertex $\alpha$ is an equivalence
class of edges, under the equivalence relation $\Lambda$ induced
by the relation $\lambda$ defined such that,
$a \lambda b$ iff $b^{-1}$ is the successor
of $a$ in some boundary. Every inner vertex 
$\alpha = (b_1,\ldots,b_m)$ can be cyclically
ordered such that $b_i$ has $b_{i-1}^{-1}$ and $b_{i+1}^{-1}$
as successors, and for a border vertex
$\alpha = (b_1,\ldots,b_m)$, the same is true for $2\leq i \leq m-1$, 
but $b_1$ only has $b_2^{-1}$ as
successor, and $b_m$ only has $b_{m-1}^{-1}$ as successor.
An edge from $\alpha$ to $\beta$ is any edge $a\in\beta$ such
that $a^{-1}\in \alpha$. For every edge $a$, we will call the
vertex that $a$ defines the {\it target\/} of $a$, and 
the vertex that $a^{-1}$ defines the {\it source\/} of $a$.
Clearly, $a$ is an edge between its source and its target.
We now define certain paths in a
cell complex, and a notion of deformation of paths.

\begin{defin}
\label{fungrp2}
{\em
Given a cell complex $K = (F, E, B)$, a {\it polygon in $K$\/} 
is any nonempty string $a_1\ldots a_m$ of edges such that
$a_i$ and $a_{i+1}^{-1}$ lead to the same vertex, or equivalently, 
such that the target of $a_i$ is equal to the source of $a_{i+1}$.
The source of the path $a_1\ldots a_m$ is the source of $a_1$
(i.e., the vertex that $a_1^{-1}$ leads to), and the
target of the path  $a_1\ldots a_m$ is the target of $a_m$
(i.e., the vertex that $a_m$ leads to).
The polygon is {\it closed\/} if its source and target coincide.
The product of two paths $a_1\ldots a_m$  and
$b_1\ldots b_n$ is defined if the target of $a_m$ is equal to
the source of $b_1$, and is the path
$a_1\ldots a_mb_1\ldots b_n$.
Given two paths $p_1 = a_1\ldots a_m$ and $p_2 = b_1\ldots b_n$
with the same source and the same target, we say that
{\it $p_2$ is an immediate deformation of $p_1$\/} if
$p_2$ is obtained from $p_1$ by either deleting
some subsequence of the form $aa^{-1}$, or deleting
some subsequence $X$ which is the boundary of some face.
The smallest equivalence relation containing the
immediate deformation relation is called  {\it path-homotopy\/}.
}
\end{defin}

\medskip
It is easily verified that path-homotopy is compatible with
the composition of paths. Then, for any vertex $\alpha_0$,
the set of equivalence classes of path-homotopic polygons
forms a group $\pi(K, \alpha_0)$.
It is also easy to see that any two groups $\pi(K, \alpha_0)$ and
$\pi(K, \alpha_1)$ are isomorphic, and that if $K_1$ and $K_2$
are equivalent cell complexes, then 
$\pi(K_1, \alpha_0)$ and $\pi(K_2, \alpha_0)$ are isomorphic.
Thus, the group $\pi(K, \alpha_0)$ only depends on the equivalence
class of the cell complex $K$. Furthermore, it can be proved that
the group $\pi(K, \alpha_0)$ is isomorphic to the fundamental group
$\pi(|K|, (\alpha_0)_g)$ associated with the geometric
realization $|K|$ of $K$ (this is proved in
Ahlfors and Sario \cite{Ahlfors}). 
It is then possible to determine
what these groups are, by considering the canonical cell complexes.

\medskip
Let us first assume that there are no borders, which 
corresponds to $q = 0$. In this case, there is only one
(inner) vertex, and all polygons are closed. For an orientable
cell complex (of type (I)), the fundamental group is the group
presented by the generators 
$\{a_1, b_1, \ldots, a_p, b_p\}$, and satisfying the single equation
\[a_1b_1a_1^{-1}b_1^{-1}\cdots a_pb_pa_p^{-1}b_p^{-1} = 1.\]
When $p = 0$, it is the trivial group reduced to $1$.
For a nonorientable cell complex (of type (II)), 
the fundamental group is the group
presented by the generators 
$\{a_1, \ldots, a_p\}$, and satisfying the single equation
\[a_1a_1\cdots a_pa_p = 1.\]

\medskip
In the presence of borders, which corresponds to $q\geq 1$,
it is easy to see that the closed polygons are products of
$a_i, b_i$, and the $d_i = c_ih_ic_i^{-1}$. For
cell complexes of type (I), these generators
satisfy the single equation
\[a_1b_1a_1^{-1}b_1^{-1}\cdots a_pb_pa_p^{-1}b_p^{-1}
d_1\cdots d_q = 1,\]
and for cell complexes of type (II), these generators
satisfy the single equation
\[a_1a_1\cdots a_pa_pd_1\cdots d_q = 1.\]
Using these equations, $d_q$ can be expressed in terms of
the other generators, and we get a free group.
In the orientable case, we get a free group with
$2q + p -1$ generators, and in the nonorientable case, we get a free
group with $p + q - 1$ generators.

\medskip
The above result shows that there are only two kinds of complexes
having a trivial group, namely for orientable complexes
for which $p = q = 0$, or $p = 0$ and $q = 1$.
The corresponding (bordered) surfaces are a {\it sphere\/}, 
and a {\it closed disk\/} (a bordered surface).
We can also figure out for which other surfaces the fundamental group
is abelian. This happens in the orientable case when
$p = 1$ and $q = 0$, a {\it torus\/}, or
$p = 0$ and $q = 2$, an {\it annulus\/}, and in the nonorientable case
when $p = 1$ and $q = 0$, a {\it projective plane\/}, 
or $p = 1$ and $q = 1$, a {\it M\"obius strip\/}.

\medskip
It is also possible to use the above results to
determine the homology groups $H_1(K)$ of the (bordered) surfaces,
since it can be shown that $H_1(K) = \pi(K,a)/[\pi(K,a),\pi(K,a)]$,
where $[\pi(K,a),\pi(K,a)]$ is the {\it commutator subgroup of
$\pi(K,a)$\/} (see Ahlfors and Sario \cite{Ahlfors}).
Recall that for any group $G$, the commutator subgroup
is the subgroup of $G$ generated by all elements of the form
$aba^{-1}b^{-1}$ (the {\it commutators\/}).
It is a normal subgroup of $G$, since
for any $h\in G$ and any $d\in [G,G]$, we have
$hdh^{-1} = (hdh^{-1}d^{-1})d$, which is also in $G$.
Then, $G/[G,G]$ is abelian, and $[G, G]$ is the smallest
subgroup of $G$ for which $G/[G,G]$ is abelian.

\medskip
Applying the above to the fundamental groups of the surfaces,
in the orientable case, we see that the commutators cause a lot
of cancellation, and we get the equation
\[d_1 + \cdots + d_q = 0,\]
whereas in the nonorientable case, we get the equation
\[2a_1 + \cdots + 2a_p + d_1 + \cdots + d_q = 0.\]
If $q > 0$, we can express $d_q$ in terms of the other
generators, and in the orientable case we get a free abelian
group with $2p + q - 1$ generators, and in the nonorientable
case a free abelian group with $p + q - 1$ generators.
When $q = 0$, in the orientable case, we get a free
abelian group with $2p$ generators, and in the nonorientable case,
since we have the equation
\[2(a_1 + \cdots + a_p) = 0,\]
there is an element of order $2$, and we get the direct sum
of a free abelian group of order $p-1$ with $\integs/2\integs$.

\medskip
Incidentally, the number $p$ is called the {\it genus\/} of a surface.
Intuitively, it counts the number of holes in the surface,
which is certainly the case in the orientable case,
but in the nonorientable case, it is considered that
the projective plane has one hole and the Klein bottle has two holes.
Of course, the genus of a surface is the number of copies of  tori
occurring in the canonical connected sum
of the surface when orientable (which, when $p = 0$, yields the sphere),
or  the number of copies of projective planes
occurring in the canonical connected sum
of the surface when  nonorientable.
In terms of the Euler-Poincar\'e characteristic, for an
orientable surface, the genus $g$ is given by the formula
\[g = (2 - \chi - q)/2,\]
and 
for a nonorientable surface, the genus $g$ is given by the formula
\[g = 2 - \chi - q,\] 
where $q$ is the number of contours.

\medskip
It is rather curious that bordered surfaces, orientable or not,
have free groups as fundamental groups
(free abelian groups for the homology groups $H_1(K)$).
It is also shown in Massey \cite{Massey87} that every bordered
surface, orientable or not, can be embedded in $\reals^3$.
This is not the case for nonorientable surfaces (with an empty border).

\medskip
Finally, we conclude with a few words about the Poincar\'e conjecture.
We observed that the only surface which is simply connected
(with a trivial fundamental group) is the sphere.
Poincar\'e conjectured in the early $1900$'s
that the same thing holds for compact
simply-connected $3$-manifolds, that is, any compact
simply-connected $3$-manifold is homeomorphic to the $3$-sphere $S^3$.

\medskip
This famous problem is still open! One of the fascinating
aspects of the Poincar\'e conjecture is that one cannot hope to have
a classification theory of compact $3$-manifolds until it is solved
(recall that $3$-manifolds can be triangulated, a result of E. Moise, 1952,
see Massey \cite{Massey87}).
What makes the Poincar\'e conjecture
even more challenging is that a generalization
of it was shown to be true by Smale for $m > 4$ in 1960, and true
for $m = 4$ by Michael Freedman in 1982.
Good luck, and let me know if you crack it!

\chapter{Topological Preliminaries}
\label{chap7}
\section{Metric Spaces and Normed Vector Spaces}
\label{sec10}
This Chapter provides a review of 
basic topological notions. 
For a comprehensive account, we highly recommend 
Munkres \cite{Munkrestop}, Amstrong \cite{Amstrong},
Dixmier \cite{Dixmier}, Singer and Thorpe \cite{Singer76},
Lang \cite{Lang97}, or Schwartz \cite{Schwartz1}. 
Most spaces considered will have a topological structure given
by a metric or a norm, and we first review these notions.
We begin with metric spaces.

\begin{defin}
\label{distdef}
{\em
A {\it metric space\/} is a set $E$ together with a function
$\mapdef{d}{E\times E}{\reals_+}$, called a {\it metric, or
distance\/},  assigning a nonnegative
real number $\dist{}{x}{y}$ to any two points $x, y\in E$,
and satisfying the following conditions for all $x,y,z\in E$:
\begin{enumerate}
\item[(D1)]
$\dist{}{x}{y} = \dist{}{y}{x}.\qquad$ \hfill (symmetry)
\item[(D2)]
$\dist{}{x}{y}\geq 0$,\ and $\dist{}{x}{y} = 0$ iff $x = y.\qquad$ 
\hfill (positivity)
\item[(D3)]
$\dist{}{x}{z} \leq \dist{}{x}{y} + \dist{}{y}{z}.\qquad$ 
\hfill (triangular inequality)
\end{enumerate} 
}
\end{defin}

\medskip
Geometrically, condition (D3) expresses the fact that
in a triangle with vertices $x, y, z$, the length of any side
is bounded by the sum of the lengths of the other two sides.
From (D3), we immediately get
\[\absval{\dist{}{x}{y}}{\dist{}{y}{z}} \leq \dist{}{x}{z}.\]

\medskip
Let us give some examples of metric spaces.
Recall that the {\it absolute value\/}
$|x|$ of a real number $x\in\reals$ is defined such
that $|x| = x$ if $x\geq 0$,  $|x| = -x$ if $x < 0$, and
for a complex number $x = a + ib$, as $|x| = \sqrt{a^2 + b^2}$.

\begin{example}
{\em
Let $E= \reals$, and $\dist{}{x}{y} = \absval{x}{y}$,
the absolute value of $x - y$. This is the so-called natural metric
on $\reals$.
}
\end{example}

\begin{example}
{\em
Let $E= \reals^n$ (or $E = \complex^n$). We have the Euclidean metric 
\[\dist{2}{x}{y} = \eudist{x}{y}{n},\]
the distance between the points $\linvec{x}{n}$ and
$\linvec{y}{n}$.
}
\end{example}

\begin{example}
{\em
For every set $E$, we can define the {\it discrete metric\/},
defined such that $\dist{}{x}{y} = 1$ iff $x\not= y$ and
$\dist{}{x}{x} = 0$. 
}
\end{example}

\begin{example}
{\em
For any $a, b\in \reals$ such that $a< b$, we define
the following sets:
\begin{enumerate}
\item
$[a,b] = \{x\in\reals\ |\ a\leq x\leq b\},\quad$ (closed interval)
\item
$]a,b[\> = \{x\in\reals\ |\ a < x < b\},\quad$ (open interval)
\item
$[a,b[\> = \{x\in\reals\ |\ a \leq x < b\},\quad$ 
(interval closed on the left, open on the right)
\item
$]a,b]\> = \{x\in\reals\ |\ a < x \leq b\},\quad$ 
(interval open on the left, closed on the right)
\end{enumerate}
Let $E = [a, b]$, and  $\dist{}{x}{y} = \absval{x}{y}$.
Then, $([a, b], d)$ is a metric space.
}
\end{example}

We will need to define the notion of proximity in order to define
convergence of limits and continuity of functions. For this,
we introduce some standard ``small neighborhoods''.

\begin{defin}
\label{balldef}
{\em
Given a metric space  $E$ with metric $d$, for every $a\in E$,
for every $\rho\in\reals$, with $\rho > 0$, the set
\[\cloball{a}{\rho} = \{x\in E\ |\ \dist{}{a}{x} \leq \rho\}\]
is called the {\it closed ball of center $a$ and radius $\rho$\/},
the set
\[\opball{a}{\rho} = \{x\in E\ |\ \dist{}{a}{x} < \rho\}\]
is called the {\it open ball of center $a$ and radius $\rho$\/},
and the set
\[S(a,\, \rho) = \{x\in E\ |\ \dist{}{a}{x} = \rho\}\]
is called the {\it sphere of center $a$ and radius $\rho$\/}.
It should be noted that $\rho$ is finite (i.e. not $+\infty$).
A subset $X$ of a metric space $E$ is {\it bounded\/}
if there is a closed ball $\cloball{a}{\rho}$ such
that $X\subseteq \cloball{a}{\rho}$.  
}
\end{defin}

\medskip
Clearly, $\cloball{a}{\rho} = \opball{a}{\rho} \cup S(a,\, \rho)$.

\medskip
In $E = \reals$ with the distance $\absval{x}{y}$,
an open ball of center $a$ and radius $\rho$ is the open
interval $]a - \rho, a + \rho[$. In $E = \reals^2$ 
with the Euclidean metric,
an open ball of center $a$ and radius $\rho$ is the set of points
inside the disk of center $a$ and radius $\rho$, excluding
the boundary points on the circle. In $E = \reals^3$
with the Euclidean metric,
an open ball of center $a$ and radius $\rho$ is the set of points
inside the sphere of center $a$ and radius $\rho$, excluding
the boundary points on the sphere. 

\medskip
One should be aware that intuition can be misleading in
forming a geometric image of a closed (or open) ball.
For example, if $d$ is the discrete metric, a closed ball 
of center $a$ and radius $\rho < 1$ consists only of its center $a$,
and a closed ball 
of center $a$ and radius $\rho \geq 1$ consists of the entire space!

\danger
If $E = [a, b]$, and  $\dist{}{x}{y} = \absval{x}{y}$,
as in example 4, an open ball $\opball{a}{\rho}$,
with $\rho < b - a$, is in fact the interval $[a, a + \rho[$,
which is closed on the left. 

\medskip
We now consider a very important special case of metric
spaces, normed vector spaces. 

\begin{defin}
\label{normspacdef}
{\em
Let $E$ be a vector space over a field $K$, where $K$
is either the field $\reals$ of reals, or the field $\complex$
of complex numbers. 
A {\it norm on $E$\/} is a function
$\mapdef{\smnorme{\;}}{E}{\reals_+}$,
assigning a nonnegative
real number $\smnorme{\novect{u}}$ to any vector $\novect{u}\in E$,
and satisfying the following conditions for all 
$\novect{x}, \novect{y}, \novect{z}\in E$:
\begin{enumerate}
\item[(N1)]
$\smnorme{\novect{x}}\geq 0$,\ and $\smnorme{\novect{x}} = 0$ 
iff $\novect{x} = \novect{0}.\qquad$ \hfill (positivity)
\item[(N2)]
$\smnorme{\lambda\novect{x}} = |\lambda|\,\smnorme{\novect{x}}.\qquad$ 
\hfill (scaling)
\item[(N3)]
$\smnorme{\novect{x} + \novect{y}}\leq \smnorme{\novect{x}} + 
\smnorme{\novect{y}}.\qquad$ 
\hfill (convexity inequality)
\end{enumerate} 
A vector space $E$ together with a norm $\smnorme{\,}$ is called
a {\it normed vector space\/}.
}
\end{defin}

\medskip
From (N3), we easily get
\[\absval{\smnorme{\novect{x}}}{\smnorme{\novect{y}}} \leq
\smnorme{\novect{x} - \novect{y}}.\]

\medskip
Given a normed vector space $E$, if we define $d$ such that
\[\dist{}{\novect{x}}{\novect{y}} = \smnorme{\novect{x} - \novect{y}},\]
it is easily seen that $d$ is a metric. Thus, every normed vector
space is immediately a metric space. Note that the metric
associated with a norm is invariant under translation, that is,
\[\dist{}{\novect{x}+\novect{u}}{\novect{y}+\novect{u}} 
= \dist{}{\novect{x}}{\novect{y}}.\]
For this reason, we can restrict ourselves to open or closed
balls of center $\novect{0}$.

\medskip
Let us give some examples of normed vector spaces.

\begin{example}
{\em
Let $E= \reals$, and $\smnorme{x} = |x|$,
the absolute value of $x$. The associated metric is 
$\absval{x}{y}$, as in example 1.
}
\end{example}

\begin{example}
{\em
Let $E= \reals^n$ (or $E = \complex^n$). There are three standard norms.
For every $\linvec{x}{n}\in E$, we have the norm
$\smnorme{\novect{x}}_1$, defined such that,
\[\smnorme{\novect{x}}_1 = |x_1| + \cdots + |x_n|,\]
we have the Euclidean norm $\smnorme{\novect{x}}_2$, defined such that,
\[\smnorme{\novect{x}}_2 = \eunorme{x}{n},\]
and the {\it sup\/}-norm $\smnorme{\novect{x}}_\infty$, defined such that, 
\[\smnorme{\novect{x}}_\infty = \max \{|x_i|\ |\ 1\leq i\leq n\}.\]
}
\end{example}

Some work is required to show the convexity inequality for
the Euclidean norm, but this can be found in any standard text.
Note that the Euclidean distance is the distance associated
with the Euclidean norm.
The following proposition is easy to show.

\begin{prop}
\label{comparnorm}
The following inequalities hold for all $\novect{x}\in\reals^n$
(or $\novect{x}\in\complex^n$):
\[\eqaligneno{
\smnorme{\novect{x}}_\infty &\leq \smnorme{\novect{x}}_1
\leq n\smnorme{\novect{x}}_\infty,\cr
\smnorme{\novect{x}}_\infty &\leq \smnorme{\novect{x}}_2
\leq \sqrt{n}\smnorme{\novect{x}}_\infty,\cr
\smnorme{\novect{x}}_2 &\leq \smnorme{\novect{x}}_1
\leq \sqrt{n}\smnorme{\novect{x}}_2.\cr
}\]
\end{prop}

\medskip
In a normed vector space, we define a closed ball or an open ball
of radius $\rho$ as a closed ball or an open ball of center
$\novect{0}$. We may use the notation $\ncloball{\rho}$
and  $\nopball{\rho}$.

\medskip
We will now define the crucial notions of open sets and closed
sets, and of a topological space.

\begin{defin}
\label{opendef}
{\em
Let $E$ be a metric space with metric $d$. A subset $U\subseteq E$
is an {\it open set\/} in $E$ if either $U = \emptyset$, or
for every $a\in U$, there is some open ball $\opball{a}{\rho}$ such that, 
$\opball{a}{\rho}\subseteq U$.%
\footnote{Recall that $\rho > 0$.}
A subset $F\subseteq E$ is a {\it closed set\/} in $E$ if
its complement $E - F$ is open in $E$.
}
\end{defin}

\medskip
The set $E$ itself is open, since for every $a\in E$, every open ball
of center $a$ is contained in $E$. 
In $E= \reals^n$, given 
$n$ intervals $[a_i, b_i]$, with $a_i < b_i$,
it is easy to show that the open $n$-cube
\[\{\linvec{x}{n}\in E\ |\ a_i < x_i < b_i,\, 1\leq i\leq n\}\]
is an open set. In fact, it is possible to find a metric
for which such open $n$-cubes are open balls!
Similarly, we can define the closed $n$-cube
\[\{\linvec{x}{n}\in E\ |\ a_i \leq  x_i \leq  b_i,\, 1\leq i\leq n\},\]
which is a closed set.

\medskip
The open sets satisfy some important properties
that lead to the definition of a topological space.

\begin{prop}
\label{metriclem}
Given a metric space $E$ with metric $d$, the family $\s{O}$
of open sets defined in Definition \ref{opendef} satisfies
the following properties:
\begin{enumerate}
\item[(O1)]
For every finite family $(U_i)_{1\leq i\leq n}$ of sets $U_i\in\s{O}$,
we have  $U_1\cap \cdots \cap U_n\in\s{O}$,
i.e. $\s{O}$ is closed under  finite intersections.
\item[(O2)]
For every arbitrary family $\famil{U}{i}{I}$
of sets  $U_i\in\s{O}$, we have  $\bigcup_{i\in I} U_i\in\s{O}$,
i.e. $\s{O}$ is closed under arbitrary unions.
\item[(O3)]
$\emptyset\in\s{O}$, and $E\in\s{O}$, i.e. $\emptyset$ and $E$
belong to $\s{O}$.
\end{enumerate}
Furthermore, for any two distinct points $a\not= b$ in $E$,
there exist two open sets $U_a$ and $U_b$ such that,
$a\in U_a$, $b\in U_b$, and $U_a \cap U_b = \emptyset$.
\end{prop}

\proof
It is straightforward. For the last point,
letting $\rho = \dist{}{a}{b}/3$ (in fact $\rho = \dist{}{a}{b}/2$ works too),
we can pick $U_a = \opball{a}{\rho}$ and $U_b = \opball{b}{\rho}$.
By the triangle inequality, we must have $U_a\cap U_b=\emptyset$.
$\square$

\medskip
The above proposition leads to the very general concept of
a topological space.

\danger
One should be careful that in general, the family of open
sets is not closed under infinite intersections.
For example, in $\reals$ under the metric $\absval{x}{y}$,
letting $U_n = ]-1/n,\, +1/n[$, each $U_n$ is open, but
$\bigcap_{n} U_n = \{0\}$, which is not open.

\section{Topological Spaces, Continuous Functions, Limits}
\label{sec11}
Motivated by Proposition \ref{metriclem}, a topological space
is defined in terms of a family of sets satisfing the properties
of open sets stated in that proposition.

\begin{defin}
\label{topspac}
{\em
Given a set $E$, a {\it topology on $E$ (or a topological
structure on $E$\/}), is defined as a family $\s{O}$ of
subsets of $E$ called {\it open sets\/}, and satisfying
the following three properties:
\begin{enumerate}
\item[(1)]
For every finite family $(U_i)_{1\leq i\leq n}$ of sets $U_i\in\s{O}$,
we have  $U_1\cap \cdots \cap U_n\in\s{O}$,
i.e. $\s{O}$ is closed under  finite intersections.
\item[(2)]
For every arbitrary family $\famil{U}{i}{I}$
of sets  $U_i\in\s{O}$, we have  $\bigcup_{i\in I} U_i\in\s{O}$,
i.e. $\s{O}$ is closed under arbitrary unions.
\item[(3)]
$\emptyset\in\s{O}$, and $E\in\s{O}$, i.e. $\emptyset$ and $E$
belong to $\s{O}$.
\end{enumerate}
A set $E$ together with a topology $\s{O}$ on $E$ is called
a {\it topological space\/}.
Given a topological space $(E, \s{O})$, 
a subset $F$ of $E$ is a {\it closed set\/} if
$F = E - U$ for some open set $U\in\s{O}$, i.e. $F$ is the
complement of some open set.

\danger
It is possible that an open set is also a closed set.
For example, $\emptyset$ and $E$ are both open and closed.
When a topological space contains a proper nonempty
subset $U$ which is both open and closed, the space
$E$ is said to be {\it disconnected\/}. 
Connected spaces will be studied in Section \ref{secsurt1}.

\medskip
A topological space $(E, \s{O})$ is 
said to satisfy the {\it Hausdorff separation axiom
(or $T_2$-separation axiom)\/} if
for any two distinct points $a\not= b$ in $E$,
there exist two open sets $U_a$ and $U_b$ such that,
$a\in U_a$, $b\in U_b$, and $U_a \cap U_b = \emptyset$.
When the $T_2$-separation axiom is satisfied, we also say that
$(E, \s{O})$ is a {\it Hausdorff space\/}.
}
\end{defin}

\medskip
As shown by Proposition \ref{metriclem}, any metric space
is a topological Hausdorff space, the family of open sets
being in fact the family of arbitrary unions of open balls.
Similarly, any normed vector space is a topological Hausdorff space,
the family of open sets being the family of arbitrary unions of open balls.
The topology $\s{O}$ consisting of all subsets of $E$ is
called the {\it discrete topology\/}.

\remark
Most (if not all) spaces used in analysis are Hausdorff
spaces. Intuitively, the Hausdorff separation axiom says that there
are enough ``small'' open sets. Without this axiom, some counter-intuitive
behaviors may arise. For example, a sequence may have more than one
limit point (or a compact set may not be closed). Nevertheless, non-Hausdorff
topological spaces arise naturally in algebraic geometry. But even there,
some substitute for separation is used.
\endremark

One of the reasons why topological spaces are important
is that the definition of a topology only involves a
certain family $\s{O}$ of sets, and not {\bf how}
such family is generated from a metric or a norm.
For example, different metrics or different norms can
define the same family of open sets. Many topological properties
only depend on the family $\s{O}$ and not on the specific
metric or norm. But the fact that a topology
is definable from a metric or a norm is important, because it usually
implies nice properties of a space. 
All our examples will be spaces whose topology is defined
by a metric or a norm.

\medskip
By taking complements, we can state properties of the closed
sets dual to those of Definition \ref{topspac}.
Thus, $\emptyset$ and $E$ are closed sets, and the
closed sets are closed under finite unions and
arbitrary intersections. It is also worth noting that
the Hausdorff separation axiom implies that for every $a\in E$,
the set $\{a\}$ is closed. Indeed, if $x\in E - \{a\}$, then
$x\not= a$, and so there exist open sets $U_a$ and $U_x$
such that $a\in U_a$, $x\in U_x$, and $U_a\cap U_x = \emptyset$.
Thus, for every  $x\in E - \{a\}$, there is an open set $U_x$
containing $x$ and contained in $E - \{a\}$, showing by (O3) that
$E - \{a\}$ is open, and thus that the set $\{a\}$ is closed. 

\medskip
Given a topological space  $(E, \s{O})$, given any subset
$A$ of $E$, since $E\in\s{O}$ and $E$ is a closed set,
the family $\s{C}_A = \{F\ |\ A\subseteq F,\, \hbox{$F$ a closed set}\}$
of closed sets containing $A$ is nonempty,
and since any arbitrary intersection of closed sets is 
a closed set, the intersection $\bigcap \s{C}_A$ of the sets
in the family  $\s{C}_A$ is the smallest closed set
containing $A$. By a similar reasoning, the union of all
the open subsets contained in $A$ is the largest open set
contained in $A$. 

\begin{defin}
\label{closurdef}
{\em
Given a topological space  $(E, \s{O})$, given any subset
$A$ of $E$, the smallest closed set  containing  $A$ is denoted
as $\adher{A}$, and is called the {\it closure, or adherence of $A$\/}.
A subset $A$ of $E$ is {\it dense in $E$\/} if $\adher{A} = E$.  
The largest open set contained in $A$ is denoted as $\interio{A}$,
and is called the {\it interior of $A$\/}.
The set $\fr{A} = \adher{A}\>\cap \adher{E-A}$, is called the 
{\it boundary (or frontier) of $A$\/}. We also denote the boundary
of $A$ as $\dBd A$.
}
\end{defin}

\remark
The notation $\adher{A}$ for the closure of a subset $A$ of $E$
is somewhat unfortunate, since  $\adher{A}$
is often used to denote the set complement of $A$ in $E$.
Still, we prefer it to more cumbersome notations such as {\rm clo}$(A)$,
and we denote the complement of $A$ in $E$ as $E - A$.
\endremark

By definition, it is clear that a subset $A$ of $E$ is closed
iff $A = \adher{A}$. 
The set $\rats$ of rationals is dense in $\reals$.
It is easily   shown that 
$\adher{A} =\> \interio{A}\cup\>\dBd A$ and
$\interio{A}\cap\>\dBd A = \emptyset$.
Another useful characterization of $\adher{A}$ is
given by the following proposition.

\begin{prop}
\label{closlem1}
Given a topological space  $(E, \s{O})$, given any subset
$A$ of $E$, the closure $\adher{A}$ of $A$ is the set of
all points $x\in E$ such that for every open set $U$ containing $x$,
then $U\cap A \not= \emptyset$.
\end{prop}

\proof If $A = \emptyset$, since $\emptyset$ is closed,
the proposition holds trivially. Thus, assume that $A\not= \emptyset$.
First, assume that $x\in \adher{A}$.
Let $U$ be any open set such that $x\in U$. 
If $U\cap A = \emptyset$, since $U$ is open,
then $E - U$ is a closed set containing $A$, and since
$\adher{A}$ is the intersection of all closed sets containing $A$,
we must have $x\in E - U$, which is impossible.
Conversely, assume that $x\in E$ is a point such that
for  every open set $U$ containing $x$,
then $U\cap A \not= \emptyset$. Let $F$ be any closed subset
containing $A$. If $x\notin F$, since $F$ is closed,
then $U = E - F$ is an open set such that $x\in U$, and
$U\cap A = \emptyset$, a contradiction. Thus, we have $x\in F$
for every closed set containing $A$, that is, $x\in\adher{A}$.
$\square$

\medskip
Often, it is necessary to consider a subset $A$ of a topological
space $E$, and to view the subset $A$ as a topological space.
The following proposition shows how to define a topology on a subset.

\begin{prop}
\label{subtoplem1}
Given a topological space  $(E, \s{O})$, given any subset
$A$ of $E$, let 
\[\s{U} = \{U\cap A\ |\ U\in\s{O}\}\]
be the family of all subsets of $A$ obtained as  the intersection
of any open set in $\s{O}$ with $A$.  
The following properties hold.
\begin{enumerate}
\item[(1)] 
The space $(A, \s{U})$ is a topological space.
\item[(2)] 
If $E$ is a metric space with metric $d$, then
the restriction $\mapdef{d_{A}}{A\times A}{\reals_+}$
of the metric $d$ to $A$ defines a metric space.
Furthermore,  the topology induced by the metric $d_A$
agrees with the topology defined by $\s{U}$, as above.
\end{enumerate}
\end{prop}

\proof Left as an exercise.
$\square$

\medskip
Proposition \ref{subtoplem1} suggests the following definition.

\begin{defin}
\label{subtopdef}
{\em
Given a topological space  $(E, \s{O})$, given any subset
$A$ of $E$, the {\it subspace topology on $A$ induced by $\s{O}$\/}
is the family $\s{U}$ of open sets defined such that
\[\s{U} = \{U\cap A\ |\ U\in\s{O}\}\]
is the family of all subsets of $A$ obtained as  the intersection
of any open set in $\s{O}$ with $A$.  
We say that $(A, \s{U})$ has the {\it subspace topology\/}.
If $(E, d)$ is a metric space, 
the restriction $\mapdef{d_{A}}{A\times A}{\reals_+}$
of the metric $d$ to $A$ is called the {\it subspace metric\/}.
}
\end{defin}

\medskip
For example, if $E = \reals^n$
and $d$ is the Euclidean metric, 
we obtain the subspace topology on the closed $n$-cube
\[\{\linvec{x}{n}\in E\ |\ a_i \leq  x_i \leq  b_i,\, 1\leq i\leq n\}.\]

\danger
One should realize that every open set $U\in\s{O}$
which is entirely contained in $A$ is also in the family $\s{U}$,
but $\s{U}$ may contain open sets that are not in $\s{O}$.
For example, if $E = \reals$ with $\absval{x}{y}$,
and $A = [a, b]$, then sets of the form $[a,c[$,
with $a < c < b$ belong to $\s{U}$, but they are not open
sets for $\reals$ under  $\absval{x}{y}$.
However, there is agreement in the following situation.

\begin{prop}
\label{subtoplem2}
Given a topological space  $(E, \s{O})$, given any subset
$A$ of $E$, if $\s{U}$ is the subspace topology, then
the following properties hold.
\begin{enumerate}
\item[(1)] 
If $A$ is an open set $A\in \s{O}$, then
every open set $U\in\s{U}$ is an open set $U\in\s{O}$.
\item[(2)] 
If $A$ is a closed set in $E$, then
every closed set w.r.t. the subspace topology is
a closed set w.r.t. $\s{O}$.
\end{enumerate}
\end{prop}

\proof Left as an exercise.
$\square$

\medskip
The concept of product topology is also useful.
We have the following proposition.

\begin{prop}
\label{prodspaclem1}
Given $n$ topological spaces  $(E_i, \s{O}_i)$,
let $\s{B}$ be the family of subsets of $E_1\times\cdots\times E_n$
defined as follows:
\[\s{B} = \{U_1\times\cdots\times U_n\ |\ U_i\in\s{O}_i,\, 1\leq i\leq n\},\]
and let $\s{P}$ be the family consisting of arbitrary unions
of sets in $\s{B}$, including $\emptyset$.
Then, $\s{P}$ is a topology on $E_1\times\cdots\times E_n$.
\end{prop}

\proof Left as an exercise.
$\square$

\begin{defin}
\label{prodspacdef}
{\em
Given $n$ topological spaces  $(E_i, \s{O}_i)$,
the {\it product topology on $E_1\times\cdots\times E_n$\/}
is the family  $\s{P}$ of subsets of $E_1\times\cdots\times E_n$
defined as follows: if
\[\s{B} = \{U_1\times\cdots\times U_n\ |\ U_i\in\s{O}_i,\, 1\leq i\leq n\},\]
then $\s{P}$ is the family consisting of arbitrary unions
of sets in $\s{B}$, including $\emptyset$.
}
\end{defin}

\medskip
If each $(E_i, \smnorme{\>}_i)$ is a normed vector
space, there are three natural norms that can be defined
on $E_1\times\cdots\times E_n$:
\[\eqaligneno{
\smnorme{\linvec{x}{n}}_1 &= \smnorme{x_1}_1 + \cdots + \smnorme{x_n}_n,\cr
\smnorme{\linvec{x}{n}}_2 &= 
\Big(\smnorme{x_1}^2_1 + \cdots + \smnorme{x_n}^2_n\Big)^{\frac{1}{2}},\cr
\smnorme{\linvec{x}{n}}_{\infty} &= 
\max\{\smnorme{x_1}_1, \ldots, \smnorme{x_n}_n\}.\cr
}\]

\medskip
It is easy to show that they all define the same
topology, which is the product topology.
One can also verify that when $E_i = \reals$, with the
standard topology induced by $\absval{x}{y}$,
the topology product on $\reals^n$ is the standard
topology induced by the Euclidean norm.

\begin{defin}
\label{equivnor}
{\em
Two metrics $d_1$ and $d_2$ on a space $E$ are {\it equivalent\/}
if they induce the same topology $\s{O}$ on $E$ (i.e., they
define the same family $\s{O}$ of open sets). Similarly,
two norms $\smnorme{\>}_1$ and $\smnorme{\>}_2$ on a space $E$ 
are {\it equivalent\/} if they induce the same topology $\s{O}$ on $E$.
}
\end{defin}

\remark
Given a topological space $(E, \s{O})$, it
is often useful, as in Proposition \ref{prodspaclem1}, to define
the topology $\s{O}$ in terms of a subfamily $\s{B}$
of subsets of $E$. We say that a family $\s{B}$ of subsets
of $E$ is a {\it basis for the topology $\s{O}$\/}
if $\s{B}$ is a subset of $\s{O}$ and 
if every open set $U$ in $\s{O}$ can be obtained as some union
(possibly infinite) of sets in $\s{B}$ (agreeing that
the empty union is the empty set). 
It is immediately verified that if a family $\s{B} = (U_i)_{i\in I}$
is a basis for the topology of $(E, \s{O})$, then
$E = \bigcup_{i\in I} U_i$, and the intersection 
of any two sets $U_i, U_j\in \s{B}$ is the union of some
sets in the family $\s{B}$ (again, agreeing that
the empty union is the empty set). 
Conversely, a family $\s{B}$ with these properties
is the basis of the topology obtained by forming arbitrary unions
of sets in $\s{B}$.
\endremark

A {\it subbasis for $\s{O}$\/}
is a family $\s{S}$ of subsets of $E$, such that the family
$\s{B}$ of all finite intersections of sets in $\s{S}$
(including $E$ itself, in case of the empty intersection)
is a basis of $\s{O}$.

\medskip
We now consider the fundamental property of continuity.

\begin{defin}
\label{continudef}
{\em
Let $(E, \s{O}_E)$ and $(F, \s{O}_F)$ be topological spaces, and
let $\mapdef{f}{E}{F}$ be a function.
For every $a\in E$, we say that {\it $f$ is continuous at $a$\/}
if for every open set $V\in\s{O}_F$ containing $f(a)$,
there is some open set $U\in\s{O}_E$ containing $a$, such that
$f(U)\subseteq V$.
We say that  {\it $f$ is continuous\/} if it is continuous
at every $a\in E$.
}
\end{defin}

\medskip
Define a {\it neighborhood of $a\in E$\/} as
any subset $N$ of $E$ containing some
open set $O\in\s{O}$ such that $a\in O$. Now, if $f$ is continuous
at $a$ and $N$ is any neighborhood of $f(a)$, there is some open
set $V\subseteq N$ containing $f(a)$,
and since $f$ is continuous at $a$,
there is some open set $U$ containing $a$, such that $f(U)\subseteq V$.
Since $V\subseteq N$, the open set $U$ is a subset of $f^{-1}(N)$
containing $a$, and $f^{-1}(N)$ is a neighborhood of $a$. Conversely, 
if $f^{-1}(N)$ is a neighborhood of $a$ whenever
$N$ is any neighborhood of $f(a)$, 
it is immediate that $f$ is continuous at $a$. Thus,
we can restate Definition
\ref{continudef} as follows:

\medskip
The function $f$ is continuous at $a\in E$
iff for every neighborhood $N$ of $f(a)\in F$, then $f^{-1}(N)$
is a neighborhood of $a$.

\medskip
It is also easy to check that $f$ is continuous on $E$ iff
$f^{-1}(V)$ is an open set in $\s{O}_E$ for every
open set $V\in\s{O}_F$.

\medskip
If $E$ and $F$ are metric spaces defined by metrics $d_1$ and $d_2$,
we can show easily that $f$ is continuous at $a$ iff

\medskip\noindent
for every $\epsilon > 0$, there is some $\eta > 0$, such that,
for every $x\in E$, 
\[\hbox{if}\>\>\dist{1}{a}{x} \leq \eta,\>\> \hbox{then}\>\> 
\dist{2}{f(a)}{f(x)} \leq \epsilon.\]

\medskip
Similarly, if $E$ and $F$ are normed vector spaces defined by norms
$\smnorme{\>}_1$ and $\smnorme{\>}_2$,
we can show easily that $f$ is continuous at $\novect{a}$ iff

\medskip\noindent
for every $\epsilon > 0$, there is some $\eta > 0$, such that,
for every $\novect{x}\in E$,
\[\hbox{if}\>\>\smnorme{\novect{x} - \novect{a}}_1 \leq \eta,\>\> 
\hbox{then}\>\> 
\smnorme{f(\novect{x}) - f(\novect{a})}_2 \leq \epsilon.\]

\medskip
It is worth noting that continuity is a topological notion,
in the sense that equivalent metrics (or equivalent norms)
define exactly the same notion of continuity.

\medskip
If  $(E, \s{O}_E)$ and $(F, \s{O}_F)$ are topological spaces, and
$\mapdef{f}{E}{F}$ is a function, for every 
nonempty subset $A\subseteq E$
of $E$, we say that {\it $f$ is continuous on $A$\/} if
the restriction of $f$ to $A$ is continuous with respect
to $(A, \s{U})$ and  $(F, \s{O}_F)$, where $\s{U}$
is the subspace topology induced by $\s{O}_E$ on $A$.

\medskip
Given a product $E_1\times\cdots\times E_n$ of topological spaces,
as usual, we let $\mapdef{\pi_i}{E_1\times\cdots\times E_n}{E_i}$
be the projection function such that,
$\pi_i\linvec{x}{n} = x_i$. It is immediately verified
that each $\pi_i$ is continuous.

\medskip
Given a topological space $(E, \s{O})$, we say that
a point $a\in E$ is {\it isolated\/} if $\{a\}$ is an open set in $\s{O}$.
Then, if  $(E, \s{O}_E)$ and $(F, \s{O}_F)$ are topological spaces, 
any function $\mapdef{f}{E}{F}$ is continuous at every isolated point
$a\in E$. In the discrete topology, every point is isolated.
In a nontrivial normed vector space $(E, \smnorme{\>}$) 
(with $E\not= \{\novect{0}\}$), no point is
isolated. To show this, we show that every open ball
$\opball{\novect{u},\rho}$ contains some vectors different from $\novect{u}$.
Indeed, since $E$ is nontrivial,
there is some $\novect{v}\in E$ such that $\novect{v}\not=\novect{0}$,
and thus $\lambda = \smnorme{\novect{v}} > 0$ (by (N1)). Let 
\[\novect{w} = \novect{u} + \frac{\rho}{\lambda + 1}\novect{v}.\]
Since $\novect{v}\not=\novect{0}$ and $\rho > 0$, we have
$\novect{w}\not= \novect{u}$. Then, 
\[\smnorme{\novect{w} - \novect{u}} = \norme{\frac{\rho}{\lambda + 1}\novect{v}}
= \frac{\rho\lambda}{\lambda + 1} < \rho,\]
which shows that $\smnorme{\novect{w} - \novect{u}} < \rho$,
for $\novect{w}\not= \novect{u}$.

\medskip
The following proposition is easily shown.

\begin{prop}
\label{composlem}
Given topological spaces  $(E, \s{O}_E)$, $(F, \s{O}_F)$, and
$(G, \s{O}_G)$, and two functions $\mapdef{f}{E}{F}$ and $\mapdef{g}{F}{G}$,
if $f$ is continuous at $a\in E$ and $g$ is continuous
at $f(a)\in F$, then $\mapdef{g\circ f}{E}{G}$ is  continuous
at $a\in E$.
Given $n$ topological spaces  $(F_i, \s{O}_i)$,
for every function $\mapdef{f}{E}{F_1\times\cdots\times F_n}$,
then $f$ is continuous at $a\in E$ 
iff every  $\mapdef{f_i}{E}{F_i}$ is continuous
at $a$, where $f_i = \pi_i\circ f$.
\end{prop}

\medskip
One can also show that in a metric space $(E, d)$,
the norm $\mapdef{d}{E\times E}{\reals}$ is continuous,
where $E\times E$ has the product topology, and that
for a normed vector space $(E, \smnorme{\>})$,
the norm $\mapdef{\smnorme{\>}}{E}{\reals}$ is continuous.

\medskip
Given a function $\mapdef{f}{E_1\times\cdots\times E_n}{F}$,
we can fix $n-1$ of the arguments, say
$a_1,\ldots, a_{i-1},a_{i+1},\ldots,a_n$, and view
$f$ as a function of the remaining argument,
\[x_i\mapsto f(a_1,\ldots, a_{i-1},x_i,a_{i+1},\ldots,a_n),\]
where $x_i\in E_i$.
If $f$ is continuous, it is clear that each $f_i$ is continuous.

\danger
One should be careful that the converse is false!
For example, consider the function $\mapdef{f}{\reals\times\reals}{\reals}$,
defined such that,
\[f(x, y) = \frac{xy}{x^2 + y^2}\quad \hbox{if $(x, y)
\not= (0,0),\>\>$ and}\quad f(0,0) = 0.\]
The function $f$ is continuous on $\reals\times\reals - \{(0,0)\}$,
but on the line $y = mx$, with $m\not= 0$,
we have $f(x, y) = \frac{m}{1+m^2}\not=0$, 
and thus, on this line, $f(x,y)$ does not approach $0$ when
$(x, y)$ approaches $(0,0)$.

\medskip
The following proposition is useful for showing that
real-valued functions are continuous.

\begin{prop}
\label{lemcont}
If $E$ is a topological space, and $(\reals, \absval{x}{y})$
the reals under the standard topology, for any two functions
$\mapdef{f}{E}{\reals}$ and $\mapdef{g}{E}{\reals}$,
for any $a\in E$, for any $\lambda\in\reals$,
if $f$ and $g$ are continuous at $a$,
then $f + g$, $\lambda f$, $f\cdot g$, are continuous at $a$,
and $f/g$ is continuous at $a$ if $g(a)\not= 0$.  
\end{prop}

\proof Left as an exercise.

\medskip
Using Proposition \ref{lemcont}, we can show easily that
every real polynomial function is continuous.

\medskip
The notion of isomorphism of topological spaces is defined as follows.

\begin{defin}
\label{homeodef}
{\em
Let $(E, \s{O}_E)$ and $(F, \s{O}_F)$ be topological spaces, and
let $\mapdef{f}{E}{F}$ be a function.
We say that {\it $f$ is a homeomorphism between $E$ and $F$\/}
if $f$ is bijective, and 
both $\mapdef{f}{E}{F}$ and $\mapdef{f^{-1}}{F}{E}$ are
continuous.
}
\end{defin}

\danger
One should be careful that a bijective continuous
function $\mapdef{f}{E}{F}$ is not necessarily an homeomorphism.
For example, if $E= \reals$ with the discrete topology,
and $F=\reals$ with the standard topology, the identity
is not a homeomorphism. Another interesting example involving
a parametric curve is given below.
Let $\mapdef{L}{\reals}{\reals^2}$ be the function, defined such that,
\[\eqaligneno{
L_1(t) &= \frac{t(1 + t^2)}{1 + t^4},\cr
L_2(t) &= \frac{t(1 - t^2)}{1 + t^4}.\cr
}\]

\medskip
If we think of $(x(t), y(t))=(L_1(t),L_2(t))$ 
as a geometric point in $\reals^2$,
the set of points $(x(t),y(t))$ obtained by letting $t$ vary 
in $\reals$ from $-\infty$ to $+\infty$, 
defines a curve having the shape of a ``figure eight'',
with self-intersection at the origin, called the ``lemniscate of
Bernoulli''. The map $L$ is continuous, and in fact bijective,
but its inverse $L^{-1}$ is not continuous. Indeed, when
we approach the origin on the branch of the curve
in the upper left quadrant (i.e., points such that, $x\leq 0$, $y\geq 0$),
then $t$ goes to $-\infty$, and when
we approach the origin on the branch of the curve
in the lower right quadrant (i.e., points such that, $x\geq 0$, $y\leq 0$),
then $t$ goes to $+\infty$.

\medskip
We also review the concept of limit of a sequence.
Given any set $E$, a  {\it sequence\/} is any function
$\mapdef{x}{\natnums}{E}$, usually denoted as $(x_n)_{n\in\natnums}$,
or $(x_n)_{n\geq 0}$, or even as $(x_n)$.

\begin{defin}
\label{convergedef1}
{\em
Given a topological space  $(E, \s{O})$, 
we say that {\it a sequence $(x_n)_{n\in\natnums}$ converges to
some $a\in E$\/} if for every open  set $U$ containing $a$,
there is some $n_0\geq 0$, such that, $x_n\in U$, for all $n\geq n_0$. 
We also say that {\it a is a limit of $(x_n)_{n\in\natnums}$\/}.
}
\end{defin}

\medskip
When $E$ is a metric space with metric $d$, it is easy to show that
this is equivalent to the fact that, 

\medskip
for every $\epsilon > 0$,
there is some $n_0\geq 0$, such that, $\dist{}{x_{n}}{a} \leq \epsilon$, 
for all $n\geq n_0$. 

\medskip
When $E$ is a normed vector space with norm $\smnorme{\>}$,
it is easy to show that
this is equivalent to the fact that, 

\medskip
for every $\epsilon > 0$,
there is some $n_0\geq 0$, such that, 
$\smnorme{\novect{x_n} - \novect{a}} \leq \epsilon$, 
for all $n\geq n_0$. 

\medskip
The following proposition shows the importance of the 
Hausdorff separation axiom.

\begin{prop}
\label{Hausdorff}
Given a topological space  $(E, \s{O})$, if
the Hausdorff separation axiom holds, then every sequence
has at most one limit.
\end{prop}

\proof Left as an exercise.

\medskip
It is worth noting that the notion of limit is  topological,
in the sense that a sequence converge to a limit $b$ iff
it converges to the same limit $b$ in any equivalent metric 
(and similarly for equivalent norms).

\medskip
We still need one more concept of limit for functions.

\begin{defin}
\label{convergedef2}
{\em
Let  $(E, \s{O}_E)$ and $(F, \s{O}_F)$ be topological spaces,
let $A$ be some nonempty subset of $E$, 
and let $\mapdef{f}{A}{F}$ be a function.
For any $a\in\adher{A}$ and any $b\in F$,
we say that {\it $f(x)$ approaches $b$ as $x$ approaches $a$ with
values in $A$\/} if for every open set $V\in\s{O}_F$
containing $b$, there is some  open set $U\in\s{O}_E$
containing $a$, such that, $f(U\cap A)\subseteq V$.
This is denoted as
\[\lim_{x \to a, x\in A} f(x) = b.\]
}
\end{defin}

\medskip
First, note that by Proposition \ref{closlem1}, since $a\in\adher{A}$,
for every open set $U$ containing $a$, we have $U\cap A\not=\emptyset$,
and the definition is nontrivial. Also, even if
$a\in A$, the value $f(a)$ of $f$ at $a$ plays no role in this definition.
When $E$ and $F$ are  metric space with metrics $d_1$ and $d_2$, 
it can  be shown easily
that the definition can be stated as follows:

\medskip
for every $\epsilon > 0$, there is some $\eta > 0$, such that,
for every $x\in A$,
\[\hbox{if}\>\>  \dist{1}{x}{a} \leq \eta,\>\>\hbox{then}\>\>
\dist{2}{f(x)}{b} \leq \epsilon.\]

\medskip
When $E$ and $F$ are normed vector spaces with norms $\smnorme{\>}_1$
and $\smnorme{\>}_2$, it can  be shown easily
that the definition can be stated as follows:

\medskip
for every $\epsilon > 0$, there is some $\eta > 0$, such that,
for every $\novect{x}\in A$,
\[\hbox{if}\>\>  \smnorme{\novect{x} - \novect{a}}_1 \leq \eta,\>\>
\hbox{then}\>\>
\smnorme{f(\novect{x}) - \novect{b}}_2 \leq \epsilon.\]

\medskip
We have the following result relating continuity at a point
and the previous notion.

\begin{prop}
\label{contconrv}
Let  $(E, \s{O}_E)$ and $(F, \s{O}_F)$ be two topological spaces,
and let $\mapdef{f}{E}{F}$ be a function.
For any $a\in E$, the function $f$ is continuous at $a$ iff
$f(x)$ approaches $f(a)$ when $x$ approaches $a$ (with values in $E$).
\end{prop}

\proof Left as a trivial exercise.

\medskip
Another important proposition relating the notion of convergence of a sequence
to continuity, is stated without proof.

\begin{prop}
\label{contconrv2}
Let  $(E, \s{O}_E)$ and $(F, \s{O}_F)$ be two topological spaces,
and let $\mapdef{f}{E}{F}$ be a function.
\begin{enumerate}
\item[(1)] 
If $f$ is continuous, then for every sequence $(x_n)_{n\in\natnums}$
in $E$, if $(x_n)$ converges to $a$, then $(f(x_n))$ converges to $f(a)$.
\item[(2)] 
If $E$ is a metric space, and $(f(x_n))$ converges to $f(a)$
whenever  $(x_n)$ converges to $a$,  for every sequence $(x_n)_{n\in\natnums}$
in $E$, then $f$ is continuous.
\end{enumerate}
\end{prop}

\remark
A special case of Definition \ref{convergedef2} 
shows up in the following case: 
$E = \reals$, and $F$ is some arbitrary topological space.
Let $A$ be some nonempty subset of $\reals$, and let
$\mapdef{f}{A}{F}$ be some function. For any $a\in A$, we say
that {\it $f$ is continuous on the right at $a$\/} if
\[\lim_{x\to a, x\in A\cap [a,\,+\infty[\>} f(x) = f(a).\]
We can define continuity on the left at $a$ in a similar
fashion.
\endremark
 
We now turn to connectivity properties of topological spaces.

\section{Connected Sets}
\label{secsurt1}
Connectivity properties of topological spaces play a very important
role in understanding the topology of surfaces. 
This section gathers the facts needed to have a good understanding
of the classification theorem for compact (bordered) surfaces.
The main references are Ahlfors and Sario \cite{Ahlfors}
and Massey \cite{Massey87,Massey}. 
For general backgroud on topology, geometry, and algebraic
topology, we also highly recommend Bredon \cite{Bredon} and
Fulton \cite{Fulton95}.

\begin{defin}
\label{connecdef1}
{\em
A topological space $(E, \s{O})$ is {\it connected\/}
if the only subsets of $E$ that are both open and closed
are the empty set and $E$ itself. Equivalently, $(E, \s{O})$
is connected if $E$ cannot be written as the union
$E = U \cup V$ of two disjoint nonempty open sets $U, V$, if
$E$ cannot be written as the union
$E = U \cup V$ of two disjoint nonempty closed sets.
A subset $S\subseteq E$ is {\it connected\/} if it is connected
in the subspace topology on $S$ induced by $(E, \s{O})$.
A connected open set is called a {\it region\/}, and
a closed set is a {\it closed region\/} if its interior is
a connected (open) set.
}
\end{defin}

\medskip
Intuitively, if a space is not connected, it is possible
to define a continuous function which is constant on
disjoint ``connected components'' and which takes possibly distinct
values on disjoint components. This can  be
stated in terms of the concept of a locally constant function.
Given two topological spaces $X, Y$, a function
$\mapdef{f}{X}{Y}$ is {\it locally constant\/} if
for every $x\in X$, there is an open set $U\subseteq X$
such that $x\in X$ and $f$ is constant on $U$.

\medskip
We claim that a locally constant function is continuous.
In fact, we will prove that $f^{-1}(V)$ is open
for every subset $V\subseteq Y$
(not just for an open set $V$).
It is enough to show that $f^{-1}(y)$ is open
for every $y\in Y$, since for every subset $V\subseteq Y$,
\[f^{-1}(V) = \bigcup_{y\in V} f^{-1}(y),\]
and open sets are closed under arbitrary unions.
However, either $f^{-1}(y) = \emptyset$ if $y \in Y - f(X)$
or $f$ is constant on $U = f^{-1}(y)$ if $y\in f(X)$
(with value $y$), and since $f$ is locally
constant, for every $x\in U$, there is some open set
$W\subseteq X$ such that $x\in W$ and $f$ is constant on $W$,
which implies that $f(w) = y$ for all $w\in W$, and thus
that $W\subseteq U$, showing that $U$ is a union of open sets,
and thus is open. 
The following proposition shows that a space is connected
iff every locally constant function is  constant.

\begin{prop}
\label{conlem1a}
A topological space is connected iff every locally constant
function is constant.
\end{prop}

\proof
First, assume that $X$ is connected. 
Let $\mapdef{f}{X}{Y}$ be a locally constant 
function to some  space $Y$, and assume that $f$ is not
constant. Pick any $y\in f(Y)$. Since $f$ is not constant,
$U_1 = f^{-1}(y) \not= X$, and of course $U_1\not= \emptyset$.
We proved just before Proposition \ref{conlem1a} that 
$f^{-1}(V)$ is open for every subset $V\subseteq Y$,
and thus $U_1 = f^{-1}(y) = f^{-1}(\{y\})$ and
$U_2 = f^{-1}(Y - \{y\})$ are both open, nonempty,
and clearly $X = U_1\cup U_2$ and $U_1$ and $U_2$ are disjoint.
This contradicts the fact that $X$ is connected, and 
$f$ must be constant.

\medskip
Assume that every locally constant function 
$\mapdef{f}{X}{Y}$ to a Hausdorff  space $Y$ is constant.
If $X$ is not connected, we can write
$X = U_1 \cup U_2$, where both $U_1, U_2$ are open, disjoint, and nonempty.
We can define the function 
$\mapdef{f}{X}{\reals}$ such that $f(x) = 1$ on $U_1$
and $f(x) = 0$ on $U_2$. Since $U_1$ and $U_2$ are open, the function
$f$ is locally constant, and yet not constant, a contradiction.
$\square$

\medskip
The following standard proposition characterizing the connected
subsets of $\reals$ can be found in most topology texts
(for example, Munkres \cite{Munkrestop}, Schwartz \cite{Schwartz1}). 
For the sake of completeness,
we give a proof.

\begin{prop}
\label{conlem1}
A subset of the real line $\reals$ is connected iff
it is an interval, i.e., of the form $[a, b]$,
$]\>a, b]$, where $a = -\infty$ is possible, 
$[a, b[\>$, where $b = +\infty$ is possible, 
or $\>]a, b[\>$, where $a = -\infty$  or
$b = +\infty$ is possible.
\end{prop}

\proof
Assume that $A$ is a connected nonempty subset of $\reals$.
The cases where $A = \emptyset$ or $A$ consists of a single point
are trivial.
We show that whenever $a, b\in A$, $a < b$, then the entire interval
$[a, b]$ is a subset of $A$. Indeed, if this was not the case,
there would be some $c\in \>]a, b[\>$ such that $c\notin A$,
and then we could write 
$A = (\>]-\infty, c[\>\>\cap A) \cup (\>]c+\infty[\>\>\cap A)$,
where $\>]-\infty, c[\>\cap A$ and $\>]c+\infty[\>\cap A$ are 
nonempty and disjoint open subsets of $A$, contradicting the
fact that $A$ is connected. It follows easily that $A$ must be an interval.

\medskip
Conversely, we show that an interval $I$ must be connected.
Let $A$ be any nonempty subset of $I$ which is both open and closed in $I$.
We show that $I = A$. Fix any $x\in A$, and consider the
set $R_x$ of all $y$ such that $[x, y]\subseteq A$.
If the set $R_x$ is unbounded, then $R_x = [x, +\infty[\>$.
Otherwise, if this set is bounded, let $b$ be its least upper bound.
We claim that $b$ is the right boundary of the interval $I$.
Because $A$ is closed in $I$, unless $I$ is open on the right
and $b$ is its right boundary, we must have $b\in A$.
In the first case, $A\cap [x, b[\> = I\cap [x, b[\> = [x, b[\>$. 
In the second case, because $A$ is also open in $I$, 
unless $b$ is the right boundary
of the interval $I$ (closed on the right),
there is some open set
$\>]b - \eta, b + \eta[\>$ contained in $A$, 
which implies that $[x, b + \eta/2] \subseteq A$,
contradicting the fact that $b$ is the least upper bound of the
set $R_x$. Thus,  $b$ must be the right boundary
of the interval $I$ (closed on the right).
A similar argument applies to the set $L_y$ of all $x$ such that
$[x, y]\subseteq A$, and either $L_y$ is unbounded, or
its greatest lower bound $a$ is the left boundary
of $I$ (open or closed on the left). In all cases, we showed that
$A = I$, and the interval must be connected.
$\square$

\medskip
A characterization on the connected subsets of $\reals^n$
is harder, and requires the notion of arcwise connectedness.
One of the most important properties of connected sets is that
they are preserved by continuous maps.

\begin{prop}
\label{concont}
Given any continuous map $\mapdef{f}{E}{F}$,
if $A\subseteq E$ is connected, then $f(A)$ is connected.
\end{prop}

\proof
If $f(A)$ is not connected, then there exist some nonempty
open sets $U, V$ in $F$ such that $f(A)\cap U$ and
$f(A)\cap V$ are nonempty and disjoint, and
\[f(A) =  (f(A)\cap U) \cup (f(A)\cap V).\]
Then, $f^{-1}(U)$ and $f^{-1}(V)$ are nonempty and 
open since $f$ is continuous, and
\[A = (A\cap f^{-1}(U))\cup (A\cap f^{-1}(V)),\]
with $A\cap f^{-1}(U)$ and  $A\cap f^{-1}(V)$ nonempty,
disjoint, and open in $A$, 
contradicting the fact that $A$ is connected.
$\square$

\medskip
An important corollary of Proposition \ref{concont}
is that for every continuous function
$\mapdef{f}{E}{\reals}$, where $E$ is a connected space,
then $f(E)$ is an interval. Indeed, this follows from
Proposition \ref{conlem1}.
Thus, if $f$ takes the values $a$ and $b$ where $a < b$, then
$f$ takes all values $c\in [a, b]$.
This is a very important property.

\medskip
Even if a topological space is not connected, it turns out
that it is the disjoint union of maximal connected
subsets, and these connected components are  closed in $E$.
In order to obtain this result, we need a few lemmas.

\begin{lemma}
\label{conlem2}
Given a topological space $E$, for any family
$(A_i)_{i\in I}$ of (nonempty) connected subsets of $E$,
if $A_i\cap A_j\not=\emptyset$ for all $i, j\in I$, then
the union $A = \bigcup_{i\in I} A_i$ of the family $(A_i)_{i\in I}$
is also connected.
\end{lemma}

\proof
Assume that  $\bigcup_{i\in I} A_i$ is not connected.
Then, there exists two nonempty  open subsets $U$ and $V$
of $E$ such that $A\cap U$ and $A\cap V$ are disjoint
and nonempty, and such that
\[A = (A\cap U) \cup (A\cap V).\] 
Now, for every $i\in I$, we can write
\[A_i = (A_i\cap U) \cup (A_i\cap V),\]
where $A_i\cap U$ and $A_i\cap V$ are disjoint, since
$A_i \subseteq A$ and $A\cap U$ and $A\cap V$ are disjoint.
Since $A_i$ is connected, either
$A_i\cap U = \emptyset$ or  $A_i\cap V = \emptyset$.
This implies that either $A_i \subseteq A\cap U$ or
$A_i \subseteq A\cap V$. However, by assumption,
$A_i\cap A_j \not= \emptyset$, for all $i, j\in I$,
and thus, either both  $A_i \subseteq A\cap U$ and
$A_j \subseteq A\cap U$, or both
$A_i \subseteq A\cap V$ and $A_j \subseteq A\cap V$,
since $A\cap U$ and $A\cap V$ are disjoint. Thus, we conclude
that either $A_i \subseteq A\cap U$ for all $i\in I$, or
$A_i \subseteq A\cap V$ for all $i\in I$.
But this proves that either
\[A = \bigcup_{i\in I} A_i  \subseteq A\cap U,\] 
or
\[A = \bigcup_{i\in I} A_i  \subseteq A\cap V,\] 
contradicting the fact that both 
$A\cap U$ and $A\cap V$ are disjoint and nonempty.
Thus, $A$ must be connected.
$\square$

\medskip
In particular, the above lemma applies when the 
connected sets in a family
$(A_i)_{i\in I}$ have a point in common.

\begin{lemma}
\label{conlem3}
If $A$ is a connected subset of a topological space $E$,
then for every subset $B$ such that
$A \subseteq B \subseteq \overline{A}$, where 
$\overline{A}$ is the closure of $A$ in $E$,
the set $B$ is connected.
\end{lemma}

\proof
If $B$ is not connected, then there are two nonempty open
subsets $U, V$ of $E$ such that $B\cap U$ and $B\cap V$
are disjoint and nonempty, and
\[B =  (B\cap U) \cup (B\cap V).\]
Since $A \subseteq B$, the above implies that
\[A =  (A\cap U) \cup (A\cap V),\]
and since $A$ is connected, either 
$A\cap U = \emptyset$, or $A\cap V = \emptyset$.
Without loss of generality, assume that 
$A\cap V = \emptyset$, which implies that 
$A \subseteq A\cap U \subseteq B\cap U$.
However, $B\cap U$ is closed in the subspace topology for $B$,
and since $B \subseteq \overline{A}$ and $\overline{A}$ is closed in $E$,
the closure of $A$ in $B$ w.r.t. the subspace topology of $B$
is clearly $B\cap\overline{A} = B$, which implies that
$B\subseteq B\cap U$ (since the closure is the smallest closed set
containing the given set).
Thus, $B\cap V = \emptyset$, a contradiction.
$\square$

\medskip
In particular, Lemma \ref{conlem3} shows that
if $A$ is a connected subset, then its closure $\overline{A}$
is also connected.
We are now ready to introduce the connected components of
a space.

\begin{defin}
\label{connecdef2}
{\em
Given a topological space $(E, \s{O})$
we say that two points $a, b\in E$ are {\it connected\/}
if there is some connected subset $A$ of $E$ such that
$a\in A$ and $b\in A$.
}
\end{defin}

\medskip
It is immediately verified that the relation
``$a$ and $b$ are connected in $E$'' is an equivalence
relation. Only transitivity is not obvious, but it follows
immediately as a special case of Lemma \ref{conlem2}.
Thus, the above equivalence relation defines a partition
of $E$ into nonempty disjoint {\it connected components\/}.
The following proposition is easily proved using Lemma \ref{conlem2}
and Lemma \ref{conlem3}.

\begin{prop}
\label{conlem4}
Given any topological space $E$, for any $a\in E$, 
the connected component containing $a$ is the largest
connected set containing $a$. The connected components of $E$
are closed.
\end{prop}

\medskip
The notion of a locally connected space is also useful.

\begin{defin}
\label{loconnecdef}
{\em
A topological space $(E, \s{O})$ is {\it locally connected\/}
if for every $a\in E$, for every neighborhood $V$ of $a$, there
is a connected neighborhood $U$ of $a$ such that $U\subseteq V$.
}
\end{defin}

\medskip
As we shall see in a moment, it would be equivalent to
require that $E$ has a basis of connected open sets.

\danger
There are connected spaces that are not locally connected,
and there are locally connected spaces that are not connected.
The two properties are independent.

\begin{prop}
\label{conlem5}
A topological space $E$ is locally connected iff
for every open subset $A$ of $E$, the connected components of $A$
are open.
\end{prop}

\proof
Assume that $E$ is locally connected. Let $A$ be any open subset
of $E$, and let $C$ be one of the connected components of $A$.
For any $a\in C\subseteq A$, there is some connected neigborhood $U$
of $a$ such that $U\subseteq A$, and since $C$ is a connected
component of $A$ containing $a$, we must have $U\subseteq C$.
This shows that for every $a\in C$, there is some open subset
containing $a$ contained in $C$, and $C$ is open.

\medskip
Conversely, assume that 
for every open subset $A$ of $E$, the connected components of $A$
are open. Then, for every $a\in E$ and every neighborhood $U$ of
$a$, since $U$ contains some open set $A$ containing $a$,
the interior $\interio{U}$ of $U$ is an open set containing $a$,
and its connected components are open. In particular, the connected
component $C$ containing $a$ is a connected open set containing $a$
and contained in $U$.
$\square$

\medskip
Proposition \ref{conlem5} shows that in a locally connected space,
the connected open sets form a basis for the topology.
It is easily seen that $\reals^n$ is locally connected.
Another very important property of surfaces, and more
generally manifolds, is to be arcwise connected.
The intuition is that any two points can be joined by
a continuous arc of curve. This is formalized as follows.

\begin{defin}
\label{carcdef}
{\em
Given a topological space $(E, \s{O})$, an {\it arc (or path)\/}
is a continuous map $\mapdef{\gamma}{[a, b]}{E}$, where $[a, b]$ is
a closed interval of the real line $\reals$. The point $\gamma(a)$
is the {\it initial point\/} of the arc, and the point $\gamma(b)$ is the
{\it terminal point \/} of the arc. We say that {\it $\gamma$ is an arc joining
$\gamma(a)$ and $\gamma(b)$\/}. An arc is a {\it closed curve\/} if
$\gamma(a) = \gamma(b)$. The set $\gamma([a, b])$
is the {\it trace\/} of the arc $\gamma$.
}
\end{defin}

\medskip
Typically, $a = 0$ and $b = 1$. In the sequel, this will be assumed.

\danger
One should not confuse an arc $\mapdef{\gamma}{[a, b]}{E}$
with its trace. For example, $\gamma$ could be constant, and thus,
its trace reduced to a single point.

\medskip
An arc is a {\it Jordan arc\/} if
$\gamma$ is a homeomorphism onto its trace.
An arc  $\mapdef{\gamma}{[a, b]}{E}$ is a {\it Jordan curve\/}
if $\gamma(a) = \gamma(b)$, and $\gamma$ is injective on $[a, b[\>$.
Since $[a, b]$ is connected, by Proposition \ref{concont},
the trace $\gamma([a, b])$ of an arc is a connected subset of $E$.

\medskip
Given two arcs $\mapdef{\gamma}{[0, 1]}{E}$ and 
$\mapdef{\delta}{[0, 1]}{E}$ such that $\gamma(1) = \delta(0)$, we can form
a new arc defined as follows.

\begin{defin}
\label{arcdef}
{\em
Given two arcs $\mapdef{\gamma}{[0, 1]}{E}$ and 
$\mapdef{\delta}{[0, 1]}{E}$ such that $\gamma(1) = \delta(0)$, we can form
their {\it composition (or product) $\gamma\delta$\/}, defined such that
\[\gamma\delta(t) = \cases{\gamma(2t) & if $0\leq t\leq 1/2$;\cr
                 \delta(2t -1) & if $1/2 \leq t\leq 1$.\cr
}\]
The {\it inverse $\gamma^{-1}$ of the arc $\gamma$\/} is the arc defined
such that $\gamma^{-1}(t) = \gamma(1-t)$, for all $t\in [0,1]$. 
}
\end{defin}

\medskip
It is trivially verified that Definition \ref{arcdef} yields
continuous arcs.

\begin{defin}
\label{arcwisedef}
{\em
A topological space $E$ is {\it arcwise connected\/} if
for any two points $a, b\in E$, there is an arc
$\mapdef{\gamma}{[0, 1]}{E}$ joining $a$ and $b$,
i.e., such that $\gamma(0) = a$ and $\gamma(1) = b$.
A topological space $E$ is {\it locally arcwise connected\/}
if for every $a\in E$, for every neighborhood $V$ of $a$, there
is an arcwise connected neighborhood $U$ of $a$ such that $U\subseteq V$.
}
\end{defin}

\medskip
The space $\reals^n$ is locally arcwise connected,
since for any open ball,
any two points in this ball are joined by a line segment.
Manifolds and surfaces are also locally arcwise connected.
It is easy to verify that Proposition \ref{concont}
also applies to arcwise connectedness.
The following theorem is crucial to the theory of manifolds
and surfaces.

\begin{thm}
\label{arcwiseth}
If a topological space $E$ is arcwise connected, then it is connected.
If a topological space $E$ is connected and locally arcwise connected,
then $E$ is arcwise connected.
\end{thm}

\proof
First, assume that $E$ is arcwise connected.
Pick any point $a$ in $E$. Since $E$ is arcwise connected, 
for every $b\in E$, there is a path 
$\mapdef{\gamma_{b}}{[0,1]}{E}$ from $a$ to $b$, and so
\[E = \bigcup_{b\in E}  \gamma_b([0,1])\]
a union of connected subsets all containing $a$. 
By Lemma \ref{conlem2}, $E$ is connected.

\medskip
Now assume that $E$ is connected and locally arcwise connected.
For any point $a\in E$, let $F_a$ be the set of all points
$b$ such that there is an arc 
$\mapdef{\gamma_b}{[0,1]}{E}$ from $a$ to $b$. 
Clearly, $F_a$ contains $a$. We show that $F_a$ is both
open and closed. For any $b\in F_a$, since $E$ is locally
arcwise connected, there is an arcwise connected neighborhood $U$
containing $b$ (because $E$ is a neighborhood of $b$).
Thus, $b$  can be joined to every point $c\in U$ by an arc,
and since by the definition of $F_a$, there is an arc from
$a$ to $b$, the composition of these two arcs yields an arc from
$a$ to $c$, which shows that $c\in F_a$. But then
$U\subseteq F_a$, and thus $F_a$ is open.
Now assume that $b$ is in the complement of $F_a$.
As in the previous case, there is some 
arcwise connected neighborhood $U$ containing $b$.
Thus, every point $c\in U$ can be joined to $b$ by an arc.
If there was an arc joining $a$ to $c$, we would
get an arc from $a$ to $b$, contradicting the fact that
$b$ is in the complement of $F_a$. Thus, 
every point $c\in U$ is in the complement of $F_a$, which
shows that $U$ is contained in the complement of $F_a$,
and thus, that the  the complement of $F_a$ is open.
Consequently, we have shown that $F_a$ is both open and closed,
and since it is nonempty, we must have $E = F_a$, which
shows that $E$ is arcwise connected.
$\square$

\medskip
If $E$ is locally arcwise connected, the above argument shows
that the connected components of $E$ are arcwise connected.

\danger
It is not true that a connected space is arcwise connected.
For example, the space consisting of the graph of the function
\[f(x) = \sin(1/x),\]
where $x > 0$,
together with  the portion of the $y$-axis, for which $-1 \leq y\leq 1$,
is connected, but not arcwise connected.

\medskip
A trivial modification of the proof of Theorem \ref{arcwiseth}
shows that in a normed vector space $E$, a connected open set 
is arcwise connected by polygonal lines (i.e., arcs consisting
of line segments). This is because in every open ball,
any two points are connected by a line segment.
Furthermore, if $E$ is finite dimensional, these polygonal lines
can be forced to be parallel to  basis vectors.

\medskip
We now consider compactness.

\section{Compact Sets}
\label{secsurt2}
The property of compactness is very important in topology
and analysis. We provide a quick review geared towards the
study of surfaces, and for details, refer the reader
to Munkres \cite{Munkrestop}, Schwartz \cite{Schwartz1}.
In this section, we will need to assume that the topological
spaces are Hausdorff spaces. This is not a luxury, as many of the
results are false otherwise.

\medskip
There are various equivalent ways of defining compactness.
For our purposes, the most convenient way involves the notion of open
cover. 

\begin{defin}
\label{compadef}
{\em
Given a topological space $E$, for any
subset $A$ of $E$, an {\it open cover $(U_i)_{i\in I}$ of $A$\/}
is a family of open subsets of $E$ such that 
$A \subseteq \bigcup_{i\in I} U_i$. An  {\it open subcover\/}
of an open cover $(U_i)_{i\in I}$ of $A$ is any subfamily
$(U_j)_{j\in J}$ which is an open cover of $A$, with
$J\subseteq I$. An open cover  $(U_i)_{i\in I}$ of $A$ is {\it finite\/}
if $I$ is finite. The topological space $E$ is {\it compact\/}
if it is Hausdorff and
for every open cover $(U_i)_{i\in I}$ of $E$, there
is a finite open subcover  $(U_j)_{j\in J}$ of $E$.
Given any subset $A$ of $E$, we say that $A$ is {\it compact\/}
if it is compact with respect to the subspace topology.
We say that $A$ is {\it relatively compact\/} if its
closure $\overline{A}$ is compact.
}
\end{defin}

\medskip
It is immediately verified that a subset $A$ of $E$ is
compact in the subspace topology relative to $A$ iff
for every open cover  $(U_i)_{i\in I}$ of $A$ by open subsets of $E$,
there is a finite open subcover $(U_j)_{j\in J}$ of $A$.
The property that every open cover contains a finite open
subcover is often called the {\it Heine-Borel-Lebesgue\/} property.
By considering complements, a Hausdorff space is compact iff
for every family $(F_i)_{i\in I}$ of closed sets, if
$\bigcap_{i\in I} F_i = \emptyset$, then 
$\bigcap_{j\in J} F_j = \emptyset$ for some finite subset
$J$ of $I$.

\danger
Definition \ref{compadef} requires that a compact space be Hausdorff.
There are books in which a compact space is not necessarily
required to be Hausdorff. Following Schwartz,
we prefer calling such a space {\it quasi-compact\/}.

\medskip
Another equivalent and useful characterization can be given
in terms of families having the finite intersection property.
A family $(F_i)_{i\in I}$ of sets has the
{\it finite intersection property\/} if
$\bigcap_{j\in J} F_j \not= \emptyset$ for every finite
subset $J$ of $I$. We have the following proposition.

\begin{prop}
\label{compac1}
A topological Hausdorff space $E$ is compact iff for
every family $(F_i)_{i\in I}$ of closed sets having  the
finite intersection property, then 
$\bigcap_{i\in I} F_i \not= \emptyset$.
\end{prop}

\proof If $E$ is compact and  $(F_i)_{i\in I}$ 
is a family of closed sets having the finite intersection property,
then $\bigcap_{i\in I} F_i$ cannot be empty, since otherwise
we would have $\bigcap_{j\in J} F_j = \emptyset$ for some finite
subset $J$ of $I$, a contradiction. The converse is equally
obvious.
$\square$

\medskip
Another useful consequence of compactness is as follows.
For any family   $(F_i)_{i\in I}$ of closed sets such that
$F_{i+1} \subseteq F_{i}$ for all $i\in I$,
if $\bigcap_{i\in I} F_i = \emptyset$,
then $F_i = \emptyset$ for some $i\in I$.
Indeed, there must be some finite subset $J$ of $I$ such that
$\bigcap_{j\in J} F_j = \emptyset$, and since 
$F_{i+1} \subseteq F_{i}$ for all $i\in I$, we must have
$F_j = \emptyset$ for the smallest $F_j$ in $(F_j)_{j\in J}$. 
Using this fact, we note that $\reals$ is {\it not\/} compact.
Indeed, the family of closed sets $([n, +\infty[\>)_{n\geq 0}$
is decreasing and has an empty intersection.

\medskip
Given a metric space, if we define a {\it bounded subset\/}
to be a subset that can be enclosed in some closed ball 
(of finite radius), 
any nonbounded subset of a metric space is not compact.
However, a closed interval $[a, b]$ of the real line is compact.

\begin{prop}
\label{compac2}
Every closed interval $[a, b]$ of the real line is compact.
\end{prop}

\proof We proceed by contradiction. Let
$(U_i)_{i\in I}$ be any open cover of $[a, b]$, and assume that
there is no finite open subcover. Let $c = (a + b)/2$.
If both $[a, c]$ and $[c, b]$ had some finite open subcover,
so would $[a, b]$, and thus, either $[a, c]$ does not have
any finite subcover, or $[c, b]$  does not have any finite open subcover.
Let $[a_1, b_1]$ be such a bad subinterval. The same argument
applies, and we split $[a_1, b_1]$ into two equal subintervals,
one of which must be bad. Thus, having  defined
$[a_n, b_n]$ of length $(b - a)/2^n$ as an interval having
no finite open subcover, splitting $[a_n, b_n]$ into
two equal intervals, we know that at least one of the two
has no finite open subcover, and we denote such a bad interval as
$[a_{n+1}, b_{n+1}]$. The sequence $(a_n)$ is nondecreasing and
bounded from above by $b$, and thus, 
by a fundamental property of the real line,
it converges to its least upper bound $\alpha$. Similarly, the
sequence $(b_n)$ is nonincreasing and
bounded from below by $a$, and thus, 
it converges to its greatest lowest  bound $\beta$. Since $[a_n,b_n]$
has length $(b - a)/2^n$, we must have $\alpha = \beta$.
However, the common limit $\alpha = \beta$ of the sequences
$(a_n)$ and $(b_n)$ must belong to some open set 
$U_i$ of the open cover, and since $U_i$ is open,
it must contain some interval $[c, d]$ containing $\alpha$.
Then, because $\alpha$ is the common limit of the sequences
$(a_n)$ and $(b_n)$, there is some $N$ such that
the intervals $[a_n, b_n]$ are all contained in the interval
$[c, d]$ for all $n\geq N$, which contradicts the fact that none
of the intervals $[a_n, b_n]$ has a finite open subcover.
Thus, $[a, b]$ is indeed compact.
$\square$

\medskip
It is easy to adapt the argument of Proposition \ref{compac2}
to show that in $\reals^m$, every closed set
$[a_1,b_1]\times\cdots\times [a_m,b_m]$ is compact.
At every stage, we need to divide into $2^m$ subpieces instead
of $2$.

\medskip
The following two propositions give very important properties of the
compact sets, and they only hold for Hausdorff spaces.

\begin{prop}
\label{compac3}
Given a topological Hausdorff space $E$, for every
compact subset $A$ and every point $b$ not in $A$, there
exist disjoint open sets $U$ and $V$ such that
$A\subseteq U$ and $b\in V$.
As a consequence, every compact subset is closed.
\end{prop}

\proof Since $E$ is Hausdorff, for every $a\in A$,
there are some disjoint open sets $U_a$ and $V_b$
containing $a$ and $b$ respectively. Thus, the family
$(U_a)_{a\in A}$ forms an open cover of $A$. Since $A$ is compact
there is a finite open subcover $(U_j)_{j\in J}$ of $A$,
where $J\subseteq A$,
and then $\bigcup_{j\in J}  U_j$ is an open set containing
$A$ disjoint from the open set $\bigcap_{j\in J} V_j$ 
containing $b$. 
This shows that every point $b$ in the complement of $A$
belongs to some open set in this complement, and thus
that the complement is open, i.e., that $A$ is closed.
$\square$

\medskip
Actually, the proof of Proposition \ref{compac3} can be used to show
the following useful property.

\begin{prop}
\label{compac4}
Given a topological Hausdorff space $E$, for every
pair of compact disjoint subsets $A$ and  $B$, there
exist disjoint open sets $U$ and $V$ such that
$A\subseteq U$ and $B\subseteq V$.
\end{prop}

\proof We repeat the argument of Proposition \ref{compac3}
with $B$ playing the role of $b$, and use Proposition \ref{compac3}
to find disjoint open sets $U_a$ containing $a\in A$
and $V_a$ containing $B$.
$\square$

\medskip
The following proposition shows that in a compact topological space,
every closed set is compact. 

\begin{prop}
\label{compac5}
Given a compact topological  space $E$, 
every closed set is compact.
\end{prop}

\proof Since $A$ is closed, $E - A$ is open,
and from any open cover $(U_i)_{i\in I}$
of $A$, we can form an open cover
of $E$ by adding $E - A$ to  $(U_i)_{i\in I}$, and since $E$ is compact,
a finite subcover $(U_j)_{j\in J}\cup \{E - A\}$ of $E$ can be extracted, 
such that $(U_j)_{j\in J}$ is a finite subcover of $A$.
$\square$

\remark
Proposition \ref{compac5} also holds for quasi-compact
spaces, i.e., the Hausdorff separation property is not needed.
\endremark

Putting Proposition \ref{compac4} and Proposition \ref{compac5} together,
we note that if $X$ is compact, then
for every pair of disjoint closed sets $A$ and $B$,
there exist disjoint open sets $U$ and $V$ such that
$A \subseteq U$ and $B \subseteq V$. We say that $X$ is a {\it normal\/}
space.

\begin{prop}
\label{compac6}
Given a compact topological  space $E$, 
for every $a\in E$, for every neighborhood $V$ of $a$, there
exists a compact neighborhood $U$ of $a$ such that $U\subseteq V$.
\end{prop}

\proof Since $V$ is a neighborhood of $a$, there is some
open subset $O$ of $V$ containing $a$. Then the complement
$K = E - O$ of $O$ is closed, and since $E$ is compact,
by Proposition \ref{compac5}, $K$ is compact.
Now, if we consider the family of all closed sets of the form
$K \cap F$, where $F$ is any closed neighborhood of $a$,
since $a\notin K$, this family has an empty intersection,
and thus there is a finite number of closed neighborhood 
$F_1, \ldots, F_n$ of $a$, such that 
$K\cap F_1\cap\cdots\cap F_n = \emptyset$.
Then, $U = F_1\cap\cdots\cap F_n$ is a compact neigborhood of $a$
contained in $O\subseteq V$.
$\square$

\medskip
It can be shown that in a normed vector space of finite dimension,
a subset is compact iff it is closed and bounded. For $\reals^n$,
this is easy. 

\danger
In a normed vector space of infinite dimension, there are closed
and bounded sets that are not compact!

\medskip
More could be said about compactness in metric spaces,
but we will only need the notion of Lebesgue number, which will
be discussed a little later.
Another crucial property of compactness is that it is preserved
under continuity.

\begin{prop}
\label{compac7}
Let $E$ be a topological space, 
and $F$ be a topological Hausdorff space.
For every compact subset $A$ of $E$, for every continuous map
$\mapdef{f}{E}{F}$, the subspace $f(A)$ is compact.
\end{prop}

\proof 
Let $(U_i)_{i\in I}$ be an open cover of $f(A)$.
We claim that $(f^{-1}(U_i))_{i\in I}$ is an open
cover of $A$, which is easily checked.
Since $A$ is compact, there is a finite open subcover
$(f^{-1}(U_j))_{j\in J}$ of $A$, and thus,  $(U_j)_{j\in J}$ 
is an open subcover of $f(A)$.
$\square$

\medskip
As a corollary of Proposition \ref{compac7}, if $E$ is compact,
$F$ is Hausdorff, and $\mapdef{f}{E}{F}$ is continuous and bijective,
then $f$ is a homeomorphism. Indeed, it is enough to show
that $f^{-1}$ is continuous, which is equivalent to showing
that $f$ maps closed sets to closed sets.
However, closed sets are compact, and Proposition \ref{compac7}
shows that compact sets are mapped to compact sets,
which, by Proposition \ref{compac3}, are closed.

\medskip
It can also be shown that if $E$ is a compact nonempty space
and $\mapdef{f}{E}{\reals}$ is a continuous function,
then there are points $a, b\in E$ such that $f(a)$ is the
minimum of $f(E)$ and $f(b)$ is the maximum of $f(E)$.
Indeed, $f(E)$ is a compact subset of $\reals$, and thus
a closed and bounded set which contains its 
greatest lower bound and its least upper bound.

\medskip
Another useful notion is that of local compactness.
Indeed, manifolds and surfaces are locally compact.

\begin{defin}
\label{loccompadef}
{\em
A  topological  space $E$ is {\it locally compact\/}
if it is Hausdorff and
for every $a\in E$, there is some compact neighborhood
$A$ of $a$.
}
\end{defin}

\medskip
From Proposition \ref{compac6}, every compact space is locally compact,
but the converse is false. It can be shown that
a normed vector space of finite dimension is locally compact.

\begin{prop}
\label{compac8}
Given a locally compact topological  space $E$, 
for every $a\in E$, for every neighborhood $N$ of $a$, there
exists a compact neighborhood $U$ of $a$, such that $U\subseteq N$.
\end{prop}

\proof For any $a\in E$, there is some compact
neighborhood $V$ of $a$. By Proposition \ref{compac6},
every neigborhood of $a$ relative to $V$ contains some
compact neighborhood $U$ of $a$ relative to $V$. But every
neighborhood of $a$ relative to $V$ is a neighborhood of
$a$ relative to $E$, and every neighborhood $N$ of $a$ in $E$
yields a neighborhood $V\cap N$ of $a$ in $V$, and thus
for every neighborhood $N$ of $a$, there
exists a compact neighborhood $U$ of $a$ such that $U\subseteq N$.
$\square$

\medskip
It is much harder to deal with noncompact surfaces (or manifolds) 
than it is to deal with
compact surfaces (or  manifolds). However, surfaces
(and manifolds) are locally compact, and it turns out that
there are various ways of embedding a locally compact Hausdorff
space into a  compact Hausdorff space.
The most economical construction consists
in  adding just one point.
This construction, known as the {\it Alexandroff compactification\/},
is technically useful, and we now describe it and sketch the proof
that it achieves its goal. 

\medskip
To help the reader's intuition, let us consider the case of the
plane $\reals^2$. If we view the plane $\reals^2$ as embedded
in $3$-space $\reals^3$, say as the $xOy$ plane of equation $z = 0$,
we can consider the sphere $\Sigma$ of radius $1$ centered on the $z$-axis
at the point $(0,0,1)$, and tangent to the $xOy$ plane at the origin
(sphere of equation $x^2 + y^2 + (z - 1)^2 = 1$).
If $N$ denotes the north pole on the sphere, i.e., the point 
of coordinates $(0,0,2)$, then any line $D$
passing through the north pole and not tangent to the sphere
(i.e., not parallel to the $xOy$ plane), intersects the $xOy$ plane
in a unique point $M$, 
and the sphere in a unique point $P$ other than the north pole $N$.
This, way, we obtain a bijection between the $xOy$ plane and the
punctured sphere $\Sigma$, i.e., the sphere with the north pole $N$ deleted.
This bijection is called a {\it stereographic projection\/}.
The Alexandroff compactification of the plane
consists in putting the north
pole back on the sphere, which amounts to adding a single point
at infinity $\infty$ to the plane. Intuitively, as we travel
away from the origin $O$ towards infinity 
(in any direction!), we tend towards
an ideal point at infinity $\infty$. Imagine that
we ``bend'' the plane so that it gets wrapped around the sphere,
according to stereographic projection. A simpler
example consists in taking a line and getting a circle
as its compactification.
The Alexandroff compactification is a generalization of these simple
constructions.

\begin{defin}
\label{Alexdef}
{\em
Let $(E, \s{O})$ be a locally compact  space.
Let $\omega$ be any point not in $E$, and let 
$E_{\omega} = E \cup \{\omega\}$.
Define the family $\s{O}_{\omega}$ as follows:
\[\s{O}_{\omega} = \s{O} \cup \{(E - K)\cup \{\omega\}\ |\
\hbox{$K$ compact in $E$}\}.\]
The pair $(E_{\omega}, \s{O}_{\omega})$ is called the
{\it Alexandroff compactification (or one point compactification)
of $(E, \s{O})$\/}.
}
\end{defin}

\medskip
The following theorem shows that $(E_{\omega}, \s{O}_{\omega})$ is 
indeed a topological space, and that it is compact.

\begin{thm}
\label{Alexthm}
Let $E$ be a locally compact topological  space.
The Alexandroff compactification $E_{\omega}$ of $E$ is a compact
space such that $E$ is a subspace of $E_{\omega}$,
and if $E$ is not compact, then $\overline{E} = E_{\omega}$.
\end{thm}

\proof The verification that $\s{O}_{\omega}$ is a family of 
open sets is not difficult but a bit tedious. Details can be
found in Munkres \cite{Munkrestop} or Schwartz \cite{Schwartz1}. 
Let us show that $E_{\omega}$ is 
compact. For every open cover $(U_i)_{i\in I}$ of $E_{\omega}$,
since $\omega$ must be covered, there is some $U_{i_{0}}$ of the form
\[U_{i_{0}} = (E - K_0)\cup \{\omega\}\]
where $K_0$ is compact in $E$. Consider the family $(V_i)_{i\in I}$
defined as follows:
\[\eqaligneno{
V_i &= U_i\quad\hbox{if}\quad U_i\in\s{O},\cr
V_i &= E - K\quad\hbox{if}\quad U_i = (E - K)\cup \{\omega\},\cr
}\]
where $K$ is compact in $E$.
Then, because each $K$ is compact and thus closed in $E$
(since $E$ is Hausdorff),
$E - K$ is open, and every $V_i$ is an open subset of $E$.
Furthermore, the family $(V_i)_{i\in (I - \{i_{0}\})}$
is an open cover of $K_0$. Since $K_0$ is compact, there is a finite
open subcover  $(V_j)_{j\in J}$ of $K_0$,
and thus,  $(U_j)_{j\in J \cup \{i_{0}\}}$ is a finite open
cover of $E_{\omega}$.

\medskip
Let us show that $E_{\omega}$ is Hausdorff. Given any two points
$a, b\in E_{\omega}$, if both $a, b\in E$, since $E$ is Hausdorff
and every open set in $\s{O}$ is an open set in $\s{O}_{\omega}$,
there exist disjoint open sets $U, V$ (in $\s{O}$)
such that $a\in U$ and $b\in V$.
If $b= \omega$, since $E$ is locally compact, there is some compact
set $K$ containing an open set $U$ containing $a$, and then,
$U$ and $V = (E - K)\cup \{\omega\}$ are disjoint open sets
(in $\s{O}_{\omega}$) such that $a\in U$ and $b\in V$.

\medskip
The space $E$ is a subspace of $E_{\omega}$ because
for every open set $U$ in $\s{O}_{\omega}$, either
$U\in \s{O}$ and $E\cap U = U$ is open in $E$, or
$U = (E - K)\cup \{\omega\}$ where $K$ is compact in $E$,
and thus $U\cap E = E - K$, which is open in $E$, since
$K$ is compact in $E$, and thus closed  (since $E$ is Hausdorff).
Finally, if $E$ is not compact, for every compact subset $K$ of $E$,
$E - K$ is nonempty, and thus, for every open set 
$U = (E - K)\cup \{\omega\}$ containing $\omega$, we have
$U\cap E \not= \emptyset$, which shows that $\omega\in \overline{E}$,
and thus that  $\overline{E} = E_{\omega}$.
$\square$

\medskip
Finally, in studying surfaces and manifolds, an important
property is the existence of a countable basis for the topology.
Indeed, this property guarantees the existence of trianguations
of surfaces, a crucial property.

\begin{defin}
\label{countbas}
{\em
A topological space is called {\it second-countable\/}
if there is a countable basis for its topology, i.e., if
there is a countable family $(U_i)_{i\geq 0}$ of open sets such that
every open set of $E$ is a union of open sets $U_i$.
}
\end{defin}

\medskip
It is easily seen that $\reals^n$ is second-countable,
and more generally, that every normed vector space of finite
dimension is second-countable. 
It can also be shown that if $E$ is a locally compact  space that
has a countable basis, then
$E_{\omega}$ also has a countable basis (and in fact, is metrizable).
We have the following properties.

\begin{prop}
\label{compac9}
Given a second-countable topological  space $E$, 
every open cover $(U_i)_{i\in I}$ of $E$ contains some countable subcover.
\end{prop}

\proof Let $(O_n)_{n\geq 0}$ be a countable basis for the
topology. Then, all sets $O_n$ contained in some $U_i$
can be arranged into a countable subsequence
$(\Omega_m)_{m\geq 0}$ of $(O_n)_{n\geq 0}$, and for every
$\Omega_m$, there is some $U_{i_{m}}$ such that
$\Omega_m \subseteq U_{i_{m}}$. Furthermore, every $U_i$ is
some union of sets $\Omega_j$, and thus, every $a\in E$
belongs to some  $\Omega_j$, which shows that
$(\Omega_m)_{m\geq 0}$ is a countable open subcover of
$(U_i)_{i\in I}$.
$\square$

\medskip
As an immediate corollary of Proposition \ref{compac9},
a locally connected second-countable space has countably many connected
components.

\medskip
In second-countable Hausdorff spaces, compactness can be characterized
in terms of accumulation points (this is also true of
metric spaces).

\begin{defin}
\label{accumdef}
{\em
Given a topological Hausdorff space $E$,
given any sequence $(x_n)$ of points in $E$,
a point $l\in E$ is an {\it accumulation point (or cluster point)\/}
of the sequence $(x_n)$ if every open set $U$ containing $l$
contains $x_n$ for infinitely many $n$.
}
\end{defin}

\medskip
Clearly, if $l$ is a limit of the sequence $(x_n)$, then
it is an accumulation point, since every open set $U$ containing
$a$ contains all $x_n$ except for finitely many $n$.

\begin{prop}
\label{compac10}
A second-countable topological Hausdorff  space $E$ is compact iff
every sequence $(x_n)$ has some accumulation point.
\end{prop}

\proof Assume that every sequence $(x_n)$ has 
some accumulation point. Let $(U_i)_{i\in I}$ be some
open cover of $E$. By Proposition \ref{compac9}, there is
a countable open subcover $(O_n)_{n\geq 0}$ for $E$.
Now, if $E$ is not covered by any finite subcover of
$(O_n)_{n\geq 0}$, we can define a sequence $(x_m)$ by induction
as follows:

\medskip
Let $x_0$ be arbitrary, and for every $m\geq 1$,
let $x_m$ be some point in $E$ not in
$O_1\cup \cdots\cup O_m$, which exists, since 
$O_1\cup \cdots\cup O_m$ is not an open cover of $E$.
We claim that the sequence $(x_m)$
does not have any accumulation point. Indeed, for every $l\in E$,
since  $(O_n)_{n\geq 0}$ is an open cover of $E$,
there is some $O_m$ such that $l\in O_m$, and by construction,
every $x_n$ with $n\geq m+1$ does not belong to $O_m$,
which means that $x_n\in O_m$ for only finitely many $n$, 
and $l$ is not an accumulation point.

\medskip
Conversely, assume that $E$ is compact, and let $(x_n)$
be any sequence. If $l\in E$ is not an accumulation point
of the sequence, then there is some open set $U_l$
such that $l\in U_l$ and $x_n\in U_l$ for only finitely many $n$.
Thus, if $(x_n)$ does not have any accumulation point,
the family $(U_l)_{l\in E}$ is an open cover of $E$,
and since $E$ is compact, it has some finite open subcover
$(U_l)_{l\in J}$, where $J$ is a finite subset of $E$. But every $U_l$ 
with $l\in J$ is such that  $x_n\in U_l$ for only finitely many $n$,
and since $J$ is finite, 
$x_n \in\bigcup_{l\in J} U_l$ for only finitely many $n$,
which contradicts the fact that $(U_l)_{l\in J}$ is an open cover of $E$,
and thus contains all the $x_n$.
Thus, $(x_n)$ has some accumulation point.
$\square$

\remark
It should be noted that the proof that
if $E$ is compact, then every sequence has some accumulation point,
holds for any arbitrary compact space (the proof does not use
a countable basis for the topology). The converse
also holds for metric spaces. We prove this converse
since it is a major property of metric spaces. It is also
convenient to have such a characterization of compactness
when dealing with fractal geometry.
\endremark

Given a metric space in which every sequence has some 
accumulation point,
we first prove the existence of a {\it Lebesgue number\/}.

\begin{lemma}
\label{BolzaWeir1}
Given a metric space $E$, if every sequence $(x_n)$
has an accumulation point, for every open cover
$(U_i)_{i\in I}$ of $E$, there is some $\delta > 0$
(a {\it Lebesgue number for $(U_i)_{i\in I}$\/}) such that,
for every open ball 
$\opball{a}{\epsilon}$ of diameter $\epsilon  \leq \delta$,
there is some open subset $U_i$ such that 
$\opball{a}{\epsilon}\subseteq U_i$.
\end{lemma}

\proof If there was no $\delta$ with the
above property, then, for every
natural number $n$, there would be some open ball
$\opball{a_n}{1/n}$ which is not contained in any open set $U_i$ of the
open cover $(U_i)_{i\in I}$. However, the sequence $(a_n)$ has
some accumulation point $a$, and since $(U_i)_{i\in I}$ is an open
cover of $E$, there is some $U_i$ such that $a \in U_i$.
Since $U_i$ is open, there is some open ball of 
center $a$ and radius $\epsilon$ contained in $U_i$.
Now, since $a$ is an accumulation point of the sequence $(a_n)$,
every open set containing $a$ contain $a_n$ for infinitely many $n$,
and thus, there is some $n$ large enough so that 
\[1/n \leq \epsilon/2\quad\hbox{and}\quad a_n \in \opball{a}{\epsilon/2},\]
which implies that 
\[\opball{a_n}{1/n}\subseteq \opball{a}{\epsilon} \subseteq U_i,\]
a contradiction. 
$\square$

\medskip
By a previous remark, since the proof of Proposition \ref{compac10}
implies that in a compact topological space, every
sequence has some accumulation point, by Lemma \ref{BolzaWeir1},
in a compact metric space, every open cover has a Lebesgue number.
This fact can be used to prove another important property of
compact metric spaces, the uniform continuity theorem.

\begin{defin}
\label{unifcont}
{\em
Given two metric spaces $(E, d_E)$ and $(F, d_F)$,
a function $\mapdef{f}{E}{F}$ is {\it uniformly continuous\/}
if for every $\epsilon > 0$, there is some $\eta > 0$, such that,
for all $a, b \in E$, 
\[\hbox{if}\quad d_E(a, b) \leq \eta\quad\hbox{then}\quad
d_F(f(a), f(b)) \leq \epsilon.\]
}
\end{defin}

\medskip
The {\it uniform continuity theorem\/} can be stated as follows.

\begin{thm}
\label{unifcontlem}
Given two metric spaces $(E, d_E)$ and $(F, d_F)$, if $E$ is compact
and $\mapdef{f}{E}{F}$ is a continuous function, then it is
uniformly continuous. 
\end{thm}

\proof Consider any $\epsilon > 0$, and 
let $(\opball{y}{\epsilon/2})_{y\in F}$ be the open cover of $F$
consisting of open balls of radius $\epsilon/2$. 
Since $f$ is continuous, the family
\[(f^{-1}(\opball{y}{\epsilon/2}))_{y\in F}\]
is an open cover of $E$. Since, $E$ is compact, by Lemma 
\ref{BolzaWeir1}, there is a Lebesgue number $\delta$ such that
for every open ball $\opball{a}{\eta}$ of diameter $\eta  \leq \delta$,
then  $\opball{a}{\eta}\subseteq f^{-1}(\opball{y}{\epsilon/2})$, for
some $y\in F$. In particular, for any $a, b\in E$ such that
$d_E(a, b) \leq \eta = \delta/2$,
we have $a, b\in \opball{a}{\delta}$, and thus
$a, b\in  f^{-1}(\opball{y}{\epsilon/2})$, which implies that
$f(a), f(b) \in \opball{y}{\epsilon/2}$. But then,
$d_F(f(a), f(b)) \leq \epsilon$, as desired.
$\square$

\medskip
We now prove another lemma needed to obtain
the characterization of compactness in metric spaces in terms
of accumulation points.

\begin{lemma}
\label{BolzaWeir2}
Given a metric space $E$, if every sequence $(x_n)$
has an accumulation point, then for every $\epsilon > 0$,
there is a finite open cover 
$\opball{a_0}{\epsilon}\cup\cdots\cup \opball{a_n}{\epsilon}$
of $E$ by open balls of radius $\epsilon$.
\end{lemma}

\proof Let $a_0$ be any point in $E$. If 
$\opball{a_0}{\epsilon} = E$, the lemma is proved. Otherwise,
assume that a sequence $(a_0,a_1,\ldots,a_n)$ has been defined,
such that 
$\opball{a_0}{\epsilon}\cup\cdots\cup \opball{a_n}{\epsilon}$
does not cover $E$. Then, there is some 
$a_{n+1}$ not in 
$\opball{a_0}{\epsilon}\cup\cdots\cup \opball{a_n}{\epsilon}$,
and either  
\[\opball{a_0}{\epsilon}\cup\cdots\cup \opball{a_{n + 1}}{\epsilon} = E,\]
in which case the lemma is proved, or we obtain a sequence
$(a_0,a_1,\ldots,a_{n+1})$ 
such that $\opball{a_0}{\epsilon}\cup\cdots\cup \opball{a_{n+1}}{\epsilon}$
does not cover $E$. If this process goes on forever, we obtain an infinite
sequence $(a_n)$ such that $d(a_m, a_n) > \epsilon$ for all
$m\not= n$.
Since every sequence in $E$ has some accumulation point,
the sequence $(a_n)$ has some accumulation point $a$.
Then, for infinitely many $n$, we must have 
$d(a_n, a) \leq \epsilon/3$, and thus, for at least
two distinct natural numbers $p, q$, we must have
$d(a_p, a) \leq \epsilon/3$ and $d(a_q, a) \leq \epsilon/3$,
which implies $d(a_p, a_q) \leq 2\epsilon/3$, contradicting
the fact that  $d(a_m, a_n) > \epsilon$ for all
$m\not= n$.
Thus, there must be some $n$ such that
\[\opball{a_0}{\epsilon}\cup\cdots\cup \opball{a_n}{\epsilon} = E.\]
$\square$

\medskip
A metric space satisfying the condition of Lemma \ref{BolzaWeir2}
is sometimes called {\it precompact\/} (or {\it totally bounded\/}).
We now obtain the {\it Weierstrass-Bolzano\/} property.

\begin{thm}
\label{BolzaWeirthm}
A metric space $E$ is compact iff
every sequence $(x_n)$ has an accumulation point.
\end{thm}

\proof We already observed that the proof of 
Proposition \ref{compac10} shows that for any compact space
(not necessarily metric), 
every sequence $(x_n)$ has an  accumulation point.
Conversely, let $E$ be a  metric space, and 
assume that every sequence $(x_n)$ has an accumulation point.
Given any open cover $(U_i)_{i\in I}$ for $E$,
we must find a finite open subcover of $E$.
By Lemma \ref{BolzaWeir1}, there is some $\delta > 0$
(a Lebesgue number for $(U_i)_{i\in I}$) such that,
for every open ball $\opball{a}{\epsilon}$ of diameter $\epsilon  \leq \delta$,
there is some open subset $U_j$ such that 
$\opball{a}{\epsilon}\subseteq U_j$. By Lemma \ref{BolzaWeir2},
for every $\delta > 0$, there is a finite open cover 
$\opball{a_0}{\delta}\cup\cdots\cup \opball{a_n}{\delta}$ 
of $E$ by open balls of
radius $\delta$. But from the previous statement,
every open ball $\opball{a_i}{\delta}$ is contained
in some open set $U_{j_{i}}$, and thus,
$\{U_{j_{1}},\ldots, U_{j_{n}}\}$ is an open cover of $E$.
$\square$

\medskip
Another very useful characterization of compact metric
spaces is obtained in terms of Cauchy sequences.
Such a characterization is quite useful in fractal geometry
(and elsewhere).
First, recall the definition of a Cauchy sequence, and
of a complete metric space.

\begin{defin}
\label{Cauchydef2}
{\rm
Given a metric space $(E, d)$, a sequence $(x_n)_{n\in\natnums}$
in $E$ is a {\it Cauchy sequence\/} if the following condition holds: 

\medskip
for every $\epsilon > 0$,
there is some $p\geq 0$, such that, for all $m, n\geq p$,
then $\dist{}{x_m}{x_n} \leq \epsilon$.

\medskip\
If every Cauchy sequence in $(E, d)$ converges, we say that
$(E, d)$ is a {\it complete metric space\/}.
}
\end{defin}

\medskip
First, let us show the following easy proposition.

\begin{prop}
\label{compl1}
Given a metric space $E$, if a Cauchy sequence $(x_n)$ 
has some accumulation point $a$, then $a$ is the limit
of the sequence $(x_n)$.
\end{prop}

\proof Since $(x_n)$ is a Cauchy sequence,
for every $\epsilon > 0$,
there is some $p\geq 0$, such that, for all $m, n\geq p$,
then $\dist{}{x_m}{x_n} \leq \epsilon/2$.
Since $a$ is an accumulation point for $(x_n)$, 
for infinitely many $n$, we have $d(x_n, a) \leq \epsilon/2$,
and thus for at least some $n \geq p$,
have $d(x_n, a) \leq \epsilon/2$. Then, for all $m\geq p$,
\[d(x_m, a)\leq  d(x_m, x_n) + d(x_n, a) \leq \epsilon,\]
which shows that $a$ is the limit of the sequence $(x_n)$.
$\square$

\medskip
Recall that a metric space is {\it precompact\/} (or {\it totally
bounded\/}) if for every $\epsilon > 0$,
there is a finite open cover 
$\opball{a_0}{\epsilon}\cup\cdots\cup \opball{a_n}{\epsilon}$
of $E$ by open balls of radius $\epsilon$. 
We can now prove the following theorem.

\begin{thm}
\label{complthm}
A metric space $E$ is compact iff it is precompact and complete.
\end{thm}

\proof Let $E$ be compact. For every $\epsilon > 0$,
the family of all open balls of radius $\epsilon$ is an open cover for $E$,
and since $E$ is compact, there is a finite subcover
$\opball{a_0}{\epsilon}\cup\cdots\cup \opball{a_n}{\epsilon}$
of $E$ by open balls of radius $\epsilon$.
Thus, $E$ is precompact. 
Since $E$ is compact, by Theorem \ref{BolzaWeirthm},
every sequence $(x_n)$ has some accumulation point.
Thus, every Cauchy sequence $(x_n)$ has some
accumulation point $a$, and by Proposition \ref{compl1},
$a$ is the limit of $(x_n)$. Thus, $E$ is complete.

\medskip
Now, assume that $E$ is precompact and complete.
We prove that every sequence $(x_n)$ has an accumulation point.
By the other direction of Theorem \ref{BolzaWeirthm}, this shows
that $E$ is compact.
Given any sequence $(x_n)$, we construct a Cauchy subsequence
$(y_n)$ of $(x_n)$ as follows:
Since $E$ is precompact, letting $\epsilon = 1$, there
exists a finite cover $\s{U}_{1}$ of $E$ by open balls of
radius $1$. Thus, some open ball $B_{o}^{1}$ in the cover $\s{U}_{1}$
contains infinitely many elements from the sequence $(x_n)$.
Let $y_0$ be any element of $(x_n)$ in $B_{o}^{1}$.
By induction, assume that a sequence of open balls 
$(B_{o}^{i})_{1\leq i \leq m}$
has been defined, such that every ball $B_{o}^{i}$
has radius $\frac{1}{2^{i}}$, contains infinitely many
elements from the sequence $(x_n)$, and
contains some $y_i$ from $(x_n)$ such that
\[d(y_i, y_{i+1}) \leq \frac{1}{2^{i}},\]
for all $i$, $0\leq i \leq m-1$.
Then, letting $\epsilon = \frac{1}{2^{m+1}}$,
because $E$ is precompact, there is some 
finite cover $\s{U}_{m+1}$ of $E$ by open balls 
of radius $\epsilon$, and thus of the open ball  $B_{o}^{m}$.
Thus, some open ball $B_{o}^{m+1}$ in the cover
$\s{U}_{m+1}$ contains infinitely many elements from the sequence $(x_n)$,
and we let $y_{m+1}$ be any element of $(x_n)$ in $B_{o}^{m+1}$.
Thus, we have defined by induction a sequence $(y_n)$
which is a subsequence of $(x_n)$, and  such that
\[d(y_i, y_{i+1}) \leq \frac{1}{2^{i}},\]
for all $i$. However, for all $m, n\geq 1$, we have
\[d(y_{m}, y_{n}) \leq d(y_{m}, y_{m+1}) + \cdots +  d(y_{n-1}, y_{n})
\leq \sum_{i = m}{n} \frac{1}{2^{i}} \leq \frac{1}{2^{m-1}},\]
and thus, $(y_n)$ is a Cauchy sequence.
Since $E$ is complete, the sequence $(y_n)$ has a limit,
and since it is a subsequence of $(x_n)$, the sequence 
$(x_n)$ has some accumulation point.
$\square$

\medskip
If $(E, d)$ is a nonempty complete metric space, every map
$\mapdef{f}{E}{E}$ for which there is some $k$ such that $0 \leq  k < 1$ and
\[d(f(x), f(y)) \leq k d(x, y)\]
for all $x, y\in E$, has the very important property that
it has a unique fixed point, that is, there is a unique $a\in E$
such that $f(a) = a$. A map as above is called a {\it contracting
mapping\/}. Furthermore, the fixed point of a contracting mapping
can be computed as the limit of a fast converging sequence.

\medskip
The fixed point property of contracting mappings is used
to show some important theorems of analysis, such as the implicit
function theorem, and the existence of solutions to certain differential
equations. It can also be used to show the existence of  fractal sets
defined in terms of iterated function systems,
a topic that we intend to discuss later on. Since the proof is
quite simple, we prove the fixed point property of contracting 
mappings. First, observe that a contracting mapping is (uniformly)
continuous.

\begin{prop}
\label{fixpt}
If $(E, d)$ is a nonempty complete metric space, every
contracting mapping $\mapdef{f}{E}{E}$ has a unique fixed
point. Furthermore, for every $x_0\in E$,
defining the sequence $(x_n)$ such that $x_{n+1} = f(x_{n})$,
the sequence $(x_n)$ converges to the unique fixed point of $f$.
\end{prop}

\proof First, we prove that $f$ has at most one fixed point.
Indeed, if $f(a) = a$ and $f(b) = b$, since
\[d(a, b) = d(f(a), f(b)) \leq k d(a, b)\]
and $0 \leq  k < 1$, we must have $d(a, b) = 0$, that is, $a = b$.

\medskip
Next, we prove that $(x_n)$ is a Cauchy sequence.
Observe that
\[\eqaligneno{
d(x_2, x_1) &\leq k d(x_1, x_0),\cr
d(x_3, x_2) &\leq k d(x_2, x_1) \leq k^2 d(x_1, x_0),\cr
\cdots      &\quad \cdots\cr
d(x_{n+1}, x_n) &\leq k d(x_n, x_{n-1}) \leq \cdots\leq k^n d(x_1, x_0).\cr
}\]
Thus, we have
\[\eqaligneno{
d(x_{n + p}, x_{n}) &\leq d(x_{n + p}, x_{n + p -1}) + 
d(x_{n + p - 1}, x_{n + p - 2}) + \cdots + d(x_{n + 1}, x_{n})\cr
&\leq (k^{p-1} + k^{p-2} + \cdots + k + 1) k^n d(x_1, x_0)\cr
&\leq \frac{k^n}{1 - k}\, d(x_1, x_0).\cr
}\]  
We conclude that $d(x_{n + p}, x_{n})$ converges to $0$ when
$n$ goes to infinity, which shows that $(x_n)$ is a Cauchy
sequence. Since $E$ is complete, the sequence $(x_n)$ has a limit
$a$. Since $f$ is continuous, the sequence $(f(x_n))$ converges
to $f(a)$. But $x_{n+1} = f(x_n)$ converges to $a$, and
so $f(a) = a$, the unique fixed point of $f$.
$\square$

\medskip
Note that no matter how the starting point $x_0$ of the sequence $(x_n)$
is chosen, $(x_n)$ converges to the unique fixed point of $f$.
Also, the convergence is fast, since 
\[d(x_n, a) \leq \frac{k^n}{1 - k}\, d(x_1, x_0).\]

\medskip
The Hausdorff distance between compact subsets of a metric space
provides a very nice illustration of some of the theorems
on complete and compact metric spaces just presented. It can also
be used to define certain kinds of fractal sets, and thus,
we indulge into a short digression on the Hausdorff distance.

\begin{defin}
\label{Hausdist}
{\em
Given a metric space $(X, d)$, for any  subset $A\subseteq X$,
for any $\epsilon \geq 0$, define the {\it $\epsilon$-hull\/} of
$A$, as the set 
\[V_{\epsilon}(A) = \{x\in X,\> \exists a\in A |\ d(a,x) \leq \epsilon\}.\] 
Given any two nonempty bounded subsets $A, B$ of $X$, define
$D(A, B)$, {\it the Hausdorff distance between $A$ and $B$\/}, as
\[D(A, B) = \inf\{\epsilon\geq 0\ |\ A\subseteq V_{\epsilon}(B)\ \hbox{and}\
B\subseteq V_{\epsilon}(A)\}.\]
}
\end{defin}

\medskip
Note that since we are considering nonempty bounded subsets,
$D(A, B)$ is well defined (i.e., not infinite). However,
$D$ is not necessarily a distance function. It is 
a distance function if we restrict
our attention to nonempty compact subsets of $X$.
We let $\s{K}(X)$ denote the set of all nonempty compact subsets
of $X$. The remarkable fact is that $D$ is a distance on $\s{K}(X)$,
and that if $X$ is complete or compact, then so it  $\s{K}(X)$.
The following theorem is taken from Edgar \cite{Edgar}.

\begin{thm}
\label{Hausdthm}
If $(X, d)$ is a metric space, then the Hausdorff distance $D$
on the set $\s{K}(X)$ of nonempty compact subsets of $X$ is a distance.
If $(X, d)$ is complete, then $(\s{K}(X), D)$ is complete, 
and if $(X, d)$ is compact, then $(\s{K}(X), D)$ is compact.
\end{thm}

\proof Since (nonempty) compact sets are bounded,
$D(A, B)$ is well defined. Clearly, $D$ is symmetric. 
Assume that $D(A, B) = 0$. Then, for every $\epsilon > 0$,
$A\subseteq V_{\epsilon}(B)$, which means that for every $a\in A$,
there is some $b\in B$ such that $d(a, b) \leq \epsilon$, and thus,
that $A\subseteq \overline{B}$. Since $B$ is closed,
$\overline{B} = B$, and we have $A\subseteq B$.
Similarly,  $B\subseteq A$, and thus, $A = B$.
Clearly, if $A = B$, we have $D(A, B) = 0$.
It remains to prove the triangle inequality.
If $B \subseteq V_{\epsilon_{1}}(A)$ and 
$C \subseteq V_{\epsilon_{2}}(B)$, then
\[V_{\epsilon_{2}}(B) \subseteq V_{\epsilon_{2}}(V_{\epsilon_{1}}(A)),\] 
and since
\[V_{\epsilon_{2}}(V_{\epsilon_{1}}(A))\subseteq 
V_{\epsilon_{1} + \epsilon_{2}}(A),\] 
we get
\[C \subseteq V_{\epsilon_{2}}(B) 
\subseteq V_{\epsilon_{1} + \epsilon_{2}}(A).\] 
Similarly, we can prove that
\[A \subseteq V_{\epsilon_{1} + \epsilon_{2}}(C),\] 
and thus, the triangle inequality follows.

\medskip
Next, we need to prove that if 
$(X, d)$ is complete, then $(\s{K}(X), D)$ is also complete. 
First, we show that if $(A_n)$ is a sequence of nonempty compact sets
converging to a nonempty compact set
$A$ in the Hausdorff metric, then
\[A = \{x\in X\ |\ \hbox{there is a sequence $(x_n)$ with $x_n\in A_n$
converging to $x$}\}.\]
Indeed, if $(x_n)$ is a sequence  with $x_n\in A_n$ converging to
$x$ and  $(A_n)$ converges to $A$,
then for every $\epsilon > 0$, there is some $x_n$ such that
$d(x_n, x) \leq \epsilon/2$, 
and there is some $a_n\in A$ such that $d(a_n, x_n) \leq \epsilon/2$,
and thus, $d(a_n, x)\leq \epsilon$, which shows that
$x\in \overline{A}$. Since $A$ is compact, it is closed, and
$x\in A$. Conversely, since $(A_n)$ converges to $A$,
for every $x\in A$, for every $n\geq 1$, there is some
$x_n\in A_n$ such that $d(x_n, x) \leq 1/n$,
and the sequence $(x_n)$ converges to $x$.

\medskip
Now, let $(A_n)$ be a Cauchy sequence in $\s{K}(X)$.
It can be proven that $(A_n)$ converges to the set
\[A = \{x\in X\ |\ \hbox{there is a sequence $(x_n)$ with $x_n\in A_n$
converging to $x$}\},\] 
and that $A$ is nonempty and compact.
To prove that $A$ is compact, one proves that it is totally bounded
and complete. Details are given in Edgar \cite{Edgar}.

\medskip
Finally, we need to prove that 
if $(X, d)$ is compact, then $(\s{K}(X), D)$ is compact.
Since we already know that $(\s{K}(X), D)$ is complete if
$(X, d)$ is, it is enough to prove that
$(\s{K}(X), D)$ is totally bounded  if $(X, d)$ is, which
is fairly easy. 
$\square$

\medskip
In view of Theorem \ref{Hausdthm} and Theorem \ref{fixpt},
it is possible to define some nonempty compact subsets
of $X$ in terms of fixed points of contracting maps.
We will see later on how this can be done in terms of iterated function
systems, yielding a large class of fractals.

\medskip
Finally, returning to second-countable spaces,
we give another characterization of accumulation 
points.

\begin{prop}
\label{compac11}
Given a second-countable  topological Hausdorff space $E$, 
a point $l$ is an accumulation point of the sequence $(x_n)$
iff $l$ is the limit of some subsequence $(x_{n_{k}})$ of $(x_n)$.
\end{prop}

\proof
Clearly, if $l$ is the limit of some subsequence $(x_{n_{k}})$ of $(x_n)$,
it is an accumulation point of $(x_n)$.

\medskip
Conversely, let $(U_k)_{k\geq 0}$ be the sequence of
open sets containing $l$, where each $U_k$ 
belongs to a countable basis of $E$, and let
$V_k = U_1\cap \cdots \cap U_k$. For every $k\geq 1$, we can
find some $n_{k}> n_{k-1}$ such that $x_{n_{k}}\in V_k$,
since $l$ is an accumulation point of $(x_n)$.
Now, since every open set containing $l$ contains some $U_{k_{0}}$, and
since $x_{n_{k}}\in U_{k_{0}}$ for all $k\geq 0$,
the sequence $(x_{n_{k}})$ has limit $l$.
$\square$

\remark 
Proposition \ref{compac11} also holds for metric spaces.
\endremark

As promised, we show how certain fractals can be defined
by iterated function systems, using Theorem \ref{Hausdthm} 
and Theorem \ref{fixpt}.

\chapter{A Detour On Fractals}
\label{chap2}
\section{Iterated Function Systems and Fractals}
\label{fracsec}
A pleasant application of the Hausdorff distance and of the fixed point
theorem for contracting mappings is a method for defining
a class of ``self-similar'' fractals. For this, we can
use iterated function systems.

\begin{defin}
\label{iresys}
{\em
Given a metric space $(X, d)$, 
an {\it iterated function system\/}, for short, an {\it ifs\/},
is a finite sequence of functions $(f_1,\ldots, f_n)$,
where each $\mapdef{f_i}{X}{X}$ is a contracting mapping.
A nonempty compact subset $K$ of $X$ is an {\it invariant set
(or attractor)\/}
for the ifs   $(f_1,\ldots, f_n)$ if
\[K = f_1(K)\cup \cdots \cup f_n(K).\]
}
\end{defin}

\medskip
The major result about ifs's is the following.

\begin{thm}
\label{fracalhm}
If $(X, d)$ is a nonempty complete metric space,
every iterated function system   $(f_1,\ldots, f_n)$
has a unique invariant set $A$ which is a nonempty compact
subset of $X$. Furthermore, for every nonempty compact subset
$A_0$ of $X$, this invariant
set $A$ if the limit of the sequence $(A_m)$, where
$A_{m+1} = f_1(A_m)\cup \cdots \cup f_n(A_m)$.
\end{thm}

\proof Since $X$ is complete, by Theorem \ref{Hausdthm}, the space
$(\s{K}(X), D)$ is a complete metric space.
The theorem will follow from Theorem \ref{fixpt},
if we can show that the map
$\mapdef{F}{\s{K}(X)}{\s{K}(X)}$ defined such that
\[F(K)  = f_1(K)\cup \cdots \cup f_n(K),\]
for every nonempty compact set $K$, is a contracting mapping.
Let $A, B$ be any two nonempty compact subsets of $X$, and
consider any $\eta \geq D(A, B)$. Since each
$\mapdef{f_i}{X}{X}$ is a contracting mapping, there is some $\lambda_i$,
with $0\leq \lambda_i < 1$, such that
\[d(f_i(a), f_i(b)) \leq \lambda_i d(a, b),\]
for all $a, b\in X$. Let $\lambda = \max\{\lambda_1, \ldots, \lambda_n\}$.
We claim that 
\[D(F(A), F(B)) \leq \lambda D(A, B).\]
For any $x\in F(A) = f_1(A)\cup \cdots \cup f_n(A)$,
there is some $a_i\in A_i$ such that $x = f_i(a_i)$,
and since $\eta = D(A, B)$, there is some $b_i\in B$ such that
\[d(a_i, b_i) \leq \eta,\]
and thus,
\[d(x, f_i(b_i)) = d(f_i(a_i), f_i(b_i)) \leq \lambda_i d(a_i, b_i)\leq
\lambda\eta.\]
This show that 
\[F(A) \subseteq V_{\lambda\eta}(F(B)).\]
Similarly, we can prove that
\[F(B) \subseteq V_{\lambda\eta}(F(A)),\]
and since this holds for all $\eta\geq D(A, B)$,
we proved that 
\[D(F(A), F(B)) \leq \lambda D(A, B)\]
where $\lambda = \max\{\lambda_1, \ldots, \lambda_n\}$.
Since $0\leq \lambda_i < 1$, we have $0\leq \lambda < 1$, and
$F$ is indeed a contracting mapping.
$\square$

\medskip
Theorem \ref{fracalhm} justifies the existence of many
familiar ``self-similar'' fractals. One of the best known fractals is 
the {\it Sierpinski gasket\/}.

\begin{example}
\label{ex1}
{\em
Consider an equilateral triangle with vertices $a, b, c$,
and let $f_1, f_2, f_3$ be the dilatations of centers $a, b, c$
and ratio $1/2$. The Sierpinski gasket is the invariant set
of the ifs $(f_1, f_2, f_3)$.
The dilations $f_1, f_2, f_3$ can be defined explicitly 
as follows, assuming that $a = (-1/2, 0)$, $b = (1/2, 0)$,
and $c = (0, \sqrt{3}/2)$. The contractions $f_a, f_b, f_c$ are specified by
\[\eqaligneno{
x' &= \frac{1}{2}x - \frac{1}{4},\cr
y' &= \frac{1}{2}y,\cr
}\]
\[\eqaligneno{
x' &= \frac{1}{2}x + \frac{1}{4},\cr
y' &= \frac{1}{2}y,\cr
}\]
and
\[\eqaligneno{
x' &= \frac{1}{2}x,\cr
y' &= \frac{1}{2}y + \frac{\sqrt{3}}{4}.\cr
}\]

We wrote a {\it Mathematica\/} program that iterates any finite
number of affine maps on any input figure consisting of
combinations of points, line segments, and polygons
(with their interior points).
Starting with the edges of the triangle $a, b, c$,  after $6$ iterations,
we get the  picture shown in Figure \ref{gasket1fig}.

\begin{figure}
\centerline{
\psfig{figure=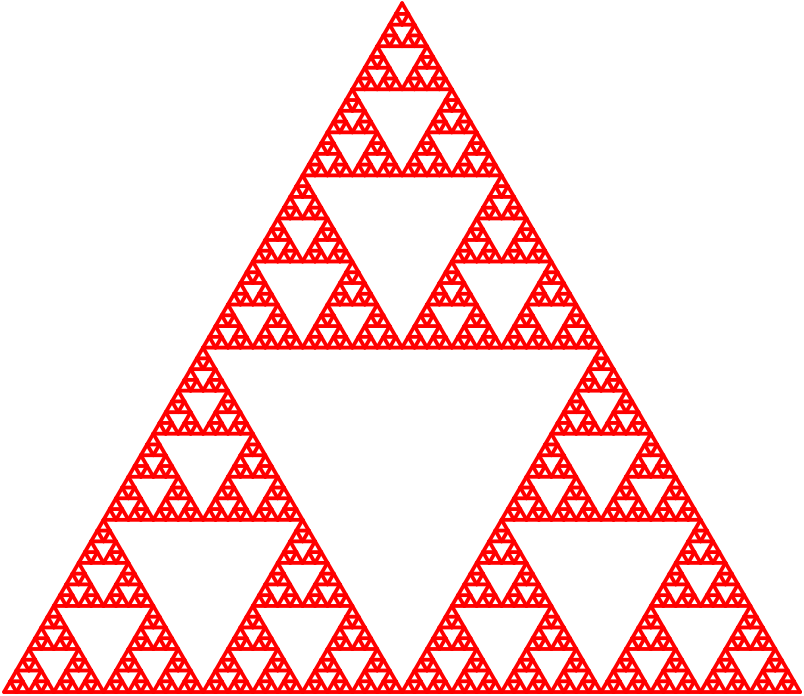,width=5truein}}
\caption{The Sierpinski gasket}
\label{gasket1fig}
\end{figure}

It is amusing that the same fractal is obtained no matter what
the initial nonempty compact figure is. It is interesting to see
what happens if we start with a solid triangle (with its
interior points). The result after $6$ iterations
is shown in Figure \ref{gasket2fig}.

\medskip
\begin{figure}
\centerline{
\psfig{figure=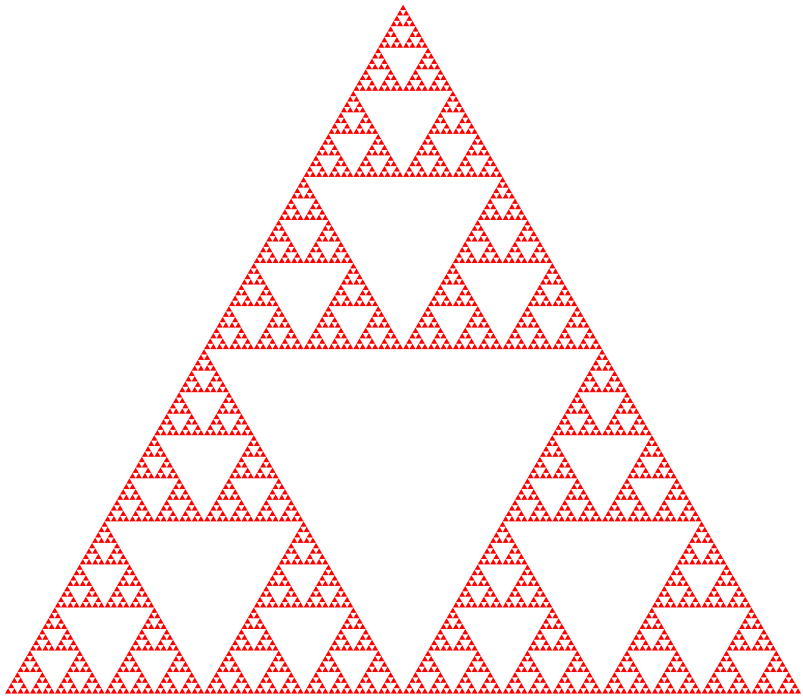,width=5truein}}
\caption{The Sierpinski gasket, version 2}
\label{gasket2fig}
\end{figure}

The convergence towards the Sierpinski gasket is very fast.
Incidently, there are many other ways of defining the
Sierpinski gasket. 
}
\end{example}

\medskip
A nice variation on the theme of the Sierpinski gasket is the
{\it Sierpinski dragon\/}. 

\begin{example}
{\em
The Sierpinski dragon  is specified by the following three
contractions:
\[\eqaligneno{
x' &= -\frac{1}{4}x - \frac{\sqrt{3}}{4}y + \frac{3}{4},\cr 
y' &= \frac{\sqrt{3}}{4}x - \frac{1}{4}y + \frac{\sqrt{3}}{4},\cr
}\]
\[\eqaligneno{
x' &= -\frac{1}{4}x + \frac{\sqrt{3}}{4}y - \frac{3}{4},\cr 
y' &= -\frac{\sqrt{3}}{4}x - \frac{1}{4}y +  \frac{\sqrt{3}}{4},\cr
}\]
\[\eqaligneno{
x' &= \frac{1}{2}x,\cr
y' &=  \frac{1}{2}y +  \frac{\sqrt{3}}{2}.\cr
}\]

The result of $7$ iterations starting from the line
segment $(-1, 0), (1, 0))$, is shown in Figure \ref{Sierpdragfig}.

\begin{figure}
\centerline{
\psfig{figure=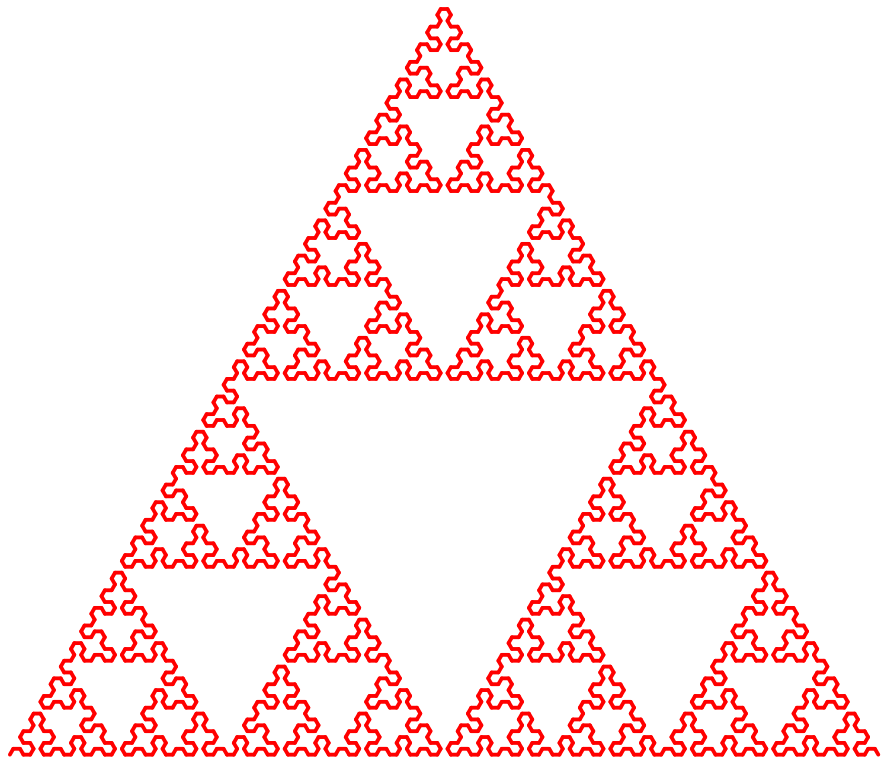,width=5truein}}
\caption{The Sierpinski dragon}
\label{Sierpdragfig}
\end{figure}

This curve converges to the boundary of the Sierpinski gasket.
}
\end{example}

A different kind of fractal
is the {\it Heighway dragon\/}. 
\begin{example}
{\em
The Heighway dragon is specified by the following two contractions:
\[\eqaligneno{
x' &= \frac{1}{2}x - \frac{1}{2}y,\cr
y' &= \frac{1}{2}x +  \frac{1}{2}y,\cr
}\]
\[\eqaligneno{
x' &= -\frac{1}{2}x  -\frac{1}{2}y,\cr
y' &= \frac{1}{2}x -\frac{1}{2}y + 1.\cr
}\]

It can be shown that for any number of iterations,
the polygon does not cross itself. This means that
no edge is traversed twice, and that if a point is traversed
twice, then this point is the endpoint of some edge.
The result of $13$ iterations, starting with the line segment
$((0, 0), (0, 1))$, is shown in Figure \ref{Highwayfig}.

\begin{figure}
\centerline{
\psfig{figure=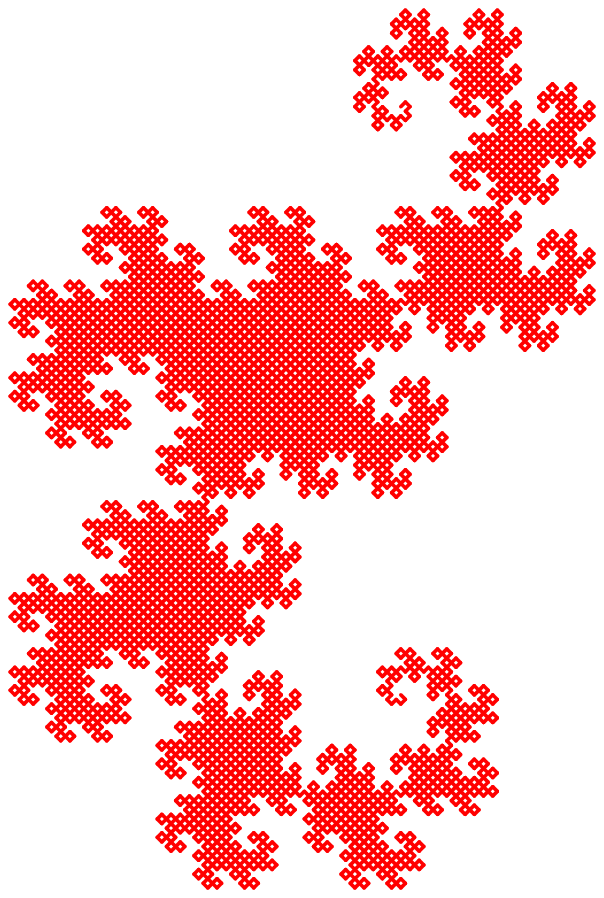,width=5truein}}
\caption{The Heighway dragon}
\label{Highwayfig}
\end{figure}

The Heighway dragon turns out to fill a closed and bounded set.
It can also be shown that the plane can be tiled with
copies of the Heighway dragon.
}
\end{example}

Another well known example is the {\it Koch curve\/}. 

\begin{example}
{\em
The Koch curve is specified by
the following four contractions:
\[\eqaligneno{
x' &= \frac{1}{3}x - \frac{2}{3},\cr 
y' &= \frac{1}{3}y,\cr
}\]
\[\eqaligneno{
x' &= \frac{1}{6}x - \frac{\sqrt{3}}{6}y - \frac{1}{6},\cr 
y' &= \frac{\sqrt{3}}{6}x + \frac{1}{6}y +  \frac{\sqrt{3}}{6},\cr
}\] 
\[\eqaligneno{
x' &= \frac{1}{6}x + \frac{\sqrt{3}}{6}y + \frac{1}{6},\cr 
y' &= -\frac{\sqrt{3}}{6}x +  \frac{1}{6}y + \frac{\sqrt{3}}{6},\cr
}\] 
\[\eqaligneno{
x' &= \frac{1}{3}x + \frac{2}{3},\cr 
y' &= \frac{1}{3}y,\cr
}\] 

The Koch curve is an example of a continuous curve which
is nowhere differentiable (because it ``wiggles'' too much).
It is a curve of infinite length.
The result of $6$ iterations, starting with the line segment
$((-1, 0), (1, 0))$, is shown in Figure \ref{Kochfig}.

\begin{figure}
\centerline{
\psfig{figure=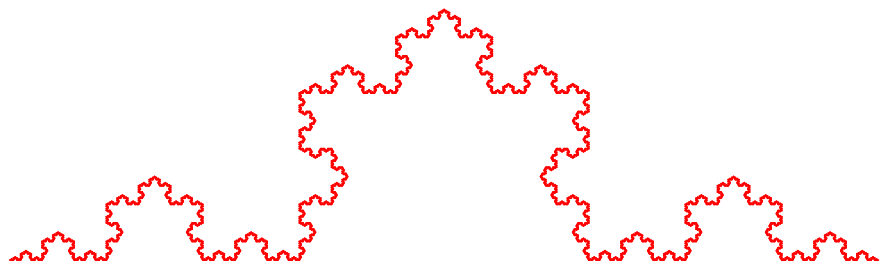,width=5truein}}
\caption{The Koch curve}
\label{Kochfig}
\end{figure}
}
\end{example}

The curve obtained by putting three Kock curves together on
the sides of an equilateral triangle is known as the 
{\it snowflake curve\/} (for obvious reasons, see below!).

\begin{example}
{\em
The snowflake curve obtained after $5$ iterations is shown in Figure
\ref{snowflakefig}.

\begin{figure}
\centerline{
\psfig{figure=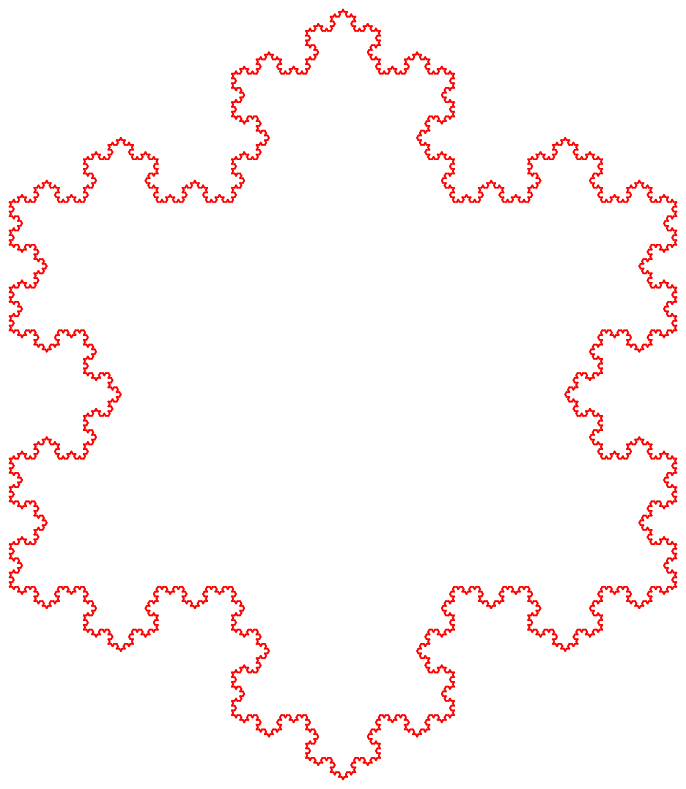,width=6.5truein}}
\caption{The snowflake curve}
\label{snowflakefig}
\end{figure}

The snowflake curve is an example of a closed curve of infinite
length bounding a finite area.
}
\end{example}

We conclude with another famous example, a variant of the
{\it Hilbert curve\/}. 

\begin{example}
{\em
This version of the  Hilbert curve is defined by the following
four contractions:
\[\eqaligneno{
x' &= \frac{1}{2}x - \frac{1}{2},\cr 
y' &= \frac{1}{2}y + 1,\cr
}\] 
\[\eqaligneno{
x' &= \frac{1}{2}x + \frac{1}{2},\cr 
y' &= \frac{1}{2}y + 1,\cr
}\] 
\[\eqaligneno{
x' &= -\frac{1}{2}y + 1,\cr 
y' &= \frac{1}{2}x + \frac{1}{2},\cr 
}\]
\[\eqaligneno{
x' &= \frac{1}{2}y - 1,\cr 
y' &= -\frac{1}{2}x + \frac{1}{2},\cr
}\]  

This continuous curve is a space-filling curve,
in the sense that its image is the entire unit square.
The result of $6$ iterations, starting with the two lines
segments $((-1, 0), (0, 1))$ and  $((0, 1), (1, 0))$, is
shown in Figure \ref{hilcurvfig}.

\begin{figure}
\centerline{
\psfig{figure=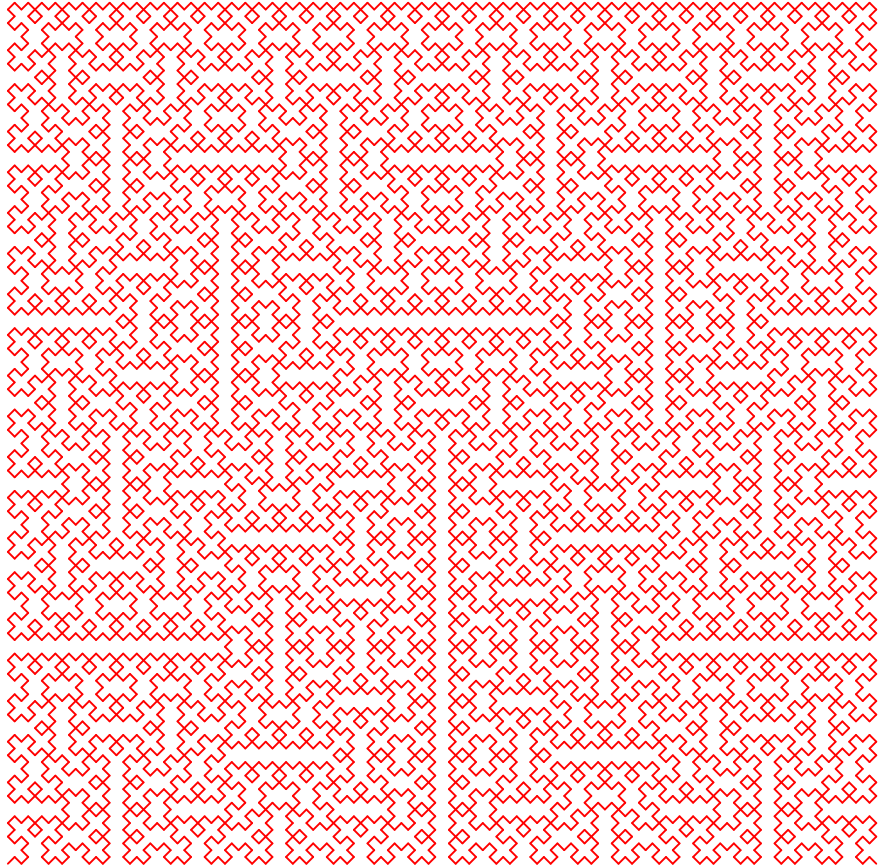,width=5truein}}
\caption{A Hilbert curve}
\label{hilcurvfig}
\end{figure}
}
\end{example}

\medskip
For more on iterated function systems and
fractals, we recommend Edgar \cite{Edgar}.


\bibliography{../basicmath/cadgeom}

\begin{thebibliography}{10}

\bibitem{Ahlfors}
Lars~V. Ahlfors and Leo Sario.
\newblock {\em Riemann Surfaces}.
\newblock Princeton Math. Series, No. 2. Princeton University Press, 1960.

\bibitem{Amstrong}
Mark~A. Amstrong.
\newblock {\em Basic Topology}.
\newblock UTM. Springer, first edition, 1983.

\bibitem{Bredon}
Glen~E Bredon.
\newblock {\em Topology and Geometry}.
\newblock GTM No. 139. Springer Verlag, first edition, 1993.

\bibitem{Dixmier}
Jacques Dixmier.
\newblock {\em General Topology}.
\newblock UTM. Springer Verlag, first edition, 1984.

\bibitem{Dold}
Albrecht Dold.
\newblock {\em Lectures on Algebraic Topology}.
\newblock Springer, second edition, 1980.

\bibitem{Edgar}
Gerald~A. Edgar.
\newblock {\em Measure, Topology, and Fractal Geometry}.
\newblock Undergraduate Texts in Mathematics. Springer Verlag, first edition,
  1992.

\bibitem{Fulton95}
William Fulton.
\newblock {\em Algebraic Topology, A first course}.
\newblock GTM No. 153. Springer Verlag, first edition, 1995.

\bibitem{Hilbert}
D.~Hilbert and S.~Cohn-Vossen.
\newblock {\em Geometry and the Imagination}.
\newblock Chelsea Publishing Co., 1952.

\bibitem{Kinsey}
L.~Christine Kinsey.
\newblock {\em Topology of Surfaces}.
\newblock UTM. Springer Verlag, first edition, 1993.

\bibitem{Lang97}
Serge Lang.
\newblock {\em Undergraduate Analysis}.
\newblock UTM. Springer Verlag, second edition, 1997.

\bibitem{Massey87}
William~S. Massey.
\newblock {\em Algebraic Topology: An Introduction}.
\newblock GTM No. 56. Springer Verlag, second edition, 1987.

\bibitem{Massey}
William~S. Massey.
\newblock {\em A Basic Course in Algebraic Topology}.
\newblock GTM No. 127. Springer Verlag, first edition, 1991.

\bibitem{Munkrestop}
James~R. Munkres.
\newblock {\em Topology, a First Course}.
\newblock Prentice Hall, first edition, 1975.

\bibitem{Munkresalg}
James~R. Munkres.
\newblock {\em Elements of Algebraic Topology}.
\newblock Addison-Wesley, first edition, 1984.

\bibitem{Rotman}
Joseph~J. Rotman.
\newblock {\em Introduction to Algebraic Topology}.
\newblock GTM No. 119. Springer Verlag, first edition, 1988.

\bibitem{Sato}
Hajime Sato.
\newblock {\em Algebraic Topology: An Intuitive Approach}.
\newblock MMONO No. 183. AMS, first edition, 1999.

\bibitem{Schwartz1}
Laurent Schwartz.
\newblock {\em Analyse I. Th\'eorie des Ensembles et Topologie}.
\newblock Collection Enseignement des Sciences. Hermann, 1991.

\bibitem{Seifert80}
H.~Seifert and W.~Threlfall.
\newblock {\em A Textbook of Topology}.
\newblock Academic Press, first edition, 1980.

\bibitem{Singer76}
Isadore~M. Singer and John~A. Thorpe.
\newblock {\em Lecture Notes on Elementary Topology and Geometry}.
\newblock UTM. Springer Verlag, first edition, 1976.

\bibitem{Thurston97}
Williams~P. Thurston.
\newblock {\em Three-Dimensional Geometry and Topology}.
\newblock Princeton Math. Series, No. 35. Princeton University Press, 1997.

\end{thebibliography}
\bibliographystyle{plain} 
\end{document}